\documentclass{article}
\usepackage[utf8]{inputenc}
\usepackage{amsmath}
\usepackage{amsfonts}
\usepackage{amssymb}
\usepackage{textcomp}
\usepackage[linesnumbered,ruled,vlined]{algorithm2e}
\DontPrintSemicolon 
\usepackage{multirow}
\usepackage{float}
\usepackage{array}
\usepackage[left=0.75in,right=0.75in,top=0.75in,bottom=0.75in]{geometry}
\usepackage{comment}
\usepackage{graphicx}
\usepackage{url}
\usepackage{natbib}
\usepackage{color}
\usepackage{pgfplots}
\usepackage{caption}
\usepackage{subcaption}
\usepackage{pdflscape}
\usepackage{booktabs} 
\usepackage{enumerate}
\usepackage{authblk}
\usepackage{blkarray}
\usepackage{bm}
\usepackage{marvosym} 
\usepackage{soul}
\usepackage{textcomp, xspace}
\usepackage{lipsum}
\usepackage{tikz}
\usetikzlibrary{matrix}
\usepackage{ifthen}
\usepackage{wrapfig}
\usepackage{eurosym}
\usepackage{enumitem}
\usepackage{booktabs}
\usepackage{lineno}
\usepackage{multirow}
\usepgfplotslibrary{dateplot}
\usepackage[colorlinks]{hyperref}
\usepackage{siunitx}
\usepackage{fontawesome}
\usepackage{hyperref}
\usepackage{setspace}

\hypersetup{citecolor=blue, urlcolor=blue, linkcolor=blue}
\newcommand{\fullcirc}{\textcolor{black}{\ensuremath{\bullet}}} 
\newcommand{\emptycirc}{{\ensuremath{\circ}}} 
\usetikzlibrary{patterns}
\newcommand{\Keywords}[1]{\par\noindent
{\small{\em \textbf{Keywords}\/}: #1}}

\definecolor{myblue}{RGB}{0, 0, 180}
\definecolor{addedcolor}{rgb}{0.0, 0.0, 0.0}

\usepackage{varwidth}
\definecolor{ForestGreen}{RGB}{34, 139, 34}

\title{Charge Schedule Optimization and Infrastructure Planning for Solar-Integrated Electric Bus Transit Systems}
\author[1]{Rito Brata Nath*}
\author[1]{Madhusudan Baldua*}
\author[2]{Vivek Vasudeva}
\author[2,3]{Tarun Rambha\,\textsuperscript{\faEnvelopeO}}

\affil[1]{\small Department of Civil Engineering, Indian Institute of Science (IISc), Bengaluru, India}
\affil[2]{\small Center for infrastructure, Sustainable Transportation and Urban Planning (C\textit{i}STUP), Indian Institute of Science (IISc), Bengaluru, India}  
\affil[3]{\small Robert Bosch Centre for Cyber Physical Systems (RBCCPS), Indian Institute of Science (IISc), Bengaluru, India. }
\date{}

\newcommand{\buysellenergyscenario}[2]{x_{b,#1,t}^{#2}}
\newcommand{\buysellenergysolarscenario}[2]{y_{b,#1,t}^{#2}}
\newcommand{\energytrack}[2]{u_{b,#1}^{#2}}
\newcommand{\energytracksolar}[3]{v_{#1 #2}^{#3}}
\newcommand{\gridtobattery}[3]{h_{#1 #2}^{#3}}
\newcommand{\maxenergylevelsolar}[1]{c_{#1}}
\newcommand*{\contractedcapVar}[1]{z_{#1}}
\newcommand{\areapanel}{a_{j}}
\newcommand{\initialbatteryenergy}{d_{j}}

\newcommand{\settimesteps}[1]{T^{#1}(b, k)}
\newcommand{\timeofoperation}[1]{T^{#1}}
\newcommand{\buscount}[2]{B^{#2}(#1,t)}
\newcommand{\numberofscenario}{W}

\newcommand{\locationindex}{j}

\newcommand{\scenario}{\omega}

\newcommand{\lasttimestep}{m_{j}^{\scenario}}

\newcommand{\location}{j_{b,k}}
\newcommand{\energyrequired}[1]{e_{b, k}^{#1}}
\newcommand{\solarenergy}[1]{g_{j,t}^{#1}}
\newcommand{\probability}{p^{\omega}}

\newcommand{\batterycost}{\pi}
\newcommand{\contractedcapcost}{\gamma}

\newcommand{\electricitygridprice}[2]{\epsilon_{#1}^{#2}}
\newcommand{\minenergylevel}{\rho_{min}}
\newcommand{\maxenergylevel}{\rho_{max}}

\newcommand*{\energytransferMinute}[1]{\mu_{#1}}
\newcommand{\efficiencysolarpanel}{\eta}
\newcommand{\costsolarpanel}{\alpha}

\newcommand*{\busIndex}{b}

\newcommand*{\levelVar}{\rho}

\newcommand*{\busSet}[1]{B^{#1}}

\newcommand*{\OpportunitiesSet}[1]{K^{#1}(b)}

\newcommand*{\tripSet}{I}
\newcommand*{\busrotationList}{V}
\newcommand*{\LocationSet}{J}
\newcommand*{\tripIndexI}{i}
\newcommand*{\nearestdepotTrip}[1]{\Delta_{#1}}
\newcommand*{\numTrips}{n}
\newcommand*{\tripIndexJ}{l}
\newcommand*{\bustrip}[2]{i_{#1}^{#2}}
\newcommand*{\numTripsBusb}{n_b}
\newcommand*{\compatibleSet}{A^{comp}}
\newcommand*{\busrotationAnotherList}{V^{temp}}
\newcommand*{\LocationAnotherSet}{J^{temp}}
\newcommand*{\EnergyScenario}{E^{\omega}}
\newcommand*{\depotSet}{D}
\newcommand*{\energyTrip}[1]{d^{\omega}_{#1}}
\newcommand*{\energyDeadhead}[2]{d^{\omega}_{#1,#2}}
\newcommand*{\TripStartStop}{i^{start}}
\newcommand*{\TripEndStop}{i^{end}}

\pgfplotsset{compat=1.18}
\begin{document}
\maketitle
\let\thefootnote\relax
\footnote{(*) These authors contributed equally to this manuscript. (\faEnvelopeO~tarunrambha@iisc.ac.in) Corresponding author}
\vspace{-11mm}
\begin{abstract}
As urban transit systems transition towards electrification, using renewable energy sources (RES), such as solar, is essential to make them efficient and sustainable. However, the intermittent nature of renewables poses a challenge in deciding the solar panel requirements and battery energy storage system (BESS) capacity at charging locations. To address these challenges, we propose a two-stage \textcolor{addedcolor}{multi-scenario} model that considers seasonality in solar energy generation while incorporating temperature-based variations in bus energy consumption and dynamic time-of-use electricity prices. Specifically, we formulate the problem as a multi-scenario linear program (LP) where the first-stage long-term variables determine the charging station power capacity, BESS capacity, and the solar panel area at each charging location. The second-stage scenario-specific variables prescribe the energy transferred to buses directly from the grid or the BESS during layovers. We demonstrate the effectiveness of this framework using data from Durham Transit Network (Ontario) and Action Buses (Canberra), where bus schedules and charging locations are determined using a concurrent scheduler-based heuristic. Solar energy data is collected from the National Renewable Energy Laboratory (NREL) database. We \textcolor{addedcolor}{solve} the multi-scenario LP using Benders' decomposition, which \textcolor{addedcolor}{performs} better than the dual simplex method, especially when the number of scenarios \textcolor{addedcolor}{is} high. With solar energy production at \textcolor{addedcolor}{select terminals}, our model estimated a cost savings of \textcolor{addedcolor}{$9.72\%$} and \textcolor{addedcolor}{$23.79\%$} for the Durham and Canberra networks, respectively. Our results also show that the scenario-based schedule adapts better to seasonal variations than a schedule estimated from average input parameters.

\vspace{3mm}
\end{abstract}

\section{Introduction}
\label{sec:intro}
Assessing and minimizing the carbon footprint is crucial \textcolor{addedcolor}{for combating} climate change. \textcolor{addedcolor}{Globally, in} 2019, the transport sector resulted in 8.7 gigatons of carbon dioxide equivalent (Gt{CO$_{2}$}-eq) in direct greenhouse gas emissions, an increase of 75\% 
 from its value in 1990. The transportation sector contributes to 23\% of global energy-related emissions, of which direct road transport \textcolor{addedcolor}{accounts} for almost 70\% of emissions \citep{jaramillo2022transport}. As per a report by the International Energy Agency \citep{iea_trucks_buses_emissions}, trucks and buses contributed to over 35\% of the direct CO$_{2}$ emissions from road transport, despite making up less than 8\% \textcolor{addedcolor}{of} vehicles. Owing to this growing concern, electrification of urban transit systems has gained significant impetus over the last decade, with electric buses (EBs) emerging as a critical component of sustainable public transportation solutions. \cite{iea_electric_bus} reported that over 220 heavy-duty electric vehicles (EVs) models entered the market in 2022. Governments of twenty-seven countries have committed to reaching 100\% zero-emission bus and truck sales by 2040 \textcolor{addedcolor}{\citep{GlobalMOU2021}}. 

Despite the environmental advantages of EBs, their widespread adoption is hindered by the limited energy density, short operational range, and extended recharging times. Mitigating range anxiety remains a significant challenge, especially in scenarios where EBs operate throughout the day. \textcolor{addedcolor}{The introduction of opportunity fast charging enables EBs to make more trips, promoting greater uptake.} Charging should be done so that bus battery levels do not violate the optimal depth-of-discharge limits (DoD) \citep{shaobo2021battery}. Effective managing of EB fleets requires solving two main problems: the \textit{Electric Vehicle Scheduling Problem} (EVSP) and the \textit{Charge Scheduling Problem} (CSP). \textcolor{addedcolor}{The EVSP extends the classical Vehicle Scheduling Problem (VSP) by assigning timetabled trips to electric vehicles while incorporating battery-related and operational constraints. On the other hand, the CSP focuses on determining feasible charging schedules for electric vehicles, while accounting for infrastructure, energy requirements, and operational constraints, often with a cost-minimization objective based on time-of-use (ToU) prices.}

\begin{figure}[t]
  \centering
  \begin{subfigure}[b]{0.49\linewidth}
    \centering\includegraphics[scale=0.22]{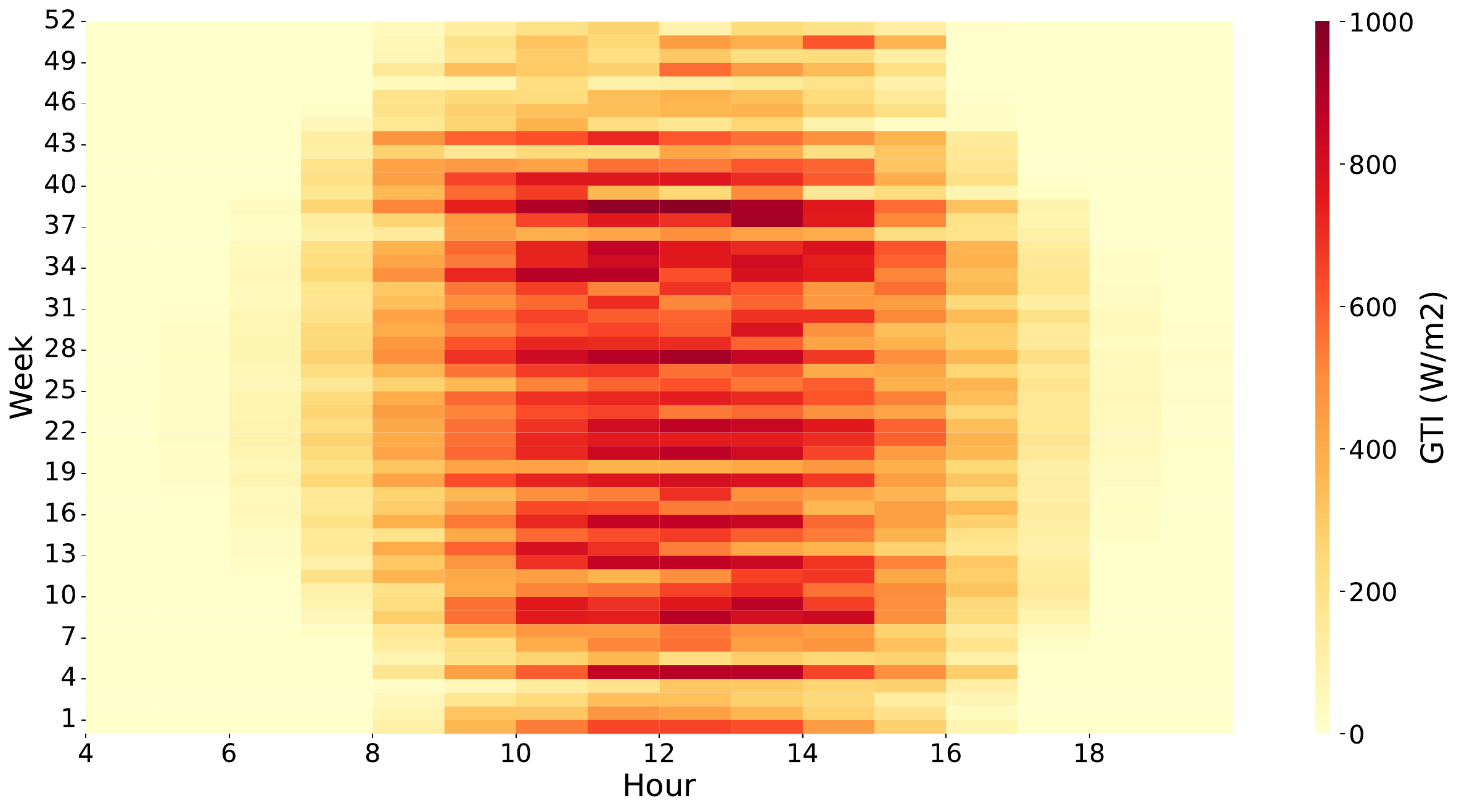}
    \caption{}
    \label{fig:gti_durham}
  \end{subfigure}
  \begin{subfigure}[b]{0.49\linewidth}
    \centering\includegraphics[scale=0.22]{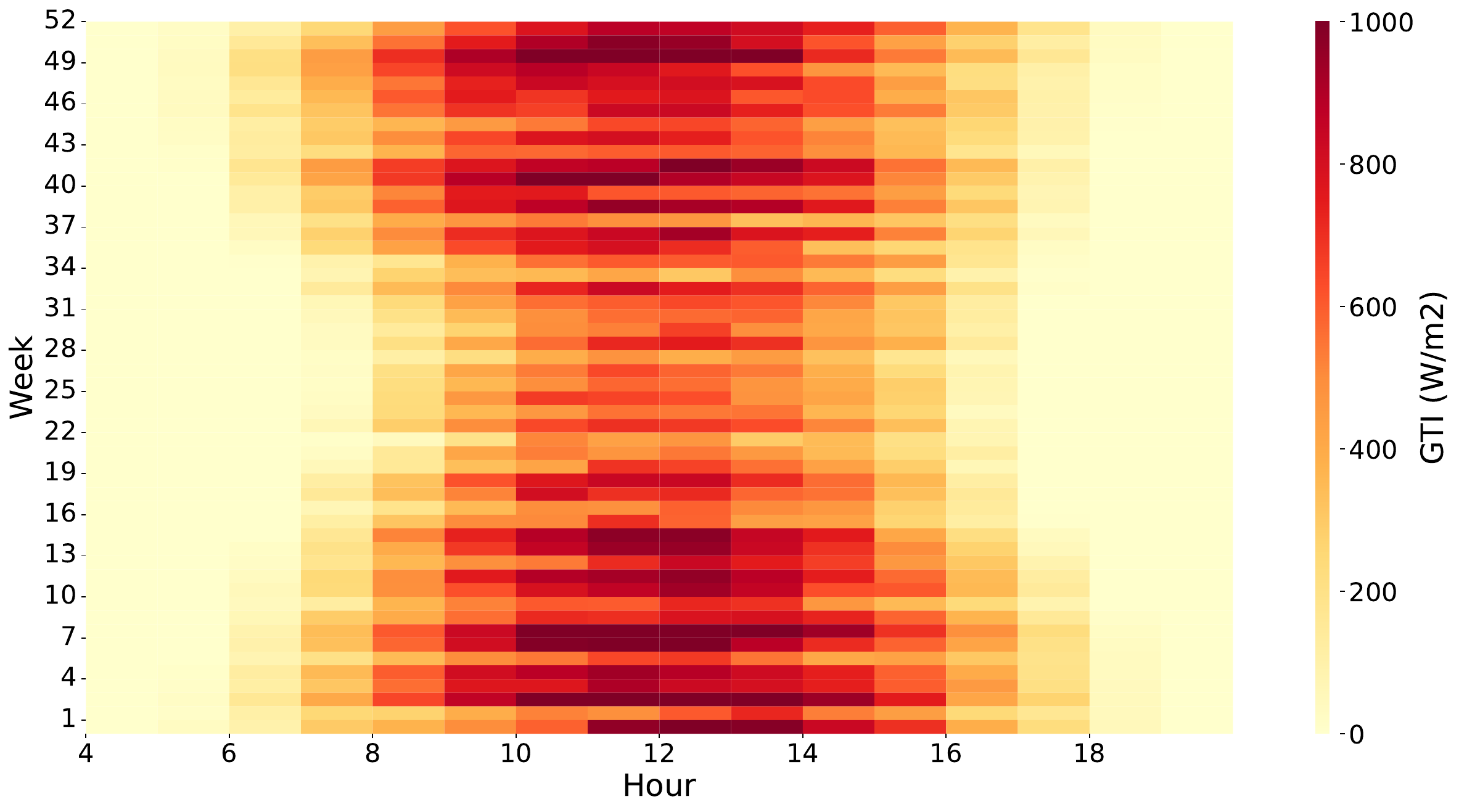}
    \caption{}
    \label{fig:gti_canberra}
  \end{subfigure}
  \caption{Weekly variations of GTI across the day for Durham (left) and Canberra (right) \textcolor{addedcolor}{in 2023 (Seasons in Durham -- Summer: Jun to Aug, Autumn: Sep to Nov, Winter: Dec to Feb, Spring: Mar to May; and seasons in Canberra -- Summer: Dec to Feb, Autumn: Mar to May, Winter: Jun to Aug, Spring: Sep to Nov)} (Source: \href{https://pvwatts.nrel.gov/}{NREL})}
\label{fig:GTI_variations}
\end{figure}

As the usage of EVs continues to grow, planners and system operators are concerned about maintaining the reliability of the power grid \citep{nrel_report}. During peak demand hours, strategic integration of charging stations with renewable energy sources (RES) such as photovoltaics (PVs) can complement the grid in an \textcolor{addedcolor}{environment-friendly} manner \citep{hassoune2017electrical, newman2020cool}. For example, many cities worldwide supplement the electricity grid with sustainable energy sources \citep{renewables_mix}. Cities such as Hobart and Canberra generate all their electricity needs from renewable energy sources. Despite this, generating additional solar energy by installing PVs locally is encouraged. This can create cheaper, self-sustaining microgrids that are resilient to cascading grid failures. However, using energy from RES presents challenges due to its stochastic nature. It is susceptible to variations in climatic conditions throughout the day, month, and year, as illustrated in Figure \ref{fig:GTI_variations} for Durham and Canberra. The plot shows the Global Tilted Irradiance (GTI), a measure of the available solar energy, across the day for different weeks of the year. One can notice that the GTI values are higher in summers between 9 am and 4 pm than at other times.

\begin{figure}[h]
    \centering
    \includegraphics[scale=0.75]{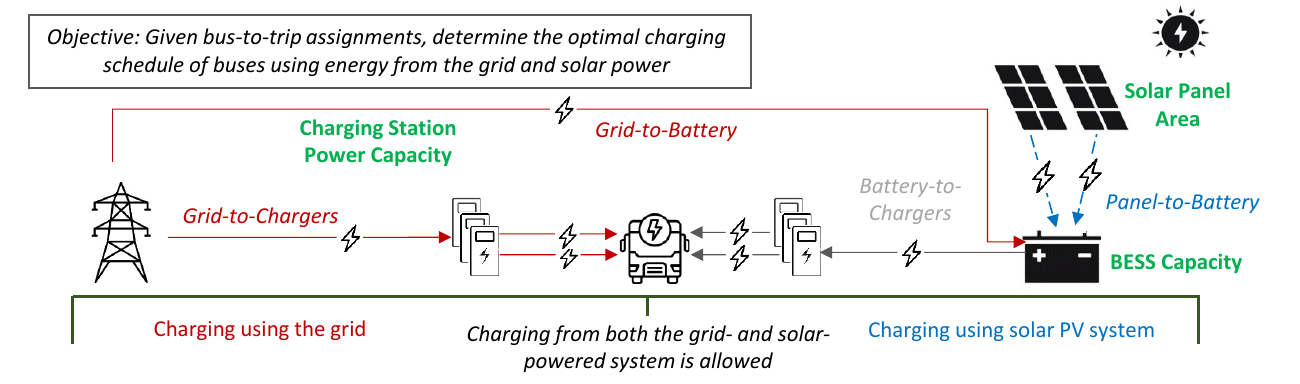}
    \caption{\textcolor{addedcolor}{Illustration of a charging depot with solar panels and battery energy storage system (BESS). Green captions denote the long-term decisions, and arrows denote the operational, short-term decisions. Red- and blue-colored labels are associated with the grid and solar-powered systems, respectively. Charging can be done using both the grid and solar PV systems}}
    \label{fig:infographic_charging_depot}
\end{figure}

From a consumption standpoint, the energy requirements of EBs are significantly influenced by driver behavior, charging strategies, road gradients, and battery life \citep{liu2021modelling, jin2022development}. Variations in ambient temperature also significantly impact trip energy consumption, implying that we can vary bus-to-trip assignments across the year for better system performance. Thus, finding solutions to the EVSP and the CSP requires addressing the planning and operational aspects of the charging infrastructure, keeping variations of supply and demand in mind. \textcolor{addedcolor}{Most existing stochastic EV charging and infrastructure planning models consider a single source of uncertainty/seasonality, such as renewable energy variability or demand fluctuations. In contrast, this study models the effect of ambient temperature on trip energy consumption and of solar irradiance on PV generation. This joint representation influences both long-term infrastructure planning and operational decisions. Ignoring temperature effects can substantially underestimate energy requirements and associated costs, while ignoring solar seasonality can underestimate grid capacity requirements.} At a planning level, we need to decide how much grid capacity to design charging depots for, what the solar panel sizes should be, and how much battery storage capacity to have. Figure \ref{fig:infographic_charging_depot} illustrates an EB depot in which buses can charge either from the grid, a solar-powered system, or both. The long-term planning decisions are marked using bold icons with green captions, and operational decisions regarding charge schedules are marked using red (grid power), blue (solar power), and gray (both grid and solar power) arrows. Depending on the realization of GTI and temperatures, in each time interval (which we refer to as a \textit{scenario}), short-term decisions on bus-to-trip assignments and charge schedules are optimized at an operational level. This scenario-based approach helps us address the following fundamental question. \textit{How beneficial is it to optimally integrate RES while operating electric transit fleets?} We address this question with a stylized, system-level model, simplifying some elements to ensure scalability. Furthermore, we demonstrate the importance of accounting for temperature-specific energy consumption in the CSP and explore solution techniques that make the framework computationally tractable for real-world networks. More specifically, our study makes the following contributions.

\begin{itemize}
    \item We propose a two-stage multi-scenario linear programming (LP) model where the first-stage variables determine the long-term decisions, such as contracted grid capacity, battery capacity to store energy, and the area of solar panels installed at each charging location. The second-stage variables are associated with a specific scenario and decide the energy the BESS draws from the grid and the energy buses take from the BESS and the grid under dynamic ToU electricity prices.
    
    \item We employ a concurrent scheduler-based (CS) heuristic that optimizes bus-to-trip assignments using the energy consumption estimates while incorporating the effect of ambient temperature. This algorithm allows us to use a different number of buses in different scenarios with a fixed set of charging locations. Charge scheduling for each scenario is solved by Benders' decomposition to improve the computational tractability for large networks. 
    
    \item We analyze the benefits of RES using a case study by constructing weekly, monthly, quarterly, and yearly scenarios for two real-world networks from Durham and Canberra with 2092 and 3911 trips, respectively. \textcolor{addedcolor}{Our model generates charge schedules with significant cost savings of $9.72\%$ and $23.79\%$ for the two networks by integrating RES at charging stations. In addition, our model shows that ignoring ambient temperature effects on energy consumption can lead to an underestimation of the total cost by $24.71\%$ and $13.34\%$ for Durham and Canberra, respectively.}

\end{itemize}

The rest of the paper is structured as follows. Section \ref{sec:litrev} presents an overview of the literature on CSP, CSP with RES, and real-world uncertainties for EBs. Section \ref{sec:preliminaries} outlines an EVSP solution, which uses a concurrent scheduler-based algorithm to generate bus-to-trip assignments.
Section \ref{sec:scenario_csp} provides a CSP formulation with seasonality in solar energy production and ambient temperature-based trip energy consumption. We also describe Benders' decomposition, which is used to solve our formulation.
In Section \ref{sec:results}, we present the results of our case studies on the Durham and Canberra transit networks and highlight valuable insights. We summarize our findings in Section \ref{sec:conc} and outline potential areas for future research. We list the abbreviations used in this paper in \ref{sec:abbreviations}. 

\section{Literature Review}
\label{sec:litrev}

Two main problems in managing EB fleet operations \textcolor{addedcolor}{are} the EVSP and the CSP. The EVSP provides a bus-to-trip assignment, whereas the CSP determines the charging schedules for EBs based on ToU electricity prices. The optimized charge schedules can also be used to decide where to open charging stations from a set of candidates. We discuss relevant literature on CSP in Section \ref{sec:csp_lit} and a summary of CSP with renewables and various sources of uncertainty in Section \ref{sec:renewables_lit}. An overview of EB fleet management problems can be found in \cite{perumal2022electric}.

\subsection{Charge Scheduling Problem}
\label{sec:csp_lit}
The charging schedule of EBs affects the operational cost, primarily due to the variability in electricity pricing and energy demand \citep{yang2017charging}. The primary goal of the CSP is to minimize the operational cost with known bus-to-trip assignments and potential charging locations. It can account for ToU electricity prices, grid capacity, charging infrastructure capacity, and integration of RES. Apart from ToU energy tariffs, many studies have examined the impact of the EV charging schedules on power grid systems \citep{alonso2014optimal, deb2018impact}.
Uncoordinated charging of energy-intensive EBs poses a significant challenge to the stability of the grid \citep{korolko2015robust,jian2017high, zhou2022robust}. Therefore, effectively managing grid constraints while meeting state-of-charge (SoC) requirements becomes crucial \citep{qin2016numerical,leou2017optimal,  he2020optimal}. 

Several studies have focused on reducing the peak load on the grid. \cite{dietmannsberger2017modelling} investigated the impacts of replacing diesel buses with EBs on the power grid system for the Hamburg bus depot. \cite{jahic2019charging} implemented greedy heuristics to minimize peak load under centralized depot charging. \cite{houbbadi2019optimal} used non-linear programming-based methods for optimizing overnight charging decisions with battery aging costs for EB fleets. \cite{abdelwahed2020evaluating} formulated discrete time-based and event-based mixed-integer linear
programs (MILPs) for optimizing the opportunity fast-charging schedule for a single depot. \cite{zhou2022robust} and \cite{he2023joint} imposed demand charges to compute precise charging costs while formulating mixed-integer second-order cone programming and mixed-integer non-linear
programming problems, respectively. Demand charges indicate the electricity billing based on the total energy requirement from the grid within 15 minutes or one hour. \cite{wang2024optimal} modeled each charger separately in their charging infrastructure and scheduling problem and proposed a bi-level optimization framework. \cite{nath2024impact} integrated location and scheduling problems by proposing split and uniform charging models based on charging priority rules while considering dynamic ToU pricing and charging station power capacity. In our study, we extend some of their CSP models by integrating solar energy sources.  

\subsection{Renewables and Uncertainties}
\label{sec:renewables_lit}
Economic and environmental benefits of EBs hinge on the effective integration of RES for either complete or partial vehicle charging \citep{yao2017optimal}. A review of models for smart charging of EVs using PVs and ToU pricing can be found in \cite{fachrizal2020smart}. Energy storage systems are a potential solution to reduce high electricity demand charges at charging locations \citep{he2019fast, zhong2024joint}. Recent studies have focused on such integrated systems. \cite{leou2017optimal} integrated wind energy sources into the power grid, and \cite{jin2022development} focused on a PV-based charging station with BESS to optimize the overall operational cost for bus and parking agencies, respectively. \cite{liao2021decentralized} explored the cost potential of EBs by considering their bus-to-grid capabilities to sell electricity back to the grid. \cite{ke2020battery} utilized a genetic algorithm to optimize charging and discharging schedules for the Penghu bus transportation system, considering factors such as solar and wind power generation, feeder load, and demand response. \cite{yusuf2023comparative} addressed uncertainties related to PV systems and building load patterns by employing 15-minute long short-term memory predictions. Their objective was to minimize the overall cost of workplace-integrated microgrids. \cite{singh2020implementation} explored a system with PVs, BESS, a diesel generator, and grid-connected charging stations to ensure constant and uninterrupted charging, whereas \cite{li2024optimizing} proposed a system in which EBs are equipped with rooftop solar panels to provide on-board energy supply. To address the uncertainty in solar energy generation, \cite{li2020energy} used chance-constrained LPs and \cite{kumar2022multiobjective} proposed a two-stage multi-objective planning framework to coordinate the scheduling of EVs, BESS, and PVs at a charging station. \cite{liu2023optimal} addressed the charging location planning problem, and \cite{liu2024electric} addressed the charge scheduling problem integrated with PV and BESS using an MILP model. \textcolor{addedcolor}{\cite{hu2025optimal} proposed a charge scheduling problem integrated with station-level photovoltaic requirements using an MILP model, which was solved using Lagrangian relaxation and a dynamic programming algorithm. \cite{zhang2025optimal} proposed a non-linear model to address the charge scheduling integrated with PV planning, and later converted it into a linear tractable model. However, they did not consider seasonality in solar energy generation.}

\cite{zhang2013charging} modeled the uncertainty in the EV arrival rates as a Markov process with variable grid power prices at different periods. Similarly, to address the energy consumption variability during bus operations, \cite{zhuang2020stochastic} developed stochastic energy management strategies at charging stations. Their approach involved a distributionally robust Markov Decision Process (MDP) with RES. \cite{huang2022two} formulated a two-stage stochastic program with a hybrid charging scheme, which was solved using a progressive hedging algorithm. \textcolor{addedcolor}{\cite{hendriks2024integrated} formulated an MILP model to minimize the grid load and number of chargers, considering uncertainties in trip energy consumption.} \cite{rafique2022two} proposed a stochastic two-stage multi-objective model for charging and discharging schedules of EBs to minimize the depot operator's energy and battery degradation costs. To enhance tractability, \cite{soares2017two} used Benders' decomposition on a two-stage stochastic model for an intelligent grid, which addresses variations in demand response, RES, EVs, and market prices. Likewise, \cite{najafi2025integrated} used Benders' decomposition to solve the charging infrastructure and scheduling problem using renewables. However, they simplified their model by combining trips based on charging time-window availability and considered only three scenarios for PV generation. 

Due to the sporadic nature of RES, determining the appropriate battery storage capacity and the solar panel capacity at a station poses a significant challenge. \textcolor{addedcolor}{Variations in renewable generation and trip energy demand strongly influence daily charging schedules.} A practical approach to dealing with these variations involves considering a range of possible scenarios and their associated probabilities of occurrence. Findings from \cite{arif2020plug} indicate that deploying PVs and BESS at charging stations using a synchronized charging schedule is a practical choice for handling peak loads. We address a similar CSP problem that can deal with variations in the energy consumption of the EBs and RES availability.

\begin{table}[t]
    \centering
    \small
    \caption{Model features in relevant works on EB charge scheduling (1: Grid Capacity, 2: Solar Infrastructure, 3: Partial Charging, 4: BESS, 5: RES, 6: Uncertainty/Seasonality, 7: Model/Solution Technique, 8: Problem Scale (we report only the number of EBs when the trip information is not available), $\emptycirc$ not clear in the paper, blank spaces represent features that are not modeled)}
    \begin{tabular}{p{35mm}p{1.7mm}p{1.7mm}p{1.7mm}p{1.7mm}p{1.7mm}p{30mm}p{37mm}p{21mm}}
    \hline
      \textbf{Ref.} & \centering \textbf{(1)} & \centering \textbf{(2)} & \centering \textbf{(3)} & \centering \textbf{(4)} & \centering \textbf{(5)} & \textbf{(6)} & \textbf{(7)} & \textbf{(8)}\\
      \hline\cite{abdelwahed2020evaluating} & \centering\fullcirc & & \centering\fullcirc  & & & & MILP & 47 EBs\\
        \cite{arif2020plug} & \centering \fullcirc & \centering \emptycirc & & \centering \fullcirc & \centering \fullcirc & & MILP and charge scheduling algorithm & 144 trips, 24 EBs\\
        \cite{hendriks2024integrated} & \centering \fullcirc & \centering  & \centering \fullcirc & \centering & \centering & \textcolor{addedcolor}{Energy consumption} & \textcolor{addedcolor}{MILP} & \textcolor{addedcolor}{18 routes, 85 EBs}\\ 
        \cite{huang2022two} & \centering \emptycirc & & & & &Energy consumption & Progressive hedging algorithm & 6 EBs \\
         \cite{ke2020battery} & \centering \fullcirc & \centering \fullcirc & \centering \fullcirc &\centering \emptycirc &\centering \fullcirc & & Genetic algorithm & 29 EBs\\
          \cite{leou2017optimal} & \centering \fullcirc & & \centering \fullcirc & & \centering \fullcirc & Charging energy & Monte Carlo simulation & 17 trips, 10 EBs \\
          \cite{liu2024electric} & \centering \fullcirc & \centering  \emptycirc & \centering \fullcirc & \centering \fullcirc & \centering \fullcirc & \textcolor{addedcolor}{Solar energy} & \textcolor{addedcolor}{MILP} & \textcolor{addedcolor}{8 bus lines, 50 EBs}\\ 
        \cite{manzolli2022electric} & \centering \fullcirc & & \centering \emptycirc & & & & MILP & 11 EBs\\
        \cite{najafi2025integrated} & \centering \fullcirc & \centering \fullcirc & \centering \emptycirc & \centering \fullcirc & \centering \fullcirc & Solar energy & Benders' decomposition & 316 trips, 10 EBs\\
         \cite{rafique2022two} & \centering \fullcirc & \centering \fullcirc & & \centering \fullcirc & \centering \fullcirc & EB schedule & MILP & 50 EBs\\
          \cite{soares2017two} & \centering \fullcirc & \centering \emptycirc & & \centering \fullcirc & \centering \fullcirc & Power system  & Benders' decomposition  & 180 EBs\\
          \cite{zhuang2020stochastic} & \centering \fullcirc & \centering \emptycirc & & \centering \fullcirc & \centering \fullcirc & Energy consumption & Distributionally robust MDP and heuristics & Not clear \\ 
          \hline
          \textbf{Our Work} & \centering \fullcirc & \centering \fullcirc & \centering \fullcirc & \centering \fullcirc & \centering \fullcirc &  Energy consumption and solar energy  & Benders' decomposition & 3911 trips, 265 EBs \\
          \hline
  
    \end{tabular}
    \label{tab:litreview}
\end{table}

Table \ref{tab:litreview} compares features of our work with existing literature. Although studies related to RES in the context of EVs are abundant, literature on EBs integrated with RES is scarce. The schedules of EVs requiring charging at a depot are challenging to predict. Hence, most studies that integrate EV operations with RES assume random arrivals of EVs at a charging station \citep{zhang2013charging, yao2017optimal}. For EBs, the timetable is fixed by transit operators, so it is possible to plan strategically and determine the exact itinerary for each bus. However, the intermittent nature of RES remains a challenge in integrated systems. To our knowledge, no study has tackled the CSP with seasonal solar energy generation and variable bus-to-trip assignments for EBs. Our study bridges this gap by proposing a multi-scenario LP formulation, which is solved using Benders' decomposition \citep{geoffrion1972generalized}. Readers can refer to \cite{rahmaniani2017benders} for a detailed review of Benders' decomposition. Most studies in Table \ref{tab:litreview} do not model solar infrastructure and battery storage capacity decisions. In addition to capturing these features, we consider partial bus charging through a split charging model using grid- and solar-powered systems. Another unique aspect of our study is the scale of the problems addressed. The largest network used in our research includes 265 EBs, serving 3911 daily trips. Furthermore, a high number of scenarios makes our problem more practical. Before delving into our problem formulation, we present some background on EVSP in Section \ref{sec:preliminaries}.

\section{Preliminaries}
\label{sec:preliminaries}

Transit operations are characterized by routes/lines and trips carried out on each route. \textcolor{addedcolor}{We assume a representative timetable for trips that repeats daily within a given planning period, consistent with standard practice in vehicle scheduling. The framework can be extended to incorporate multiple timetable types, such as weekends or different seasonal timetables, either by solving separate instances or by embedding them in a multi-period model. While absolute costs may differ across schedules, the proposed framework remains applicable and can be used to evaluate the impacts of RES integration and temperature effects under alternative timetable structures. To maintain tractability, we focus on a single schedule (e.g., weekday operations) in this study.} 

When the fleet contains EBs, operators must also identify charging locations and charging schedules and assign buses to trips. Two typical charging methods \textcolor{addedcolor}{are} overnight depot charging and within-day opportunity charging. We assume that overnight depot charging replenishes bus batteries to the same energy level daily, keeping charging schedules periodic and simplifying managerial decisions. \textcolor{addedcolor}{Opportunity fast charging is done during the day to meet the energy needs of buses and to ensure that the battery levels of buses stay higher than the minimum threshold \citep{hu2022joint}.} Terminals where many trips originate or end are chosen as depots for overnight charging. Chargers are made available at such locations, and we assume that buses charge overnight at the depots after their daily operations. Depot chargers can also be used for opportunity charging during the day. In addition, non-depot terminals serve as candidate locations for installing opportunity chargers. Buses are assumed to start their journey at a depot, serve one or more trips during the day, and return to the same depot at the end of the day. We assign buses to trips based on trip compatibility and energy requirements. Two trips $i_1$ and $i_2$ are compatible if a bus can serve trip $i_2$ after serving trip $i_1$. This is possible if the end time of trip $i_1$ plus the time to deadhead from the end-stop of $i_1$ to the start-stop of $i_2$ is less than or equal to the start time of trip $i_2$. \textcolor{addedcolor}{We define the set of compatible trips as: $
A^{\text{Comp}} = \left\{ (i_1,i_2) \;\middle|\; t^{end}_{i_1} + t_{i_1,i_2} \le t^{start}_{i_2}, i_1 \in I, i_2 \in I, i_1 \neq i_2 \right\}$, where \( t^{end}_{i_1} \) denotes the end time of trip \( i_1 \), \( t_{i_1,i_2} \) is the deadhead time between trips \( i_1 \) and \( i_2 \), \( t^{start}_{i_2} \) denotes the start time of trip \( i_2 \), and $I$ is the set of all trips.}

\begin{table}[t]
 \small
  \centering
  \caption{ Notation used in the CS algorithm}
    \begin{tabular}{lp{130mm}}
    \hline
    \textbf{Notation} & \textbf{Description} \\
    \hline
    \textbf{Index:} \\ \hline
    $b$ & Bus index \\
    $\scenario$ & Index of a scenario \\
    \hline \textbf{Sets/Lists:} \\
    \hline
    $\tripSet$ & Set of trips (an individual trip is indexed by $i$) \\
    $\depotSet$ & Set of depots \\
    $\numberofscenario$ & Set of scenarios \\
    $\busrotationList^{\scenario}$ & List of rotations of all the EBs for scenario $\scenario$ \\
     $\busrotationList$ & List of rotations of all the EBs for all scenarios, i.e., $(V_{\scenario})_{\scenario \in \numberofscenario}$\\
     $\busrotationList_b^{\scenario}$ & Rotation of $\busIndex^\text{th}$ electric bus in scenario $\scenario$ \\
     $\LocationSet$ & Set of charging locations for all scenarios, including the depots \\
     $\LocationSet^{\scenario}$ & Set of charging locations for scenario $\scenario$ including the depots \\
     $\compatibleSet$ & Set of compatible trip pairs \\
     $\EnergyScenario$ & List of energy requirements for service and deadheading trips for scenario $\scenario$, i.e., $[(\energyTrip{\tripIndexI})_{i \in I}, \, (\energyDeadhead{\tripIndexI}{\tripIndexJ})_{\tripIndexI \in \tripSet, \tripIndexJ \in \depotSet}, (\energyDeadhead{\tripIndexI}{\tripIndexJ})_{\tripIndexI \in \depotSet, \tripIndexJ \in \tripSet}, (\energyDeadhead{\tripIndexI}{\tripIndexJ})_{\tripIndexI \in \tripSet, \tripIndexJ \in \tripSet}]$ \\
     \textcolor{addedcolor}{$V_b^{temp}$} & \textcolor{addedcolor}{Temporary rotation list of $b^{\text{th}}$ electric bus} \\
    \textcolor{addedcolor}{$J^{temp}$} & \textcolor{addedcolor}{Set for temporarily storing the charging locations, including the depots} \\
    \hline
    \textbf{Data:} \\
    \hline
    \\[-3mm]
     $\bustrip{\busIndex}{\tripIndexJ}$ & $\tripIndexJ^\text{th}$ trip performed by bus $\busIndex$ \\
     $\nearestdepotTrip{\tripIndexI}$ & Nearest depot location from the starting bus stop location of trip $\tripIndexI \in \tripSet$ \\
     $\numTrips$ & Number of trips in the timetable \\
     $\numTripsBusb$ & Number of trips performed by bus $b$\\
     $\energyTrip{\tripIndexI}$ & Energy consumed (in kWh) on trip $\tripIndexI$ in scenario $\scenario$\\
     $\energyDeadhead{\tripIndexI}{\tripIndexJ}$ & Energy (in kWh) for deadheading from the end-stop of $\tripIndexI$ to the start-stop of $\tripIndexJ$ in scenario $\scenario$ \\
     $\TripStartStop$ & Starting stop of trip $\tripIndexI$\\     $\TripEndStop$ & Ending stop of trip $\tripIndexI$\\
     $\levelVar$ & Current energy level (in kWh) of a bus \\
     $\delta^{\omega}_{\bustrip{\busIndex}{\tripIndexJ}, \bustrip{\busIndex}{\tripIndexJ+1}}$ & Time available for charging between the end of $\tripIndexJ^\text{th}$ trip and start of $(\tripIndexJ+1)^\text{th}$ trip in bus $b$'s rotation in scenario $\omega$\\
     \textcolor{addedcolor}{$r_i^\omega$} & \textcolor{addedcolor}{Energy consumption rate (in kWh/min) of trip $i$ under scenario $\omega$} \\
    \textcolor{addedcolor}{$i_b^{l,end}$} & \textcolor{addedcolor}{End stop of $l^{\text{th}}$ trip of bus $b$} \\[-3mm]\\
    \hline 
    \textbf{Parameters:} \\
    \hline
    $\maxenergylevel$ & Maximum battery capacity (kWh) of the bus up to which it can be charged \\
  $\minenergylevel$ & Minimum battery capacity (kWh) of the bus below which it cannot be discharged\\
  $\energytransferMinute{}$ & Maximum amount of energy (kWh/min) that can be transferred to the bus either from the power grid or from the BESS in any time step\\
\hline
    \end{tabular}%
  \label{tab:notations_cs}%
\end{table}%

A \textit{scenario} or a time interval can be based on different levels of aggregation such as weeks, months, quarters, or a year and is characterized by two properties -- GTI values and ambient temperatures across different times of the day. The GTI values are used to decide the amount of energy generated from solar power in each scenario. Ambient temperature variations, on the other hand, influence the energy consumption for service and deadheading trips. Due to the differences in energy consumption for the same trip across scenarios, we also employ different EVSP assignments for different scenarios. The literature on the electric bus-to-trip assignment problem or the EVSP is abundant. EVSP is an extension of the MDVSP that can be formulated either as a multi-commodity flow model or a set partitioning model \citep{ribeiro1994column, hadjar2006branch} and is an NP-hard problem \citep{bertossi1987some}. Researchers have tackled the EVSP using MILP models and column generation for small-sized networks \citep{adler2017vehicle}. Heuristics such as genetic algorithms \citep{yao2020optimization}, simulated annealing \citep{zhou2020collaborative}, and adaptive large neighborhood search \citep{wen2016adaptive} have also been used to solve large instances. 

\begin{algorithm}[t]
\caption{\textsc{ConcurrentScheduler(CS)}}
\label{alg:cs}
\KwIn{$\tripSet$, $\depotSet$, $\numberofscenario$, $\EnergyScenario\, \forall \, \scenario \in \numberofscenario$}
\KwOut{$\LocationSet, \busrotationList^{\scenario} \, \forall \, \scenario \in W$}

\For{$\scenario \in \numberofscenario$}{
$\busrotationList^{\scenario}, \LocationSet^{\scenario} \gets \emptyset, \depotSet$ \;
Rearrange $\tripSet$ in the ascending order of the trip start times\;
Pick the first trip $\tripIndexI$ from the sorted list of $\tripSet$\;
$\busrotationList_1^{\scenario} \gets [\nearestdepotTrip{\tripIndexI}, \tripIndexI, \nearestdepotTrip{\tripIndexI}]$ and add $\busrotationList_1^{\scenario}$ to $\busrotationList^{\scenario}$\;

{\color{blue}{\tcp{\textsf{Add trips to existing rotations}}}}
\For{$\tripIndexJ = 2$ \KwTo $\numTrips$}{
    $\tripIndexI_{\tripIndexJ} \gets$ The $\tripIndexJ^\text{th}$ element from the sorted list of $\tripSet$\;
    
    \For{$b = 1$ \KwTo $|\busrotationList^{\scenario}|:$ $(\bustrip{b}{\numTripsBusb}, \tripIndexI_\tripIndexJ) \in \compatibleSet$}
        {
            $\busrotationAnotherList_b \leftarrow \busrotationList_b^{\scenario}$ and insert $\tripIndexI_{\tripIndexJ}$ to the list of service trips in $\busrotationAnotherList_b$ \;
            $\textsf{insertTrip}, \LocationAnotherSet \gets \textsc{IsRotationChargeFeasible}(\busrotationAnotherList_b, \LocationSet^{\scenario}, \EnergyScenario)$\;
            \If{$\mathsf{insertTrip}$}{
                $\busrotationList_b^{\scenario} \gets \busrotationAnotherList_b$ and $\LocationSet^{\scenario} \gets \LocationAnotherSet$\;
                \textbf{break}
            }
        }
    {\color{blue}{\tcp{\textsf{Create new rotations}}}}
    \If{$\tripIndexI_\tripIndexJ$ is not assigned to any of the buses used so far}{
        Use a new bus $|\busrotationList^{\scenario}| + 1$, update its rotation to $[\nearestdepotTrip{\tripIndexI_\tripIndexJ}, \tripIndexI_\tripIndexJ, \nearestdepotTrip{\tripIndexI_\tripIndexJ}]$, and add it to $\busrotationList^{\scenario}$\;
    }   
}    
}
$\LocationSet \gets \bigcup_{\scenario \in \numberofscenario} \LocationSet^{\scenario}$ \;

\end{algorithm}

Our scenario-specific bus schedule or \textit{rotations} and the initial set of charging locations are generated using a CS-based heuristic \citep{bodin1978ucost, nath2024impact}, which only focuses on the feasibility of bus charging and does not specifically integrate solar-powered systems. \textcolor{addedcolor}{The bus-to-trip assignment is generated exogenously using the CS heuristic and is not jointly optimized with charging decisions. While a fully integrated formulation could improve coordination, it results in a large-scale mixed-integer stochastic program that is computationally intractable for realistic instances. We therefore adopt a sequential approach that prioritizes tractability.} We present a modified version of the approach discussed in \citep{nath2024impact} in Algorithm \ref{alg:cs}. Refer to Table \ref{tab:notations_cs} for notation. For a given scenario $\scenario$, the algorithm starts by initializing the set of charging locations to the set of depots (line 2). The service trips are then sorted in ascending order of their start times (line 3). A new bus rotation is created by assigning the first trip to the bus, which we assume is always charge feasible (line 5). The algorithm prioritizes inserting trips from the sorted list of service trips to a bus rotation in use (lines 7--13) based on trip compatibility (line 8) and charging level feasibility (line 10). The feasibility of the bus energy levels is verified using the \textsc{IsRotationChargeFeasible} function, which takes \textcolor{addedcolor}{as inputs} the current bus rotation, charging locations, and the scenario-specific energy consumption. 

We present a pseudocode for this function in Algorithm \ref{alg:feasibility_check}. Charge feasibility checks and modifications to the charging locations are performed based on a \textit{Charge-and-Go} (CAG) strategy \citep{nath2024impact} where we prioritize charging a bus at the end-stop of a trip instead of the start-stop of the subsequent trip. \textcolor{addedcolor}{While this approach restricts charging flexibility and may lead to sub-optimal solutions, it helps us formulate a tractable CSP. Allowing charging decisions to be made freely at both trip ends would require additional binary decision variables that result in a complex mixed-integer linear program. Note that this strategy only affects cases where charging infrastructure is available at both the end-stop of a trip and the start-stop of the next trip.} If the trip insertion (lines 11--13) is not feasible, we create a new bus rotation (lines 14--15) and assign the trip to it. The scenario-specific bus rotations and initial charging locations are inputs to our CSP model, which will be explained in Section \ref{sec:scenario_csp}. \textcolor{addedcolor}{To ensure that our results are not driven by this modeling choice for EVSP, we also tested an alternative initial solution obtained via a local search-enhanced EVSP approach, and observed that while absolute costs vary slightly, the key insights remain consistent (see Section \ref{sec:local_search}).} \textcolor{addedcolor}{In experiments on small instances, an integrated formulation produced a prohibitively large number of variables and constraints, resulting in memory issues. This motivated us to use a scalable heuristic for schedule generation.} The union of all scenario-specific charging locations is treated as the set of initial charging locations for the CSP (line 16).

\begin{figure}[H]
    \centering
    \includegraphics[width=0.9\textwidth]{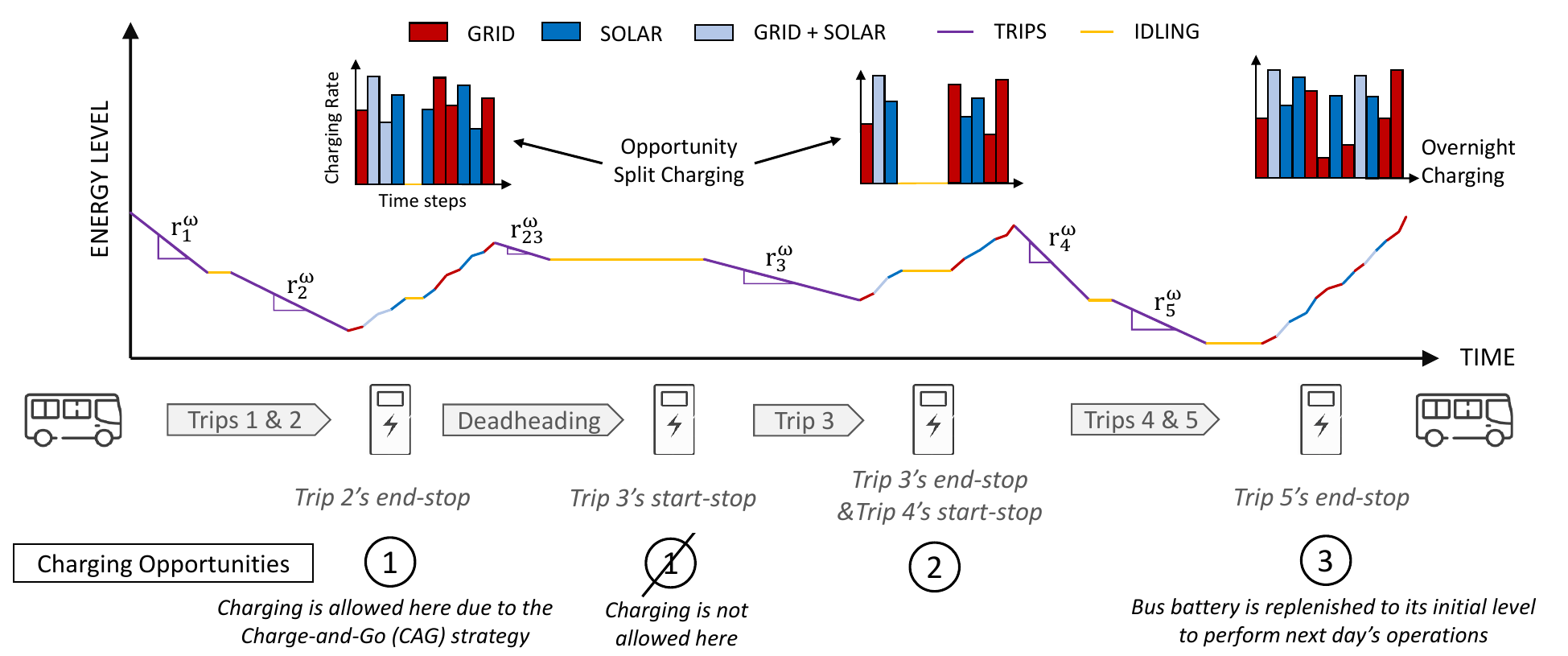}
    \caption{Charging operations for a bus showing energy levels during trips (purple), charging through grid (red), solar (dark blue) or both (light blue), and idling (yellow)}
    \label{fig:busoperation}
\end{figure}

Figure $\ref{fig:busoperation}$  illustrates the operation of a bus with five service trips, a deadheading trip between two service trips, and charging opportunities $1$, $2$, and $3$. A \textit{charging opportunity} arises if a bus has a layover at a charging station (which could also be a depot) either at the end of a trip or at the beginning of a trip. During this layover, the bus may charge, deadhead, or remain idle. Charging opportunities typically occur at a terminal equipped with charging stations or at the depot during overnight hours.

The figure illustrates that charging opportunities 1 and 2 occur between service trips, whereas overnight charging occurs during charging opportunity 3. The energy levels for the service trips are marked in purple; charging activities are shown in red (grid-powered), dark blue (solar-powered), and light blue (both grid- and solar-powered); and idling after trips or during charging opportunities is depicted in yellow. Note that solar-powered charging indicates that the bus is solely charged using solar energy, and the BESS has no residual energy left from the grid. 

The figure also illustrates the CAG charging strategy for a bus where the operator encounters a dilemma to either charge at the end-stop of trip 2 or the start-stop of trip 3 when charging stations are available at both places. The CAG strategy assumes that charging opportunity $1$ is assigned to trip 2's end-stop \textcolor{addedcolor}{and not to} trip 3's start-stop. \textcolor{addedcolor}{Note that the CAG strategy intentionally restricts spatial flexibility by prioritizing charging at trip end-stops, thereby avoiding additional charging location selection decisions.} \textcolor{addedcolor}{We make this choice to help retain a high-resolution (one-minute) time discretization in the CSP to accurately capture charge opportunities during layovers in schedules and ToU electricity price variations.}

The bus charges at charging opportunity 1 and deadheads to the start-stop of the next trip, arriving exactly at its departure time. This assumption helps eliminate additional decisions on when buses should deadhead to the next trip's start location and allows casting the CSP as an LP instead of an MILP, thus improving its tractability. We assume that charging rates may vary over time, as shown in the plots at the upper level of the figure. In other words, the energy supplied during a charging opportunity may be split non-uniformly. The energy consumption rate for service and deadheading trips (see $r_1^{\scenario}, r_{2}^{\scenario}, r_{23}^{\scenario}, \ldots, r_5^{\scenario}$) can vary depending on the average ambient temperature during the trip, its length, and travel time. Note that the slopes are related to $\energyTrip{\tripIndexI}$ and $\energyDeadhead{\tripIndexI}{\tripIndexJ}$ mentioned in Table \ref{tab:notations_cs}. In practice, the energy consumption rates for a service/deadheading trip need not be constant. We simplify it in our model since these rates have no bearing on the CSP decisions. 

\begin{algorithm}[h]
\caption{\textsc{IsRotationChargeFeasible}}\label{alg:feasibility_check}
\KwIn{$\busrotationAnotherList_b, \textcolor{addedcolor}{\LocationSet^{\scenario}}, \EnergyScenario$}
\KwOut{True/False, \textcolor{addedcolor}{$\LocationSet^{temp}$}}
    $\levelVar \gets \maxenergylevel - \energyDeadhead{\bustrip{\busIndex}{0}}{\bustrip{\busIndex}{1}} - \energyTrip{\bustrip{\busIndex}{1}}$, $\textcolor{addedcolor}{\LocationSet^{temp} \gets} \textcolor{addedcolor}{\LocationSet^{\scenario}}$\;
    \lIf{$\levelVar < \minenergylevel$}{
        \Return \textbf{false}, \textcolor{addedcolor}{$\LocationSet^{temp}$}}
    \For{$\tripIndexJ = 1$ \KwTo $(\numTripsBusb - 1)$}{
{\color{blue}{\tcp{\textsf{Current trip end is a charging station}}}}
        \If{$\bustrip{\busIndex}{\tripIndexJ, end} \in \textcolor{addedcolor}{\LocationSet^{temp}}$}{
        $\levelVar \gets \min\{\maxenergylevel, \levelVar + \energytransferMinute{} \, \delta^{\omega}_{\bustrip{\busIndex}{\tripIndexJ}, \bustrip{\busIndex}{\tripIndexJ+1}}\} - \energyDeadhead{\bustrip{\busIndex}{\tripIndexJ}}{\bustrip{\busIndex}{\tripIndexJ+1}}$\;
            }
        {\color{blue}{\tcp{\textsf{Only next trip start is a charging station}}}} 
        
        \If{$\bustrip{\busIndex}{\tripIndexJ, end} \notin \textcolor{addedcolor}{\LocationSet^{temp}}$, $\bustrip{\busIndex}{\tripIndexJ+1, start} \in \textcolor{addedcolor}{\LocationSet^{temp}}$}{
        $\levelVar \gets \min\{\maxenergylevel, \levelVar - \energyDeadhead{\bustrip{\busIndex}{\tripIndexJ}}{\bustrip{\busIndex}{\tripIndexJ+1}} + \energytransferMinute{} \,\delta^{\omega}_{\bustrip{\busIndex}{\tripIndexJ}, \bustrip{\busIndex}{\tripIndexJ+1}}\}$\; 
        }
        {\color{blue}{\tcp{\textsf{Current trip end and next trip start are not charging stations}}}}

        \If{$\bustrip{\busIndex}{\tripIndexJ, end} \notin \textcolor{addedcolor}{\LocationSet^{temp}}, \bustrip{\busIndex}{\tripIndexJ+1, start} \notin \textcolor{addedcolor}{\LocationSet^{temp}}, \levelVar  - \energyDeadhead{\bustrip{\busIndex}{\tripIndexJ}}{\bustrip{\busIndex}{\tripIndexJ+1}} - \energyTrip{\bustrip{\busIndex}{\tripIndexJ+1}} \geq \minenergylevel$}{
        $\levelVar \gets \levelVar  - \energyDeadhead{\bustrip{\busIndex}{\tripIndexJ}}{\bustrip{\busIndex}{\tripIndexJ+1}}$\;}
        \If{$\bustrip{\busIndex}{\tripIndexJ, end} \notin \textcolor{addedcolor}{\LocationSet^{temp}}, \bustrip{\busIndex}{\tripIndexJ+1, start} \notin \textcolor{addedcolor}{\LocationSet^{temp}}, \levelVar  - \energyDeadhead{\bustrip{\busIndex}{\tripIndexJ}}{\bustrip{\busIndex}{\tripIndexJ+1}} - \energyTrip{\bustrip{\busIndex}{\tripIndexJ+1}} < \minenergylevel$}{
        {\color{blue}{\tcp{\textsf{Open a new charging station at the current trip end to meet the charging requirement}}}}   
        $\levelVar \gets \min\{\maxenergylevel, \levelVar + \energytransferMinute{} \,\delta^{\omega}_{\bustrip{\busIndex}{\tripIndexJ}, \bustrip{\busIndex}{\tripIndexJ+1}}\} - \energyDeadhead{\bustrip{\busIndex}{\tripIndexJ}}{\bustrip{\busIndex}{\tripIndexJ+1}}$\;
        {\color{blue}{\tcp{\textsf{Update the charging locations if the next trip can be performed after charging}}}}
        \lIf{$\levelVar - \energyTrip{\bustrip{\busIndex}{\tripIndexJ+1}} \geq \minenergylevel$}{$\textcolor{addedcolor}{\LocationSet^{temp}} \gets$ $\textcolor{addedcolor}{\LocationSet^{temp}}~\cup$ \{$\bustrip{\busIndex}{\tripIndexJ, end}$\}}
        }
$\levelVar \gets \levelVar - \energyTrip{\bustrip{\busIndex}{\tripIndexJ+1}}$\;

\lIf{$\levelVar < \minenergylevel$}{
        \Return \textbf{false}, $\textcolor{addedcolor}{\LocationSet^{temp}}$}}   
$\levelVar \gets \levelVar -\energyDeadhead{\bustrip{\busIndex}{\numTripsBusb}}{\bustrip{\busIndex}{\numTripsBusb+1}}$\;
    
\lIf{$\levelVar < \minenergylevel$}{
        \Return \textbf{false}, $\textcolor{addedcolor}{\LocationSet^{temp}}$
}
\Return \textbf{true}, $\textcolor{addedcolor}{\LocationSet^{temp}}$
\end{algorithm}

In the concurrent scheduler-based algorithm, checking the charge level feasibility for buses while inserting a trip to an existing bus rotation is paramount for generating feasible rotations. This check is performed by Algorithm \ref{alg:feasibility_check}, which takes the bus rotation, current charging locations, and the energy consumption of service and deadheading trips for the given scenario as inputs. The algorithm first checks if the bus can be charged at the end-stop of the current trip, provided a charging station is located there (lines 4--5). If not, it proceeds to check if the start-stop of the next trip is a charging location and whether the bus can charge there (lines 6--7). If neither of the above two conditions are met, the bus proceeds to perform the next trip, provided it does not run out of energy (lines 8--9). If not, a charging location is opened at the current trip end according to the CAG strategy (lines 10--12). \textcolor{addedcolor}{Opening a charging location at the current trip end is sufficient to ensure the feasibility of the subsequent trip; hence, the algorithm does not consider adding charging locations at the end stop of the current trip and the start stop of the next trip.} A charging station is opened only if the bus can complete the next trip by charging at the maximum rate during the layover (line 12). \textcolor{addedcolor}{Here, the available charging duration during the layover excludes the deadhead travel time between the end of the current trip and the start of the next trip.} The algorithm also checks the energy level of the bus after the completion of the next trip (lines 13--14) and the deadheading trip to the depot (lines 15--16). \textcolor{addedcolor}{We ensure that buses have sufficient energy to reach the depot after completing their final trip; therefore, no additional charging locations are introduced beyond the last trip.} The algorithm returns a boolean variable indicating the bus rotation's charge feasibility and the charging stations' updated configuration. The pseudocode presented here is slightly abridged for ease of explanation. In the experiments, when terminal stops are clustered, we allow buses to deadhead to a nearby charging location within a cluster.

\begin{figure}[h]
    \centering
    \includegraphics[width=0.95\textwidth]{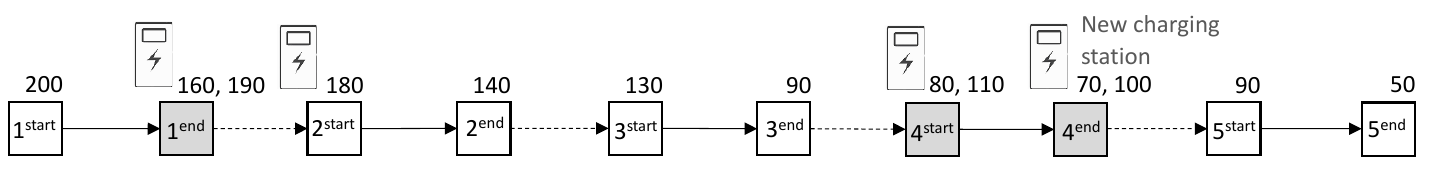}
    \caption{\textcolor{addedcolor}{Example to illustrate Algorithm 2. The numbers next to nodes indicate the charge levels upon reaching them (before and after charging in case of two entries). Gray nodes represent locations where the bus charges according to the CAG strategy}}
    \label{fig:example_is_rotation_charge_feasible}
\end{figure}

Consider the bus rotation in Figure \ref{fig:example_is_rotation_charge_feasible} \textcolor{addedcolor}{to illustrate the workings of Algorithm \ref{alg:feasibility_check}}. The rotation has five service trips $\{1, \ldots, 5\}$. Assume that (a) the minimum and maximum charge levels for the bus are $50$ kWh and $250$ kWh, respectively, (b) each trip consumes $40$ kWh, (c) each deadheading move between two consecutive trips consumes $10$ kWh, (d) the bus starts with an initial charge of $200$ kWh, (e) the EVSP solution allows increasing the charge level by up to $30$ kWh at a charging location, and (f) the existing charging locations are at $1^{end}$, $2^{start}$, and $4^{start}$. The number above the nodes indicates the charge level of the bus when it arrives at that location \textcolor{addedcolor}{and the charge level after finishing charging, if applicable}. When the bus reaches $1^{end}$, its charge level is $160$ kWh. Since charging is available at $1^{end}$, the algorithm requires the bus to be charged at that location, bringing the charge level to $190$ kWh. No charging is done at $2^{start}$ because of the CAG strategy. Note that other bus rotations can use the charging station at $2^{start}$. When the bus reaches $4^{start}$ after two service trips and three deadheading moves, the charge level reduces to $80$ kWh. Since a charging facility is available at $4^{start}$, the charge level gets updated to $110$ kWh at $4^{start}$ and $70$ kWh when the bus reaches $4^{end}$. \textcolor{addedcolor}{As} the bus requires 50 kWh to deadhead and finish trip $5$, the charge level would drop below the minimum threshold. Hence, the algorithm suggests opening a new charging station at $4^{end}$, making the bus rotation charge feasible.

\section{Scenario-based Charge Scheduling Problem} 
\label{sec:scenario_csp}

This section presents an LP model for the scenario-based CSP with dynamic ToU electricity prices, temperature-based energy consumption, and renewable energy production at charging locations. \textcolor{addedcolor}{We note that the EVSP and CSP are solved sequentially; while the initial bus schedules and charging locations are generated using a heuristic, they do not prescribe a specific charging behavior. Instead, the CSP optimally determines charging decisions given these inputs, retaining flexibility in charging profiles and timing within the feasible schedule.} For each scenario, we divide a day into one-minute time steps. \textcolor{addedcolor}{Charging decisions are modeled at this discrete time-slot resolution, which may allow non-continuous charging within a charging opportunity. The same modeling approach is adopted for BESS charging and discharging, where power levels can vary across consecutive time steps. Enforcing continuous charging can be sub-optimal and would require additional constraints that would affect tractability.} A charging opportunity can encompass multiple time steps in which the bus can be charged from the grid, BESS, or both, and the energy drawn from either source is considered a decision variable. We also decide the energy transferred from the grid to the BESS at every time step. The scenario-based LP formulation is presented in Section \ref{sec:formulation}, and the Benders' decomposition approach is outlined in Section \ref{sec:benders}. We make the following assumptions in the CSP, some of which can be trivially relaxed. We also discuss how we obtain estimates for the generated solar energy and the trip energy consumption values under different scenarios.

\textbf{Assumptions:}
\begin{itemize}
\item The CSP is optimized by the bus agency, which \textcolor{addedcolor}{bears} the initial investment cost for the contracted station power capacity, solar panels, BESS, and associated inventory.
\item \textcolor{addedcolor}{We assume that the agency already owns the EB fleet and hence fleet size is predetermined; therefore, bus acquisition costs are not included in the model.}
\item Availability of a sufficient number of chargers is assumed at the charging locations, i.e., we do not explicitly model charger assignment. Instead, we focus on estimating the aggregate power needs at a charging location. Such locations are assumed to be equipped with smart charging technology, which allows split charging \citep{sadeghian2022comprehensive, maia2024ev}. 
\end{itemize}

\begin{table}[H]
 \small
  \centering
  \caption{Notation used in the CSP formulation}
\begin{tabular}{p{3.2cm}p{13.0cm}} 
\hline
       \textbf{{Index:}} & \textbf{{Description}}  \\ 
  \hline
    $t$ & Index for a within-day time step\\
    $k$ & Index for charging opportunity\\
    $\locationindex$ & Index for location\\
    \hline
  \textbf{{Sets:}}  \\ 
  \hline
    $\busSet{\scenario}$ & Set of buses for scenario $\scenario$ \\
     $\timeofoperation{\scenario}$ & Set of time steps for scenario $\scenario$, where each day's operations are discretized at a one-minute granularity. \\
    $\OpportunitiesSet{\scenario}$ & Set of charging opportunities for bus $b$ for scenario $\scenario$ \\
    $\buscount{\locationindex}{\scenario}$ & Set of buses present at location $\locationindex$ in time step $t$ for scenario $\scenario$ \\
    $\settimesteps{\scenario}$ & Set of time steps for which bus $b$ can charge at its $k^\text{th}$ charging opportunity for scenario $\scenario$  either from the power grid or BESS\\
    \hline 

     \multirow{1}{*}{\textbf{{Decision Variables:}}} \\
    \hline
    
    $\buysellenergyscenario{\locationindex}{\scenario}$ & Energy (in kWh/min) transferred from the grid to the bus $b$ at charging location $\locationindex$ in time step $t$ for scenario $\scenario$ \\
    $\buysellenergysolarscenario{\locationindex}{\scenario}$ & Energy (in kWh/min) transferred from the BESS to the bus $b$ at charging location $\locationindex$ in time step $t$ for scenario $\scenario$ \\
    $\gridtobattery{\locationindex,}{t}{\scenario}$ & Energy transferred (in kWh/min) from the grid to the BESS at charging location $\locationindex$ in time step $t$ for scenario $\scenario$\\
    $\contractedcapVar{\locationindex}$ & Charging station power capacity (in kW) at location $\locationindex$ \\
    $\areapanel$ & Total solar panel area (in \si{m^{2}}) installed at location $\locationindex$ \\
$\maxenergylevelsolar{\locationindex}$ & BESS capacity (in kWh) at charging location $\locationindex$ \\
$\initialbatteryenergy$ & Energy level (in kWh) of the BESS at location $\locationindex$ at the start and end of daily operations across all scenarios \\ 
\hline
     \multirow{1}{*}{\textbf{{Auxiliary Variables:}}} \\
    \hline
     $\energytrack{k}{\scenario}$ & Energy level (in kWh) of bus $b$ before the start of its $(k + 1)^\text{th}$ charging opportunity for scenario $\scenario$ \\
  $\energytracksolar{\locationindex,}{t}{\scenario}$ & Energy level (in kWh) of the BESS at charging location $\locationindex$ at the end of time step $t$ for scenario $\scenario$\\    
  \hline
  \textbf{{Data:}} \\ 
  \hline
     $\location$ & Location of bus $b$ at its $k^\text{th}$ charging opportunity\\
   $\solarenergy{\scenario}$ & GTI values (in kWh/\si{m^2}) at charging location $\locationindex$ in time step $t$  $\in$  $\timeofoperation{\scenario}$ for scenario $\scenario$\\
     $\energyrequired{\scenario}$ & Energy required (in kWh) by bus $b$ after its $k^\text{th}$ charging opportunity to complete the subsequent service and deadheading trips until its $(k+1)^\text{th}$ charging opportunity for scenario $\scenario$. This is calculated using $\energyTrip{\tripIndexI}$ and $\energyDeadhead{\tripIndexI}{\tripIndexJ}$. \\ 
    $\probability$ & Probability of the occurrence of scenario $\scenario$\\ 
  $\lasttimestep$ & Time step corresponding to the end of daily operations at location $\locationindex$ for scenario $\scenario$ (to track BESS’s energy level). We set this value to $1440$. \\
    \hline
    \textbf{{Parameters:}}  \\
    \hline
  $\electricitygridprice{t}{\scenario}$ & Grid ToU electricity price per unit energy (in \$/kWh) in time step $t$ for scenario $\scenario$\\
  $\batterycost$ & Unit cost of \textcolor{addedcolor}{BESS} (in \$/kWh) \\
  $\contractedcapcost$ & Unit cost of charging station power capacity (in \$/kW) at potential charging locations \\
  $\costsolarpanel$ & Cost of solar panel per unit area (in $\$/\si{m^{2}}$) \\
  $\efficiencysolarpanel$ & Efficiency (\%) of solar panels \\   
  $DoD$ & DoD of BESS \\
  $\textcolor{addedcolor}{\vartheta}$ & \textcolor{addedcolor}{Amortized Interest Rate (\%)}\\
\hline
\end{tabular}
\label{tab:notation_charging}%
\end{table}

\begin{itemize}
\item Buses are charged from the BESS or the grid. The BESS can charge and discharge simultaneously through a hybrid inverter \citep{bess_components_1}. Commercial BESS units can be of different capacities \citep{canadian_solar, bess_capacity_2}, and their charging/discharging speeds depend on their C-rate \citep{bess_transfer_1}. The maximum energy that can be transferred per minute from the grid and the BESS to a bus is assumed to be the same \citep{vivas2022battery}. The BESS is assumed to be charged both from solar power and \textcolor{addedcolor}{from} the grid \citep{coelho2022unified}. \textcolor{addedcolor}{We do not track the source of energy (solar or grid) discharged from the BESS; this classification is used only for visualization while presenting results.} BESS cannot be fully discharged as we assume a $90\%$ DoD \citep{tejero2023analysis} for batteries. DoD indicates the percentage of a battery's total capacity that can be used before recharging. 
\item All buses are homogeneous and have sufficient charging opportunities to complete daily operations. They can be partially charged multiple times with a piecewise-linear charging profile, as seen in Figure \ref{fig:busoperation}. To allow continuity in operations, we assume that a bus starts and ends its schedule at the same depot. 
\item We assume that the fleet size prescribed by the CS-based algorithm is adequate for serving trips across all scenarios. \textcolor{addedcolor}{We emphasize that fleet size is treated as an exogenous input in our framework, derived from the upstream vehicle scheduling step for each scenario. This modeling choice reflects the perspective of a transit operator who has already procured a fleet and seeks to optimize charging decisions by integrating renewable energy generation.} As the energy consumption per trip varies across scenarios, the number of buses required may also differ. Buses may be redistributed across depots for different scenarios (i.e., across weeks, months, etc.) through rebalancing if needed.
\textcolor{addedcolor}{\item While the number of buses may vary across scenarios due to differences in energy consumption and operational requirements, these variations are not interpreted as acquisition decisions but as scenario-specific realizations of fleet utilization. Incorporating fleet sizing as a decision variable requires a fully integrated EVSP-CSP model with additional integer variables and capital cost components. We leave this extension as a direction for future work.}
\item Sufficient space for solar panel installation is assumed at the charging locations. If the panel area requirement is too high, panels can be installed on the rooftops of nearby buildings, or additional constraints can be imposed in the model. The energy transfer from the panel to BESS happens instantaneously \citep{mannepalli2022allocation}.
\textcolor{addedcolor}{\item Battery degradation is not modeled explicitly. Incorporating cycle-dependent degradation would require tracking battery usage and the evolution of state-of-health, thereby introducing additional nonlinearities and significantly increasing the complexity of the optimization model. The inclusion of detailed degradation dynamics is therefore left for future work.}
\end{itemize}

\textbf{Energy Consumption Estimates:} 
The energy consumption per trip is estimated using a linear regression equation from \cite{JI2022100069} as follows:
    \begin{align}
    \label{trip_energy_consumption}
        \text{Energy per trip $i$ (in kWh)} = \exp{\big({\alpha_{0}} + {\alpha_{1}} \log{L_{i}} +{\alpha_{2}}\log{M_{i}} +{\alpha_{3}}\log{t_{i}} +{\alpha_{4}} {|\overline{T}_{i} - \overline{T}^{*}_{g}|\big)}}
    \end{align}
where $L_{i}$ is the length (in km) of trip $i$, $M_{i}$ is the mass (in kg) of the EB serving trip $i$, $t_{i}$ is the travel time (in minutes) of trip $i$, $\overline{T}_{i}$ (in $^{\circ}C$) is the average ambient temperature during trip $i$, $\overline{T}^{*}_{g}$ (in $^{\circ}C$) is the optimum working temperature of the EB, and ${\alpha_{0}}, \ldots, {\alpha_{4}}$ are estimated linear regression parameters. The ambient temperature for a trip is calculated by averaging the historical temperature values at its start and end stops. We used the same value of $M_{i}$ across all trips. Additional ridership data, if available, can also be directly plugged into this equation to obtain a more realistic estimate of the energy consumed during a trip. The effects of travel time uncertainty due to congestion and other factors that may influence EB energy consumption are beyond the scope of this research.

\textcolor{addedcolor}{Unlike standard formulations that assume constant or linearly approximated energy consumption, equation \eqref{trip_energy_consumption} introduces a non-linear dependence on ambient temperature. This has important implications for infrastructure planning. First, it leads to asymmetric energy consumption patterns, where deviations from the optimal temperature $\overline{T}^{*}_{g}$, both in colder and hotter conditions, can significantly increase energy requirements compared to a non-linear temperature-independent model, obtained by setting $\alpha_4 = 0$ in equation \eqref{trip_energy_consumption}. In contrast, linear approximations based solely on distance traveled may misestimate energy requirements, leading to over-investment in charging infrastructure and energy resources. For example, using a fixed energy consumption rate of 1.5 kWh/km yields average energy consumption estimates that are more than twice those obtained with the proposed temperature-dependent formulation.}

\textbf{Solar Energy Estimates:}
The solar energy per minute is estimated using the following equation from \cite{nrel}: 
\begin{align}
\label{solar_energy_production}
\text{Energy per minute (in kWh/min)} &= \text{GTI} \times \text{Panel Area} \times \frac{\eta} {100} \times \frac{1}{60}  
\end{align}
where GTI is in kW/\si{m^2}, Panel Area is in \si{m^2}, and $\eta$ is in \%.

\subsection{Formulation}
\label{sec:formulation}

Table $\ref{tab:notation_charging}$ summarizes the notation used in the CSP model. The long-term decision variables determine the power capacity ($\contractedcapVar{\locationindex}$), solar panel area ($\areapanel$), the BESS capacity ($\maxenergylevelsolar{\locationindex}$), and the initial energy level for the BESS ($\initialbatteryenergy$) for every charging location $\locationindex$. The scenarios-specific daily operational decisions include determining the energy transferred from the grid ($\buysellenergyscenario{\locationindex}{\scenario}$) and BESS ($\buysellenergysolarscenario{\locationindex}{\scenario}$) to every bus $b$, and the energy transferred from the grid to BESS ($\gridtobattery{\locationindex,}{t}{\scenario}$) at charging location $\locationindex$ and time step $t$ within each day for scenario $\scenario$. Additionally, auxiliary variables are used to track the energy level ($\energytrack{k}{\scenario}$) of bus $\busIndex$ \textcolor{addedcolor}{before the start of its $(k + 1)^{\text{th}}$} charging opportunity, and the energy level of BESS ($\energytracksolar{\locationindex,}{t}{\scenario}$) at charging location $\locationindex$ and \textcolor{addedcolor}{at the end of} time step $t$ for scenario $\scenario$.
 
\begin{figure}[h]

    \centering
    \includegraphics[width=\textwidth]{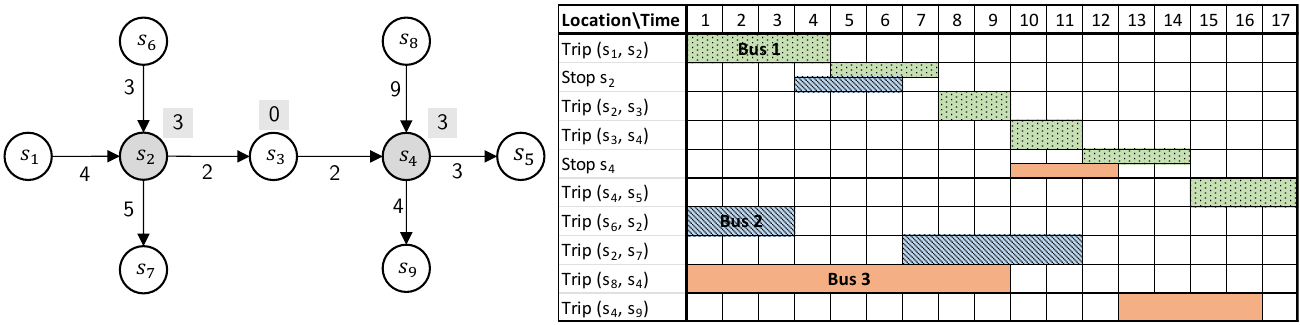}
    \caption{\textcolor{addedcolor}{Example network (left panel) and corresponding Gantt chart (right panel), illustrating bus-to-trip assignments, which act as the input data for the CSP}}
    \label{fig:example_problem_csp_network}
\end{figure}

\begin{table}[ht]
    \centering
    \caption{Buses present at a charging location over time and time steps for charging opportunities in the example shown in Figure \ref{fig:example_problem_csp_network}}
    \begin{tabular}{lccccccccc}
    \hline
    \textbf{Location} $(j)$ & $s_{2}$ &  $s_{2}$ & $s_{2}$ & $s_{2}$ & $s_{4}$ & $s_{4}$ & $s_{4}$ & $s_{4}$ & $s_{4}$ \\
    \hline
    Time $(t)$ & 4 & 5 & 6 & 7 & 10 & 11 & 12 & 13 & 14\\
    $\buscount{\locationindex}{\scenario}$ & \{2\} & \{1,2\} & \{1,2\} & \{1\} & \{3\} & \{3\} & \{1,3\} & \{1\} & \{1\} \\
    \hline
    \end{tabular}
    \vspace{0.5mm}
    \begin{tabular}{lcccc}
    \hline
    \textbf{Bus} $(b)$ & 1  &  1 & 2 & 3 \\
    \hline
    Charging Opportunity $(k)$ & 1 & 2 & 1 & 1\\
    $\settimesteps{\scenario}$ & \{5,6,7\} & \{12,13,14\} & \{4,5,6\} & \{10,11,12\} \\
    \hline
    \end{tabular}
    \label{tab:combined_example_problem_csp}
\end{table}

To illustrate the sets used in the CSP formulation, consider the network shown in the left panel of Figure \ref{fig:example_problem_csp_network}, which depicts eight trips covered by three buses. The stops $s_1\, \ldots, s_8$ indicate the end terminals of the trips. The numbers next to the links represent the trip travel times, and the numbers next to the nodes in boxes denote the layover times between trips, respectively. Note that the layover times for different trips ending at the same stop are identical. Additionally, we assume that the three buses start their journey simultaneously. Charging facilities are present at stops $s_{2}$ and $s_{4}$, marked in gray. Bus 1 covers four trips---$(s_1, s_2)$, $(s_2, s_3)$, $(s_3, s_4)$, and $(s_4, s_5)$; bus 2 covers two trips---$(s_6, s_2)$ and $(s_2, s_7)$; and bus 3 covers two trips---$(s_8, s_4)$ and $(s_4, s_9)$. Bus 1 has two charging opportunities at $s_2$ and $s_4$. Buses 2 and 3, on the other hand, have a single charging opportunity at $s_2$ and $s_4$, respectively. The bus schedule is shown in the Gantt chart in the right panel of Figure \ref{fig:example_problem_csp_network} for a scenario $\scenario$. For the CSP, we need to keep track of events that occur at the charging stations. Table \ref{tab:combined_example_problem_csp} provides the information on the set of buses present at each location in different time steps, $\buscount{\locationindex}{\scenario}$, and the time steps corresponding to different charging opportunities for each bus, $\settimesteps{\scenario}$.

The objective function \eqref{eqn:objective} for the CSP consists of \textcolor{addedcolor}{five} terms. The first three terms represent the per day \textcolor{addedcolor}{amortized} cost of the BESS capacity, the contracted station power capacity, and the installed solar panel area, respectively. These decisions are related to the bus agency's long-term initial investments. The last \textcolor{addedcolor}{two} terms are expected values of scenario-specific terms: \textcolor{addedcolor}{the first for the energy drawn from the grid to the buses, and the final for the energy transferred from the grid to the BESS.} \textcolor{addedcolor}{Energy stored in the BESS may originate from solar or grid sources. Solar energy costs are captured through the amortized PV and BESS investment costs, while grid-sourced energy is priced at the time of charging buses/BESS.} The objective denotes the daily average cost. The occurrence of each scenario is equally likely since we assume that the scenarios correspond to time intervals of equal duration. We use a life-cycle of 12 years for BESS \citep{bess_4, rallo2020lithium} and the bus battery charging systems \citep{lajunen2018lifecycle}, and 30 years for the solar panels \citep{pv_life}, which is reflected in the constants used in the objective function.

 \textcolor{addedcolor}{
 \begin{align}
   &\min \,  \, 
   \frac{\batterycost}{365}\frac{\vartheta(1+\vartheta)^{12}}{(1+\vartheta)^{12}-1} \, \sum\limits_{\locationindex \in J} \maxenergylevelsolar{\locationindex} + \frac{\contractedcapcost }{365}\frac{\vartheta(1+\vartheta)^{12}}{(1+\vartheta)^{12}-1} \, \sum_{\locationindex \in J} \contractedcapVar{\locationindex} + \frac{\costsolarpanel}{365}\frac{\vartheta(1+\vartheta)^{30}}{(1+\vartheta)^{30}-1} \, \sum\limits_{\locationindex \in J} \areapanel + \nonumber \\
   &\sum\limits_{\scenario \in \numberofscenario}\probability \bigg(\sum_{b \in \busSet{\scenario}}\sum_{\substack{k \,\in \, \OpportunitiesSet{\scenario} \cup \{|\OpportunitiesSet{\scenario}| + 1\}, \\t \in \settimesteps{\scenario}, \, j'=\location}} \electricitygridprice{t}{\scenario}\, \buysellenergyscenario{\locationindex'}{\scenario} + \sum\limits_{\locationindex \in J} \sum\limits_{t \in \timeofoperation{\scenario}} \electricitygridprice{t}{\scenario} \gridtobattery{\locationindex,}{t}{\scenario} \bigg)
   \label{eqn:objective}
\end{align}}
{\allowdisplaybreaks
\begin{alignat}{2}
\label{eqn:energylevel}
\text{s.t.} \; &\energytrack{k}{\scenario} = \energytrack{k-1}{\scenario} - \energyrequired{\scenario}  + \nonumber\\ &\sum_{t \in \settimesteps{\scenario}}(\buysellenergyscenario{\locationindex'}{\scenario} + \buysellenergysolarscenario{\locationindex'}{\scenario})  \quad && \forall \; \scenario \, \in \, \numberofscenario ,\, b \, \in \, \busSet{\scenario},\, k \in \OpportunitiesSet{\scenario} \textcolor{addedcolor}{\cup \{|\OpportunitiesSet{\scenario}| + 1\}}, \,\locationindex'=\location\\
\label{eqn:maxenergy}
&\energytrack{k-1}{\scenario} + \sum_{t \in \settimesteps{\scenario}}(\buysellenergyscenario{\locationindex'}{\scenario} + \buysellenergysolarscenario{\locationindex'}{\scenario}) \leq \maxenergylevel  \quad &&\forall \; \scenario \, \in \, \numberofscenario , \, b \, \in \, \busSet{\scenario},\, k \in \OpportunitiesSet{\scenario}\textcolor{addedcolor}{\cup \{|\OpportunitiesSet{\scenario}| + 1\}},\, \locationindex'=\location\\
\label{eqn:solarenergycharging}
&\energytracksolar{\locationindex,}{t}{\scenario} = \energytracksolar{\locationindex,}{t-1}{\scenario} + \frac{\efficiencysolarpanel\,\solarenergy{\scenario}\,\areapanel}{60} + \gridtobattery{\locationindex,}{t}{\scenario} - \nonumber\\ &\sum_{b\in \buscount{\locationindex}{\scenario}}\buysellenergysolarscenario{\locationindex}{\scenario}  \quad &&\forall \; \scenario \, \in \, \numberofscenario , \locationindex \, \in \, J, \;  t \, \in \, \timeofoperation{\scenario}\\
\label{eqn:maxenergysolar}
&\maxenergylevelsolar{\locationindex} \geq \energytracksolar{\locationindex,}{t-1}{\scenario} + \frac{\efficiencysolarpanel \,\solarenergy {\scenario}\,\areapanel}{60} + \gridtobattery{\locationindex,}{t}{\scenario} \quad &&\forall \; \scenario \, \in \, \numberofscenario , \locationindex \, \in \, J, \;  t \, \in \, \timeofoperation{\scenario}\\
\label{eqn:depthofdischarge}
&\energytracksolar{\locationindex,}{t}{\scenario} \geq (1 - DoD) \maxenergylevelsolar{\locationindex} \quad &&\forall \; \scenario \, \in \, \numberofscenario , \locationindex \, \in \, J, \;  t \, \in \, \timeofoperation{\scenario}\\
\label{eqn:contractedcapacity}
& \sum_{b \in \buscount{\locationindex}{\scenario}} 60 \, \buysellenergyscenario{\locationindex}{\scenario} + 60 \, \gridtobattery{\locationindex,}{t}{\scenario} \leq \contractedcapVar{\locationindex} \quad &&\forall \;\scenario \, \in \, \numberofscenario , \,   \locationindex \, \in \, J, \;  t \, \in \, \timeofoperation{\scenario}\\ 
\label{eqn:allowedtransfer}
&\buysellenergyscenario{\locationindex'}{\scenario} + \buysellenergysolarscenario{\locationindex'}{\scenario} \leq \energytransferMinute   \quad &&\forall \;\scenario \, \in \, \numberofscenario , \,  b \, \in \,\busSet{\scenario},  k \,\in \, \OpportunitiesSet{\scenario}\textcolor{addedcolor}{\cup \{|\OpportunitiesSet{\scenario}| + 1\}},\,t \,\in \, \settimesteps{\scenario} \,, \, \locationindex'=\location\\
\label{eqn:minenergy}
&\energytrack{k}{\scenario} \geq \minenergylevel \quad &&\forall \; \scenario \, \in \, \numberofscenario , \, b \, \in \, \busSet{\scenario}, \, k \in \OpportunitiesSet{\scenario}\textcolor{addedcolor}{\cup \{|\OpportunitiesSet{\scenario}| + 1\}}\\ 
\label{eqn:initialenergy}
&\energytrack{0}{\scenario} + e^{\scenario}_{b,0} = \maxenergylevel \quad &&\forall \;\scenario \, \in \, \numberofscenario , \,   b \, \in \, \busSet{\scenario}\\
\label{eqn:energy_consistency_buses} 
&\energytrack{0}{\scenario} = \energytrack{|\OpportunitiesSet{\scenario}| + 1}{\scenario}\quad &&\forall \;\scenario \, \in \, \numberofscenario , \,   b \, \in \, \busSet{\scenario}\\
\label{eqn:initialenergysolar}
& \energytracksolar{\locationindex,}{0}{\scenario} = \energytracksolar{\locationindex,}{\lasttimestep}{\scenario} = \initialbatteryenergy && \forall \; \scenario \, \in \, \numberofscenario , \,  \locationindex \, \in \, J\\
\label{eqn:nonnegativeconstraint}
&\buysellenergyscenario{\locationindex'}{\scenario} \geq 0, \; \buysellenergysolarscenario{\locationindex'}{\scenario} \geq 0  && \forall \; \scenario \, \in \, \numberofscenario , \, b \, \in \, \busSet{\scenario}, \, k \,\in \, \OpportunitiesSet{\scenario}\cup \{|\OpportunitiesSet{\scenario}| + 1\},\,t \,\in \, \settimesteps{\scenario} \,, \, \locationindex'=\location\\
\label{eqn:energylevelsolarnonnegative}
& \areapanel \geq 0, \; \maxenergylevelsolar{\locationindex} \geq 0, \;  \contractedcapVar{\locationindex} \geq 0, \; \initialbatteryenergy \geq 0 && \forall \; \locationindex \, \in \, J\\
& \gridtobattery{\locationindex,}{t}{\scenario} \geq 0 && \forall \; \scenario \, \in \, \numberofscenario , \, \locationindex \in \LocationSet, \,t \,\in \timeofoperation{\scenario} \label{eqn:auxvar}
\end{alignat}
 
\textcolor{addedcolor}{Charging opportunity index $|\OpportunitiesSet{\scenario}| + 1$ represents the overnight charging opportunities for bus $b$ under scenario $\scenario$. (Charging opportunity indexing in Figures $\ref{fig:busoperation}$ and $\ref{fig:example_problem_csp_network}$ starts from $1$, and the overnight charging activity is not shown.)} \textcolor{addedcolor}{Note that $e^{\scenario}_{b,0}$ and $e^{\scenario}_{b,|\OpportunitiesSet{\scenario}| + 1}$ denote the energy required by bus $b$ after its overnight charging opportunity to complete the subsequent service and deadheading trips until its first charging opportunity.} \textcolor{addedcolor}{In our formulation}, all the constraints except \eqref{eqn:energylevelsolarnonnegative} are dependent on the scenarios. Constraint \eqref{eqn:energylevel} tracks the energy level of each bus at the end of every charging opportunity.  \textcolor{addedcolor}{Tracking energy levels at every time step is unnecessary and would significantly increase the number of decision variables.} Constraint \eqref{eqn:maxenergy} ensures that buses do not charge above their battery's maximum capacity. Constraint \eqref{eqn:solarenergycharging} tracks the energy level of the BESS for all locations, calculated using the energy drawn from the solar panels and the grid and the energy transferred to buses. Constraint \eqref{eqn:maxenergysolar} models the maximum capacity of the BESS. \textcolor{addedcolor}{We note that the energy discharged from the BESS to buses within the same time step is not included in this constraint, as we assume that energy is first stored in the BESS from solar and grid sources before any discharge occurs. Therefore, the capacity limit is enforced on the pre-discharge energy level, and no additional discharge term is required.} Constraint \eqref{eqn:depthofdischarge} sets the lower bound of the energy level of BESS based on its DoD. Constraint \eqref{eqn:contractedcapacity} determines the power capacity of charging locations by calculating the peak energy requirements. \textcolor{addedcolor}{All power variables are converted to energy (kWh) using the one-minute time step, and hence all state variables represent energy, ensuring consistent energy balance calculations in \eqref{eqn:solarenergycharging}–\eqref{eqn:contractedcapacity}.} We ensure unit consistency in constraint \eqref{eqn:contractedcapacity} by applying a factor of 60 to convert kWh/min to kW. Constraint \eqref{eqn:allowedtransfer} restricts the maximum energy transferred to the bus from the power grid or BESS in any time step. Constraint \eqref{eqn:minenergy} ensures that the SoC of the bus battery never drops below a minimum level. \textcolor{addedcolor}{Constraint \eqref{eqn:initialenergy} sets the initial energy level of each bus $b$ before starting their daily operations to 85\% of battery's capacity to implement recurring schedules in each scenario.} \textcolor{addedcolor}{We assume a homogeneous fleet and fix the initial energy level of all buses to 85\% of their battery capacity, independent of the scenario.} \textcolor{addedcolor}{Constraint \eqref{eqn:energy_consistency_buses}, together with \eqref{eqn:energylevel}, conserves the energy levels of the buses after overnight charging.} Constraint \eqref{eqn:initialenergysolar} similarly equates the starting energy level of batteries on a given day to that on the subsequent day. \textcolor{addedcolor}{This ensures consistent initialization across scenarios.} \textcolor{addedcolor}{We enforce identical starting and ending energy levels across scenarios to ensure cyclic operation of the system. Allowing scenario-dependent terminal energy levels for buses and BESS would violate the continuity/cyclicity of daily operations assumed in the model.} Finally, constraints \eqref{eqn:nonnegativeconstraint}, \eqref{eqn:energylevelsolarnonnegative}, and \eqref{eqn:auxvar} require variables to be non-negative.

\textcolor{addedcolor}{Scenario-based optimization improves long-term decisions under varying solar energy availability and temperature-dependent energy consumption.} However, as we increase the number of scenarios, the number of variables and constraints increases, making it difficult to solve the LP for large networks using the simplex method. Since the constraint matrix has a block diagonal structure, this problem is a good candidate for applying Benders' decomposition.

\subsection{Benders' Decomposition}
\label{sec:benders}
\textcolor{addedcolor}{Benders’ decomposition is particularly suited to our formulation due to the clear separation between first-stage infrastructure decisions and scenario-specific operational subproblems with a block-diagonal structure. In contrast, other decomposition methods or column generation are more effective for path-based VSP formulations and are not naturally aligned with our CSP framework.}
The formulation discussed previously in \eqref{eqn:objective}--\eqref{eqn:auxvar} has a two-stage structure in which the contracted grid capacity, solar-powered BESS capacity, the area of solar panels, and the energy levels of the BESS \textcolor{addedcolor}{and buses at the start of operations} are decided in the first stage. Fine-grained charging decisions are made in the second stage, adapting to the scenario-specific GTI and temperature values. The constraints for this problem can be represented using \eqref{eq:benders}, which has a block diagonal shape with non-overlapping blocks for each scenario, except for the first-stage variables that appear in every scenario. 
\begin{equation}
    \begin{bmatrix}
        A_{00} &        &        & & & \\
        A_{10} & A_{11} &        & & & \\
        A_{20} &        & A_{22} & & & \\
        \vdots &        &        & \ddots & & \\
        A_{\phi 0} &        &        &        & & A_{\phi \theta}
    \end{bmatrix}
    \begin{bmatrix}
        \lambda_0 \\
        \lambda_1 \\
        \lambda_2 \\
        \vdots \\
        \lambda_\theta
    \end{bmatrix} = 
    \begin{bmatrix}
        \beta_0 \\
        \beta_1 \\
        \beta_2 \\
        \vdots \\
        \beta_\phi
    \end{bmatrix}\label{eq:benders}
\end{equation}

Benders' decomposition or the L-shaped method is ideal for improving the tractability of two-stage recourse models, especially when the number of scenarios is high \citep{BirgLouv97}. In this approach, one creates a  \textit{restricted main/parent problem} containing only constraints with first-stage variables. The objective for the restricted main problem includes the long-term variables and a term for each \textit{sub-problem} or \textit{child problem}, which indicates the current assessment (or a lower bound) of the sub-problem's objective. A solution to the restricted main problem allows us to fix the first-stage decisions, which results in an independent sub-problem for each scenario. However, a scenario-specific sub-problem might not be feasible for the current values of the first-stage variables. In such a case, one can use duality theory to obtain a feasibility cut from the sub-problem and add it to the restricted main problem. Alternatively, the main problem's assessment of the objective function value of the scenario-specific sub-problem could be less than the actual objective value of the sub-problem. In this case, we add an optimality cut from the sub-problem to the restricted main problem. The method iteratively solves the restricted main problem and the sub-problems until the main problem does not lead to a feasibility or an optimality cut from any sub-problem.

To apply Benders' decomposition for the CSP, we can formulate the objective function of the main problem as \eqref{eqn:first_stage_objective} where $\zeta^{\scenario}$ is the assessment of the objective value of the second stage for scenario $\scenario$. Thus, $\sum_{\scenario \in \numberofscenario}\probability\zeta^{\scenario}$ denotes the average contribution of the second-stage decisions \textcolor{addedcolor}{to} the overall objective. In our formulation, the first-stage problem contains only non-negativity constraints \eqref{eqn:energylevelsolarnonnegative} in the first iteration. Additional location-specific capacity or budgetary constraints involving the first-stage variables can be easily incorporated if needed.
\textcolor{addedcolor}{{\setlength{\belowdisplayskip}{-5pt}
\setlength{\belowdisplayshortskip}{-5pt}
\begin{align}
   &\min{ \, \left(\frac{\batterycost}{365}\frac{\vartheta(1+\vartheta)^{12}}{(1+\vartheta)^{12}-1} \, \sum\limits_{\locationindex \in J} \maxenergylevelsolar{\locationindex} + \frac{\contractedcapcost }{365}\frac{\vartheta(1+\vartheta)^{12}}{(1+\vartheta)^{12}-1} \, \sum_{\locationindex \in J} \contractedcapVar{\locationindex} + \frac{\costsolarpanel}{365}\frac{\vartheta(1+\vartheta)^{30}}{(1+\vartheta)^{30}-1} \, \sum\limits_{\locationindex \in J} \areapanel \right) + \sum_{\scenario \in \numberofscenario}\probability\zeta^{\scenario}}
   \label{eqn:first_stage_objective}
\end{align}}}

After solving \eqref{eqn:first_stage_objective} subject to \eqref{eqn:energylevelsolarnonnegative}, the values of the first-stage variables (represented with additional bars) are fixed in each scenario-based LP as shown below for a scenario $\scenario$.
\textcolor{addedcolor}{
\begin{align}
   &\min{ \,
   \sum_{b \in \busSet{\scenario}}\sum_{\substack{k \,\in \, \OpportunitiesSet{\scenario}\cup \{|\OpportunitiesSet{\scenario}| + 1\}, \\t \in \settimesteps{\scenario}, \, j'=\location}}\electricitygridprice{t}{\scenario}\, \buysellenergyscenario{\locationindex'}{\scenario} + \sum\limits_{\locationindex \in J} \sum\limits_{t \in \timeofoperation{\scenario}} \electricitygridprice{t}{\scenario} \gridtobattery{\locationindex,}{t}{\scenario}} 
   \label{eqn:second_stage_objective}
\end{align}}
\begin{alignat}{3}
\label{eqn:solarenergycharging_benders}
\text{s.t.} \; & \eqref{eqn:energylevel}, \eqref{eqn:maxenergy}, \eqref{eqn:allowedtransfer}, \eqref{eqn:minenergy}, \eqref{eqn:initialenergy},
\eqref{eqn:energy_consistency_buses},
\eqref{eqn:nonnegativeconstraint}, \eqref{eqn:auxvar} \, \text{for scenario} \, \scenario  \nonumber \\ &\energytracksolar{\locationindex,}{t}{\scenario} - \energytracksolar{\locationindex,}{t-1}{\scenario}  + \sum_{b\in \buscount{\locationindex}{\scenario}}\buysellenergysolarscenario{\locationindex}{\scenario} - \gridtobattery{\locationindex,}{t}{\scenario} = \frac{\efficiencysolarpanel\,\solarenergy{\scenario}\,{\Bar{a}_{\locationindex}}}{60} \quad &&\forall \; \locationindex \, \in \, J, \;  t \, \in \, \timeofoperation{\scenario}  \\
\label{eqn:maxenergysolar_benders}
& \energytracksolar{\locationindex,}{t-1}{\scenario} + \gridtobattery{\locationindex,}{t}{\scenario} \leq  {\Bar{c}_{\locationindex}} - \frac{\efficiencysolarpanel \,\solarenergy {\scenario}\,{\Bar{a}_{\locationindex}}}{60} \quad &&\forall \; \locationindex \, \in \, J, \;  t \, \in \, \timeofoperation{\scenario}  \\
\label{eqn:depthofdischarge_benders}
&\energytracksolar{\locationindex,}{t}{\scenario} \geq (1 - DoD) {\Bar{c}_{\locationindex}} \quad &&\forall \; \locationindex \, \in \, J, \;  t \, \in \, \timeofoperation{\scenario}  \\
\label{eqn:contractedcapacity_benders}
 & \sum_{b \in \buscount{\locationindex}{\scenario}} 60 \, \buysellenergyscenario{\locationindex}{\scenario} + 60 \, \gridtobattery{\locationindex,}{t}{\scenario} \leq{\Bar{z}_{\locationindex}} \quad &&\forall \;  \locationindex \, \in \, J, \;  t \, \in \, \timeofoperation{\scenario} \\
\label{eqn:initialenergysolar_benders}
& \energytracksolar{\locationindex,}{0}{\scenario} = \energytracksolar{\locationindex,}{\lasttimestep}{\scenario} = \Bar{d}_{\locationindex} && \forall \,  \locationindex \, \in \, J 
\end{alignat}

\textcolor{addedcolor}{
\textbf{Update of $\zeta^{\scenario}$ and Benders Cuts:}
For each scenario $\scenario$, the corresponding second-stage LP is solved. Let $\boldsymbol{\psi}$ denote the vector of all first-stage variables, $\boldsymbol{\Phi}^{\scenario}$ denote the scenario-specific right-hand-side vector, and $\boldsymbol{\Gamma}^{\scenario}$ denote the matrix of coefficients associated with the first-stage variables. The subproblem constraints include $\leq$, $=$, and $\geq$ type constraints, whose associated dual variables are non-positive, unrestricted in sign, and non-negative, respectively. Based on the solution of the subproblem, either an optimality cut or a feasibility cut is generated and added to the first-stage problem.
}
\vspace{-5mm}
\textcolor{addedcolor}{
\paragraph*{Optimality cuts:}
If the subproblem corresponding to scenario $\scenario$ is feasible, let $\boldsymbol{\lambda}^{\scenario}$ denote an optimal dual solution (with sign conventions as above). The corresponding optimality cut is
\begin{align}
\zeta^{\scenario}
\ge
(\boldsymbol{\lambda}^{\scenario})^{\top}
\left(
\boldsymbol{\Phi}^{\scenario}
-
\boldsymbol{\Gamma}^{\scenario}\boldsymbol{\psi}
\right)
\end{align}}
\vspace{-10mm}
\textcolor{addedcolor}{{\setlength{\belowdisplayskip}{-5pt}
\setlength{\belowdisplayshortskip}{-5pt}
\paragraph*{Feasibility cuts:}
If the subproblem is infeasible, let $\hat{\boldsymbol{\lambda}}^{\scenario}$ denote an extreme ray of the dual feasible cone (again with appropriate sign conventions). The corresponding feasibility cut is
\begin{align}
(\hat{\boldsymbol{\lambda}}^{\scenario})^{\top}
\left(
\boldsymbol{\Phi}^{\scenario}
-
\boldsymbol{\Gamma}^{\scenario}\boldsymbol{\psi}
\right)
\le 0
\end{align}
}}

If a scenario-based LP is infeasible, a feasibility cut is added to the first-stage problem. On the other hand, if the LP is feasible but has an objective \eqref{eqn:second_stage_objective} that is not close enough to $\zeta^{\scenario}$, then an optimality cut is added to the first-stage problem. As discussed earlier, utilizing this iterative approach improves tractability. The following section presents a case study to demonstrate its advantages.

\section{Experiments and Results}
\label{sec:results}
Our experiments were run on two real-world bus transit networks using a Dell Precision workstation with \textcolor{addedcolor}{Intel(R) Core(TM) i9-14900K} CPU @ \textcolor{addedcolor}{3.20 GHz} and \textcolor{addedcolor}{128 GB of RAM}. The LP models were implemented in Python 3.10 with CPLEX 22.1.1. We used CPLEX's built-in Benders' decomposition implementation with our annotations for the variables and benchmarked it against CPLEX's (dual simplex) LP solver. The maximum time limit for both models was set to \textcolor{addedcolor}{36 hours} for each instance. 

\subsection{Dataset and Parameters}
We tested our models on the Durham Transit Network, Ontario, Canada, and Action Buses, Canberra, Australia. The General Transit Feed Specification (GTFS) data \textcolor{addedcolor}{from 2023} was collected from \href{https://transitfeeds.com/}{Transit Feeds}. Our model requires agency-specific input data, such as fleet characteristics, potential charging locations, and depots. However, in the absence of such data, we made reasonable assumptions, where necessary, to demonstrate the advantages of our framework. The actual operations in these cities may, however, differ. Using the calendar file and service IDs, we filtered the highest number of service trips occurring on a weekday and assumed that these trips repeat daily throughout the year. Since EBs require overnight charging, we only considered trips that started before 10 pm for Durham and before 11 pm for the Canberra network to allow buses sufficient time to recharge their batteries before the next day's operation. Alternatively, one could use extra buses to complete the overnight trips and recharge them during the daytime. 

Distances between bus stops were calculated using  \href{https://project-osrm.org/}{Open Source Routing Machine}, which was also used to calculate the trip and deadheading distances. For deadheading, a constant speed of 30 km/hr was assumed. We considered the specifications of Ride's \href{https://ride.co/products/k9m/}{K9M-40} EB model for our study. To estimate the energy consumed during a trip using \eqref{trip_energy_consumption}, we assumed that buses carry 50\% of their passenger capacity. With actual ticketing data, a more realistic estimate of trip energy consumption can be achieved. However, the curb weight of an EB is much higher than the total passenger weight \citep{JI2022100069}. Hence, approximating passenger loads is unlikely to affect the results. We formed clusters of terminal stops within a radius of 500 m to limit the number of charging locations. For each cluster, we chose the stop with the maximum number of trips passing through it as a candidate charging location. \textcolor{addedcolor}{To assess the impact of the clustering assumption, we conducted an additional sensitivity analysis by increasing the radius from 500 m to 750 m. The results of this experiment are included in Section \ref{sec:clustering_sensitivity}.} 

\begin{table}[H]
 \small
  \centering
  \caption{Parameter values}
\begin{tabular}{p{2.1cm}p{2.35cm}p{3cm}p{7.95cm}} 
    \hline
    \textbf{Parameter} & \textbf{Value} & \textbf{Reference}& \textbf{Remark/Assumptions} \\
    \hline
        $\maxenergylevel$ & $266.05$ kWh & \href{https://ride.co/products/k9m/}{K9M-40}, \cite{stumpe2021study} & 85\% battery capacity\\
  $\minenergylevel$ & $46.95$ kWh & \href{https://ride.co/products/k9m/}{K9M-40}, \cite{sadati2019operational} & 15\% battery capacity\\
  $\efficiencysolarpanel$ & $20$\% & \cite{tiano2020evaluation} & Photovoltaic panel efficiency \\
  $\energytransferMinute{}$ & $2.5$ kWh/min &  \citep{jahic2019charging, vivas2022battery} & Allowed energy transfer rate to buses\\
  $\batterycost$ & $500$ \$/kWh & \cite{nrel_report_bess} (pg.5) & 12 years life span \citep{bess_4, rallo2020lithium}\\
  $\contractedcapcost$ & $654$ \$/kW &  \cite{dirks2022integration} & 12 years life span \citep{lajunen2018lifecycle}\\
  $\costsolarpanel$ & \textcolor{addedcolor}{$256.11$ \$/\si{m^{2}}} & \citep{nrel_report_panel} (pg.28, pg.32) & 30 years life span \citep{pv_life} (panels, inverter, maintenance) \\ 
  $DoD$ & 90\% & \cite{tejero2023analysis} & DoD of BESS \\
  \textcolor{addedcolor}{$\vartheta$} & \textcolor{addedcolor}{3.5\%} & \cite{harris2020probabilistic} & \textcolor{addedcolor}{Interest rate for amortization}\\
    Mass of EB & 16,121.14 kg & \href{https://ride.co/products/k9m/}{K9M-40}&  Curb weight and 50\% passenger load\\
  $\alpha_{0}, \alpha_{1}, \alpha_{2}, \alpha_{3}, \alpha_{4}$ & $-$8.11, 0.55, 0.78, 0.35, 0.008 
  &\cite{JI2022100069} \multirow{2}{*}{ \hspace{5em} {}} &\multirow{2}{*}{Estimated parameters for the energy consumption profile}\\
  $\overline{T}^{*}_{g}$ & $23.3^{\circ}$C & \cite{JI2022100069} & \\
\hline
\end{tabular}
\label{tab:parameters}%
\end{table}

We gathered time-stamped GTI and temperature data at transit stops from PVWatts, a module developed by \href{https://pvwatts.nrel.gov/}{NREL}, which also required the latitude and longitude of the charging station, tilt, and array azimuth angle for installed solar panels. We assumed that the tilt value equals a location's latitude and the array azimuth angle is $180^\circ$. We calculated the power generation per unit area for Durham and Canberra using data from \cite{canada_ontario_solar} and \cite{aus_gov_solar}, respectively, and used these values to convert the solar panel costs to \$/\si{m^{2}} based on \cite{nrel_report_panel}'s \textcolor{addedcolor}{recommended costs per power generation}. The data used in our experiments are summarized in Table \ref{tab:parameters} and were converted to per-day equivalents during implementation to scale long- and short-term model parameters. Table \ref{tab:ToU_pricing} presents the ToU electricity prices obtained from Alectra Utilities \citep{can_tou} (Durham) and Australian Energy Regulator \citep{aus_tou} (Canberra). \textcolor{addedcolor}{The ToU tariffs used for the case studies correspond to publicly available rates for small business customers; additional pricing components such as demand charges can be incorporated if such data is available.} \textcolor{addedcolor}{All monetary values reported in this study are in USD; values originally in CAD and AUD are converted using standard exchange rates.}

\textcolor{addedcolor}{
While electricity costs may include additional components such as global adjustments and peak demand charges, these are often region-specific and strongly influenced by external market conditions and geopolitical disruptions. For example, the Russia--Ukraine conflict led to significant volatility in global energy and electricity markets, affecting fuel prices, wholesale electricity prices, and associated adjustment mechanisms \citep{roeger2022gas, inacio2023assessing}. Incorporating such exogenous effects would substantially increase model complexity and reduce generalizability. Therefore, we use published ToU electricity prices directly and do not explicitly model unpredictable market shocks or policy-dependent price adjustments.
}

\begin{table}[H]
    \centering
    \caption{ToU electricity pricing (Sources: \citep{can_tou, aus_tou})}
    \begin{tabular}{cm{4cm}p{3.7cm}p{3.7cm}p{3cm}}
    \hline
        \textbf{Network} & \textbf{Period} & \textbf{Peak Hours} & \textbf{Mid Peak/Shoulder Hours} & \textbf{Off Peak/Overnight Hours}  \\
        \hline
        \multirow{3}{4em}{Durham} &
         May--October & 11 am--5 pm & 7 am--11 am, 5 pm--7 pm & 7 pm--7 am\\
        &November--April & 7 am--11 am, 5 pm--7 pm &  11 am--5 pm & 7 pm--7 am\\
            &Rate (\$/kWh) (May--Oct) & 0.1219 & 0.0817 & 0.0583 \\
            &Rate (\$/kWh) (Nov--Apr) & 0.1059 & 0.0817 & 0.0509 \\
            \hline
            \multirow{2}{4em}{Canberra} &  All year & 7 am--5 pm &5 pm--10 pm& 10 pm--7 am\\
            & Rate (\$/kWh) & 0.2439 & 0.1993 & 0.1647 \\
            \hline
    \end{tabular}
    \label{tab:ToU_pricing}
\end{table}

 \begin{figure}[H]
    \centering
    \begin{subfigure}{0.47\textwidth}
        \centering
        \includegraphics[scale=0.308]{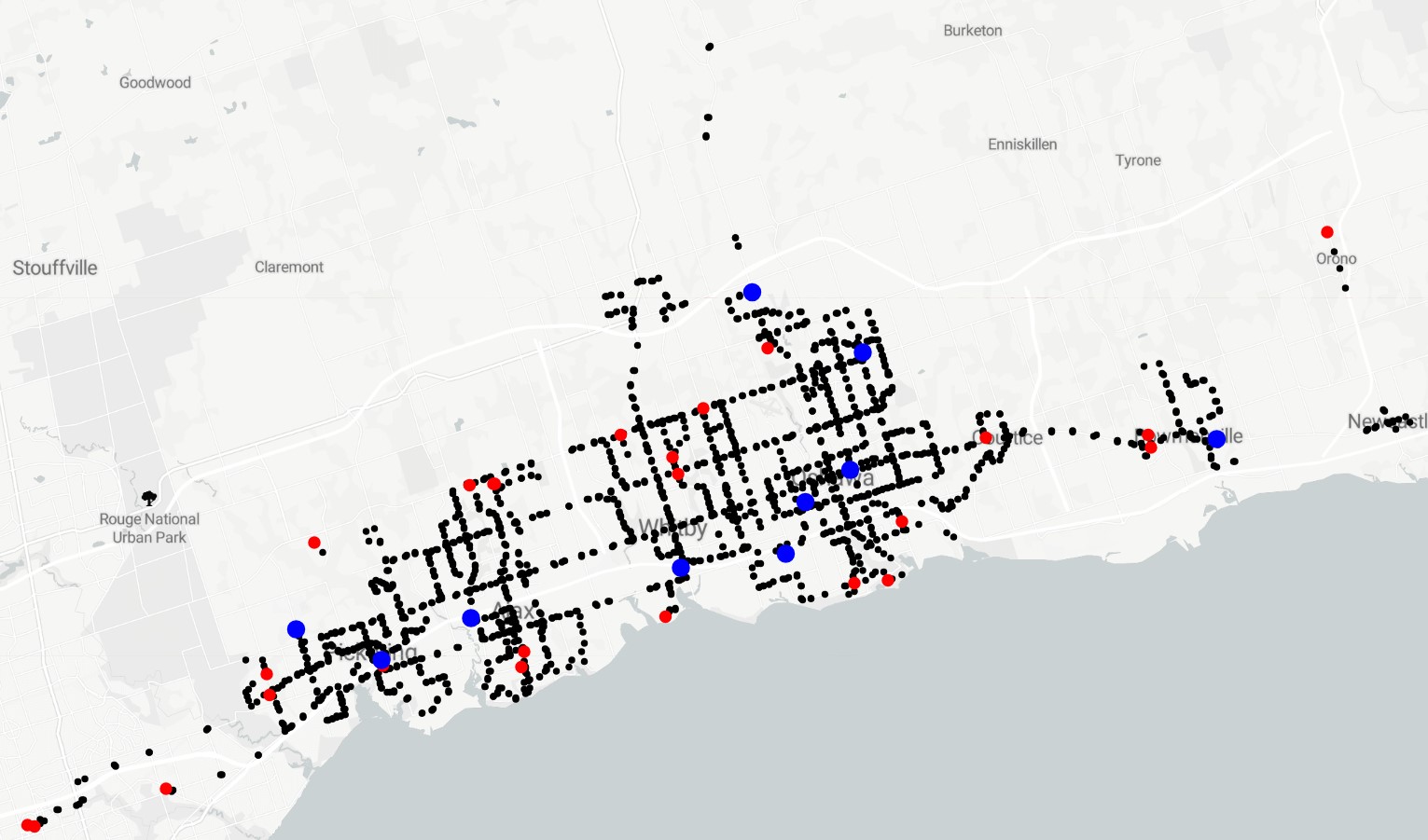}
        \caption{Durham transit network}
        \label{fig:Durham network}
    \end{subfigure}
    \hfill
    \begin{subfigure}{0.51\textwidth}
        \centering
        \includegraphics[scale=0.3]{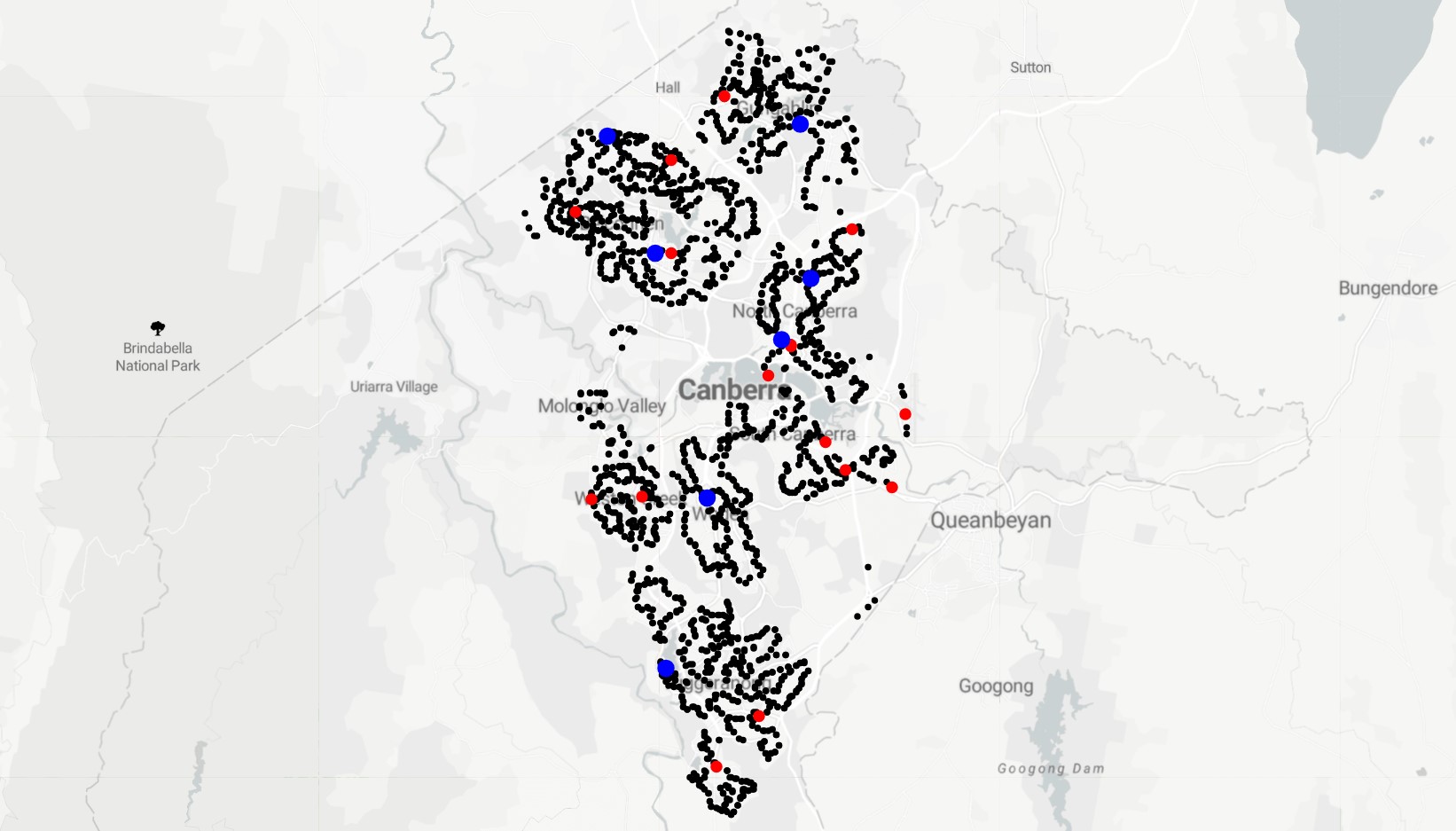}
        \caption{Canberra transit network}
        \label{fig:Canberra Network}
    \end{subfigure}
    \caption{\textcolor{addedcolor}{Spatial distribution of transit network components showing stops (black), depots (blue), and opportunity charging locations (red), where depots also function as charging stations in the model}}
        \label{fig:network_visualize}
\end{figure}

\subsection{Case Study}
Our study considered 1, 4, 12, and 52 scenarios with GTI and ambient temperature data averaged across yearly, quarterly, monthly, and weekly time scales, respectively. Temperature and GTI values were available at an hourly granularity; therefore, average values were computed for each hour across the days in the scenario, and the same hourly value was assigned to every minute within that hour. \textcolor{addedcolor}{Alternatively, in the absence of hourly data, solar energy generation can be estimated using the peak GTI value and effective solar hours \citep{footprinthero_peak_sun_hours}.} We also analyzed the benefits of introducing renewable energy sources and studied the effect of temperature-specific energy consumption estimates on fixed and operating costs. While implementing the CSP solution, if a transit operator encounters a week with weather conditions that deviate from those in the model, they can either reoptimize the CSP based on short-term weather forecasts or apply the solution from a closely matching scenario. \textcolor{addedcolor}{We also discuss the results of two worst-case scenarios based on the lowest solar energy generation and the highest energy consumption in Section \ref{sec:worst_case_analysis}.}
\vspace{-2mm}
\subsubsection{Network Details}
Information regarding stops, terminals, depots, routes, trips, and charging locations (obtained from Algorithm \ref{alg:cs}) for Durham and Canberra networks is given in Table \ref{tab:networl_details}. Figure \ref{fig:network_visualize} shows the charging location information for these networks.

\begin{table}[H]
    \centering
    \caption{Network details}
    \begin{tabular}{ccccccc}
    \hline
           \textbf{Network} & \textbf{Stops} & \textbf{Terminal Stops} & \textbf{Depots Used} & \textbf{Charging Locations} & \textbf{Routes} & \textbf{Trips}  \\
       \hline
       Durham  & 1890 &52 &10 & 36 & 32 & 2092 \\
       Canberra & 2423 & 67 & 7 & 24 & 65 & 3911 \\
       \hline
    \end{tabular}
    \label{tab:networl_details}
\end{table}
 
\begin{figure}[H]
  \centering
  \begin{subfigure}[b]{0.46\linewidth}
    \centering\includegraphics[scale=0.205]{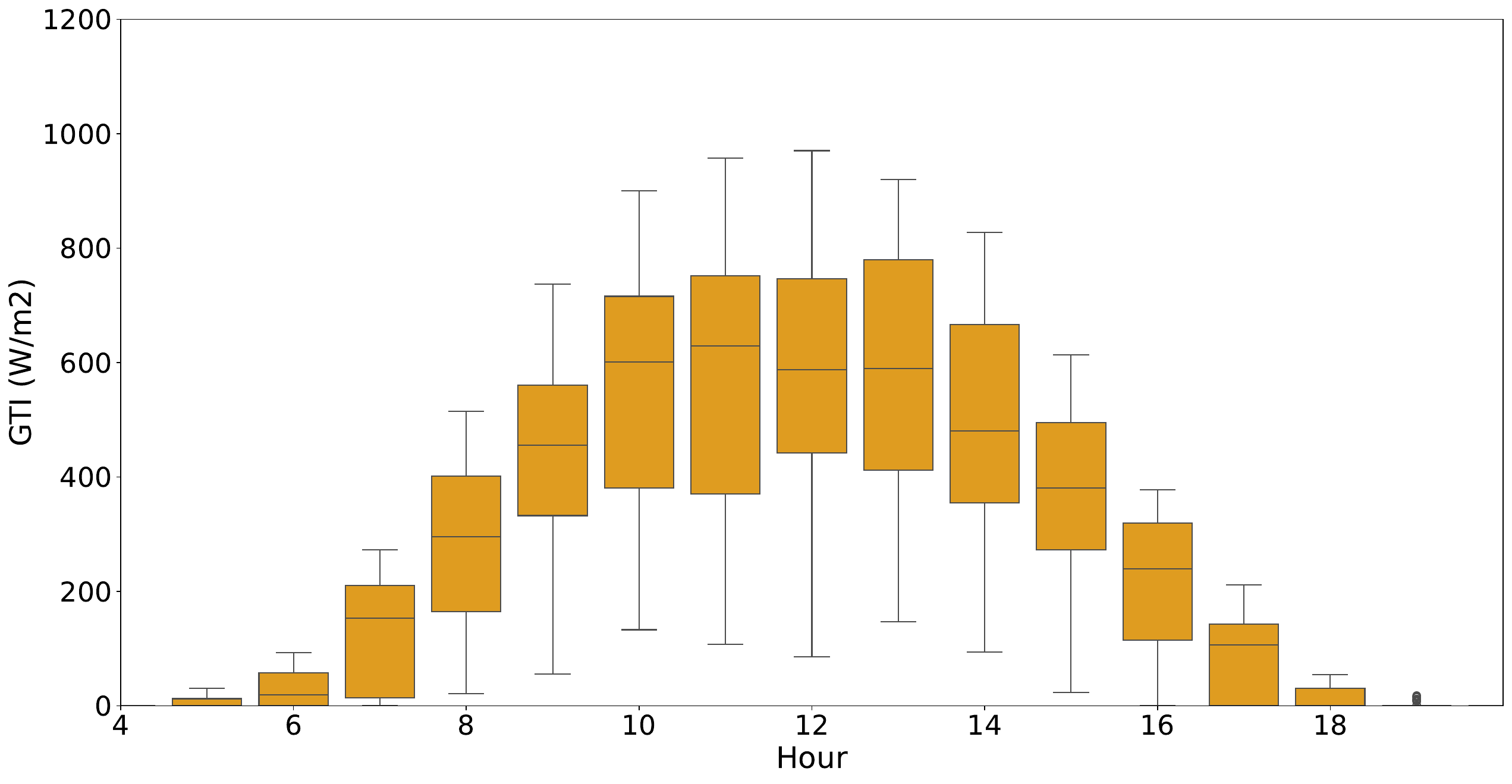}
    \caption{Durham}
    \label{fig:gti_box_durham}
  \end{subfigure}
\hfill
  \begin{subfigure}[b]{0.5\linewidth}
    \centering\includegraphics[scale=0.205]{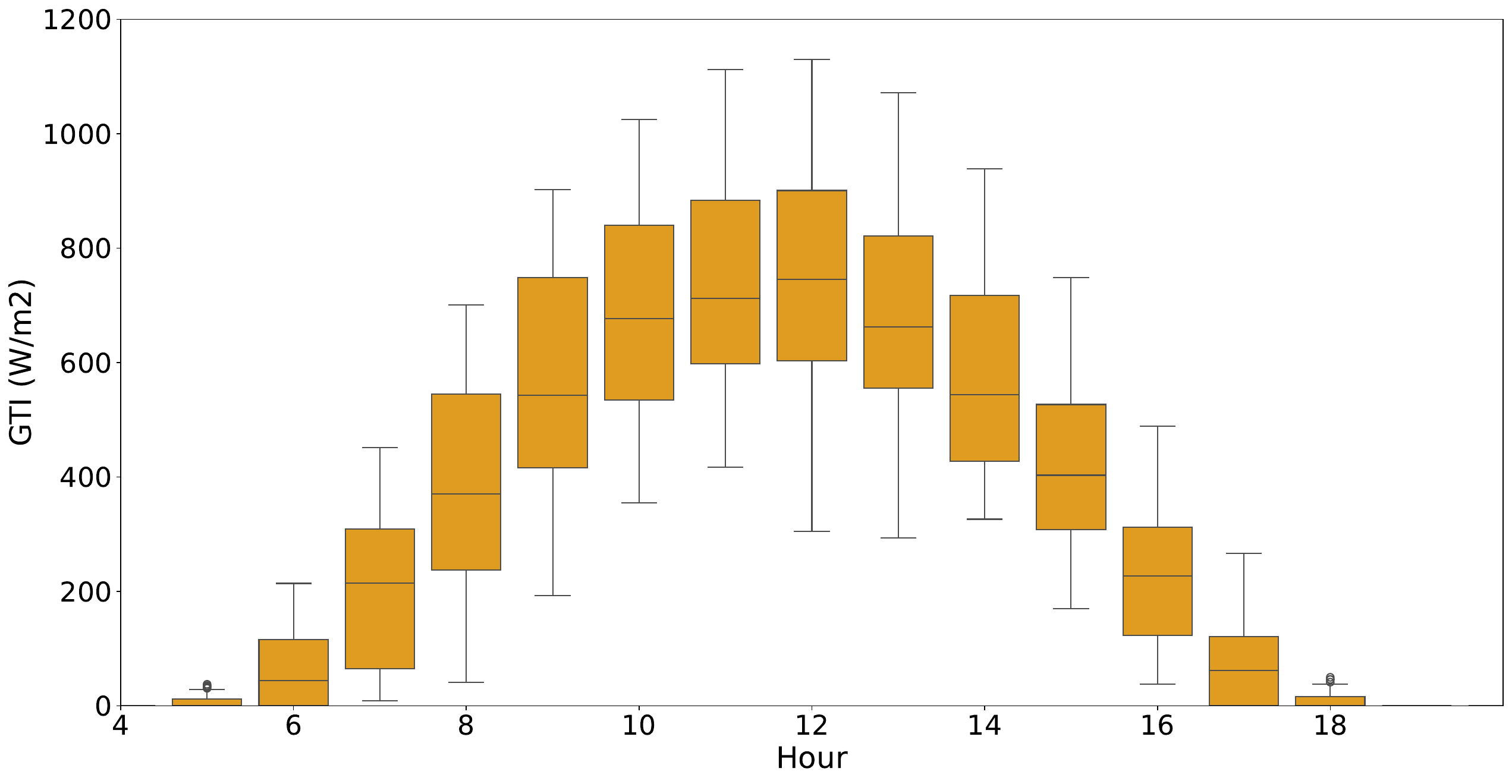}
    \caption{Canberra}
    \label{fig:gti_box_canberra}
  \end{subfigure}
  \caption{\textcolor{addedcolor}{Hourly box plot of Global Tilted Irradiance (GTI) showing intra-day variability and seasonal dispersion patterns in 2023}}
\label{fig:GTI__box_variations}
\end{figure}
 
\begin{figure}[t]
    \centering
    \begin{subfigure}{0.32\textwidth}
        \centering
        \includegraphics[scale=0.265]{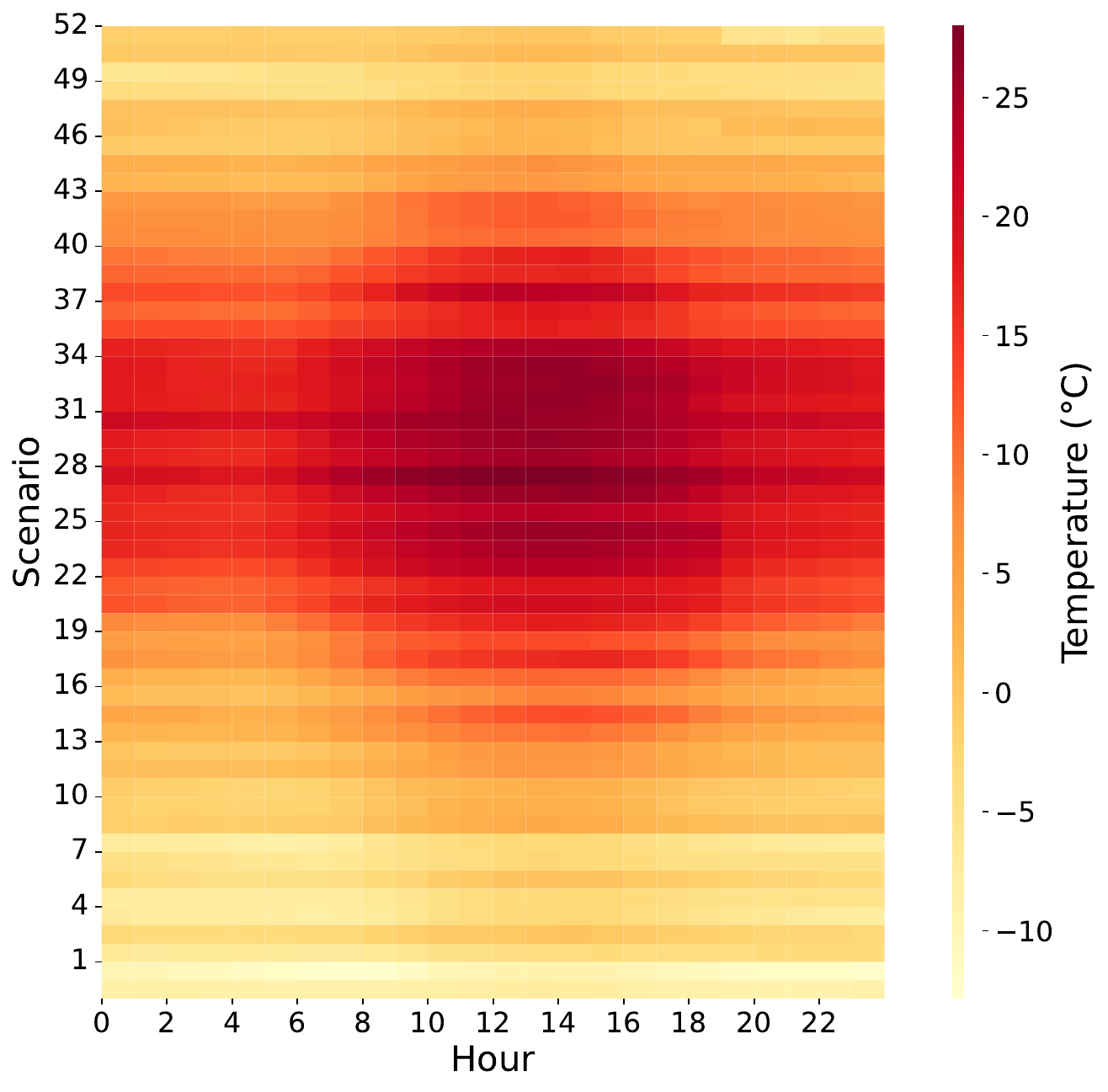}
        \caption{Durham (Temperature variation)}
        \label{fig:temperature_variation_durham}
    \end{subfigure}
    \hfill
    \begin{subfigure}{0.34\textwidth}
        \centering
        \includegraphics[scale=0.265]{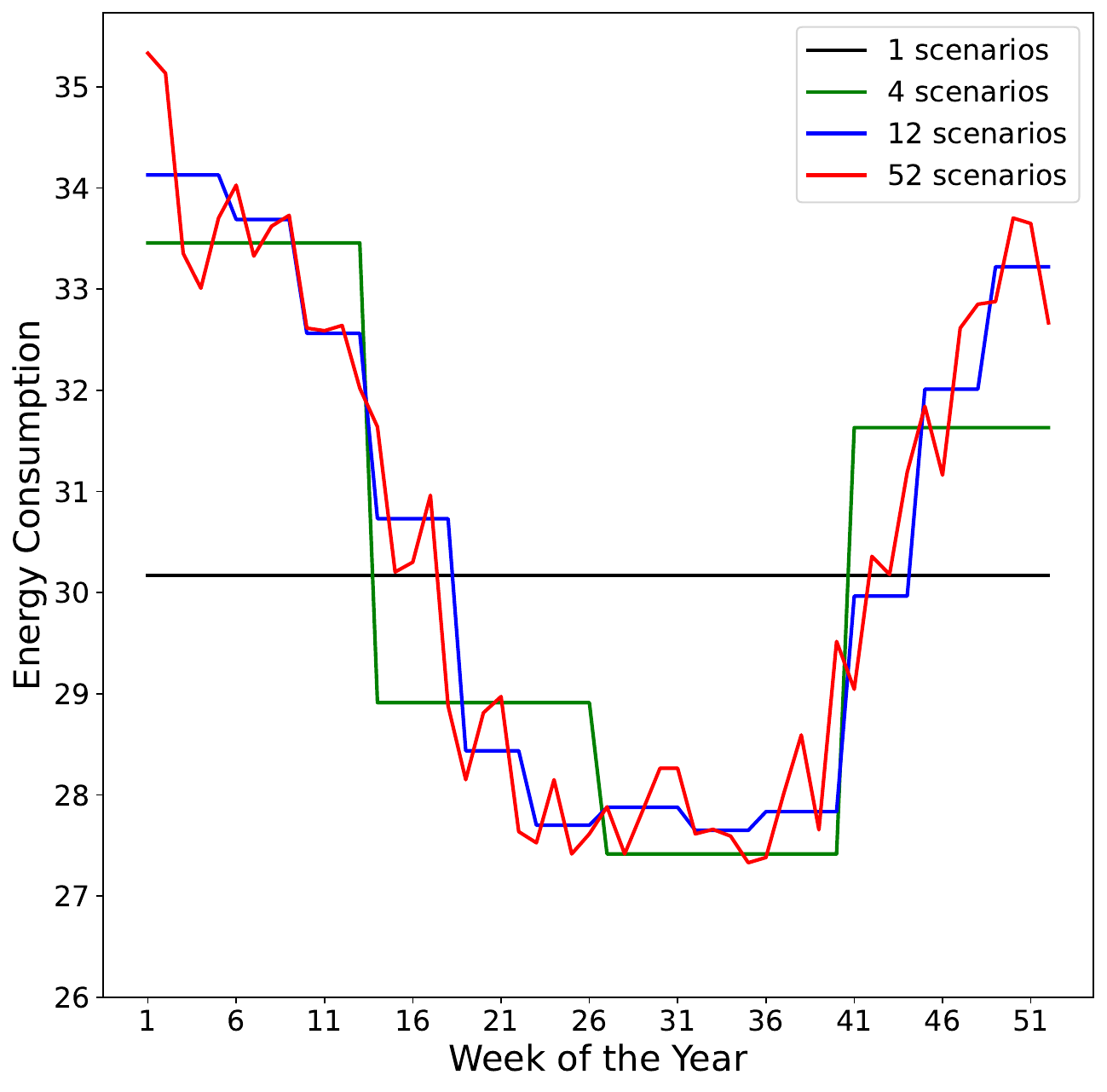}
        \caption{Durham (Energy required for trip 181)}
        \label{fig:energy_requirement_durham}
    \end{subfigure}
    \hfill
    \begin{subfigure}{0.32\textwidth}
        \centering
        \includegraphics[scale=0.265]{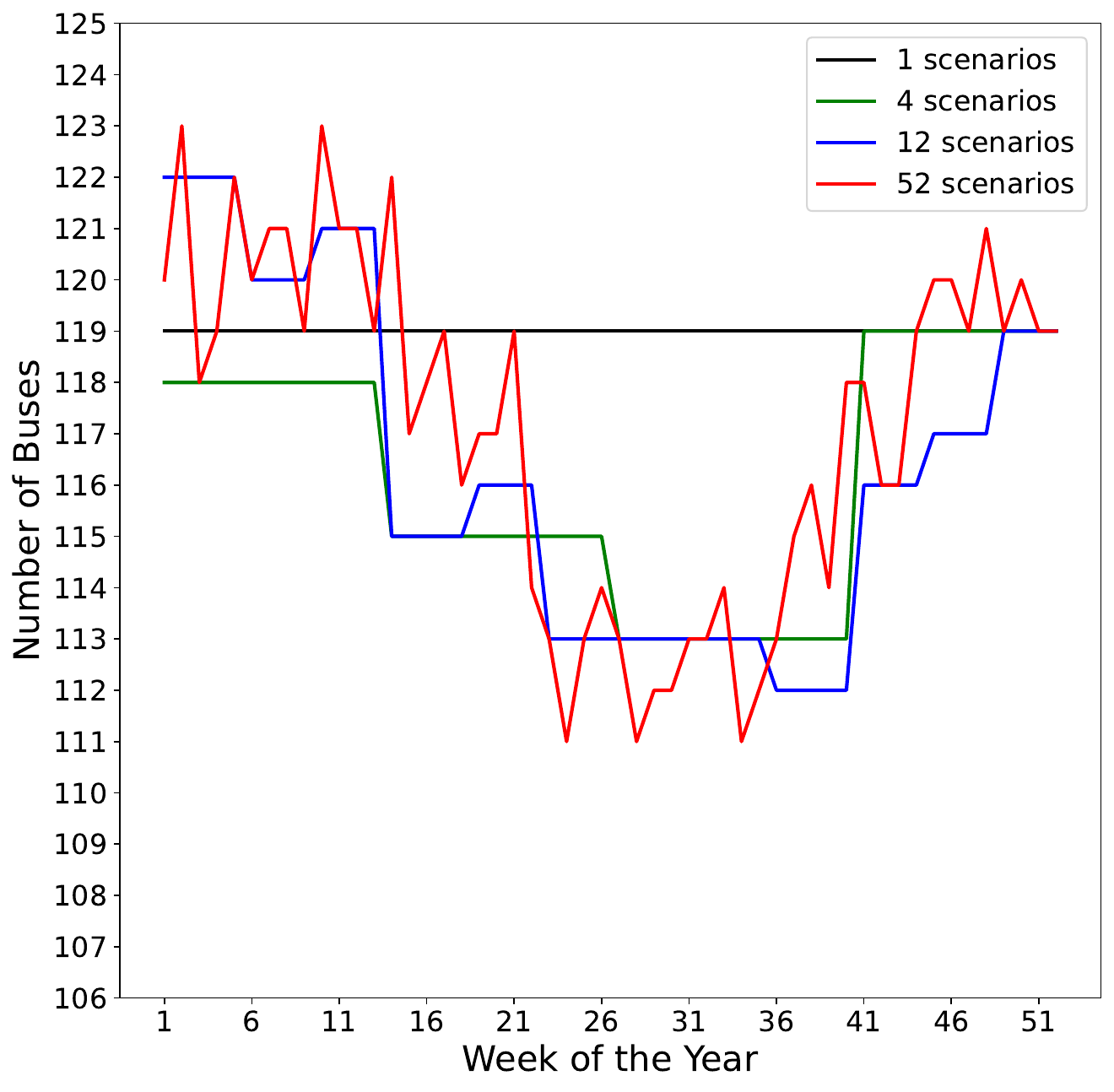}
        \caption{Durham (Number of buses)}
        \label{fig:bus_requirement_durham}
    \end{subfigure}
    \hfill
    \begin{subfigure}{0.32\textwidth}
        \centering
        \includegraphics[scale=0.265]{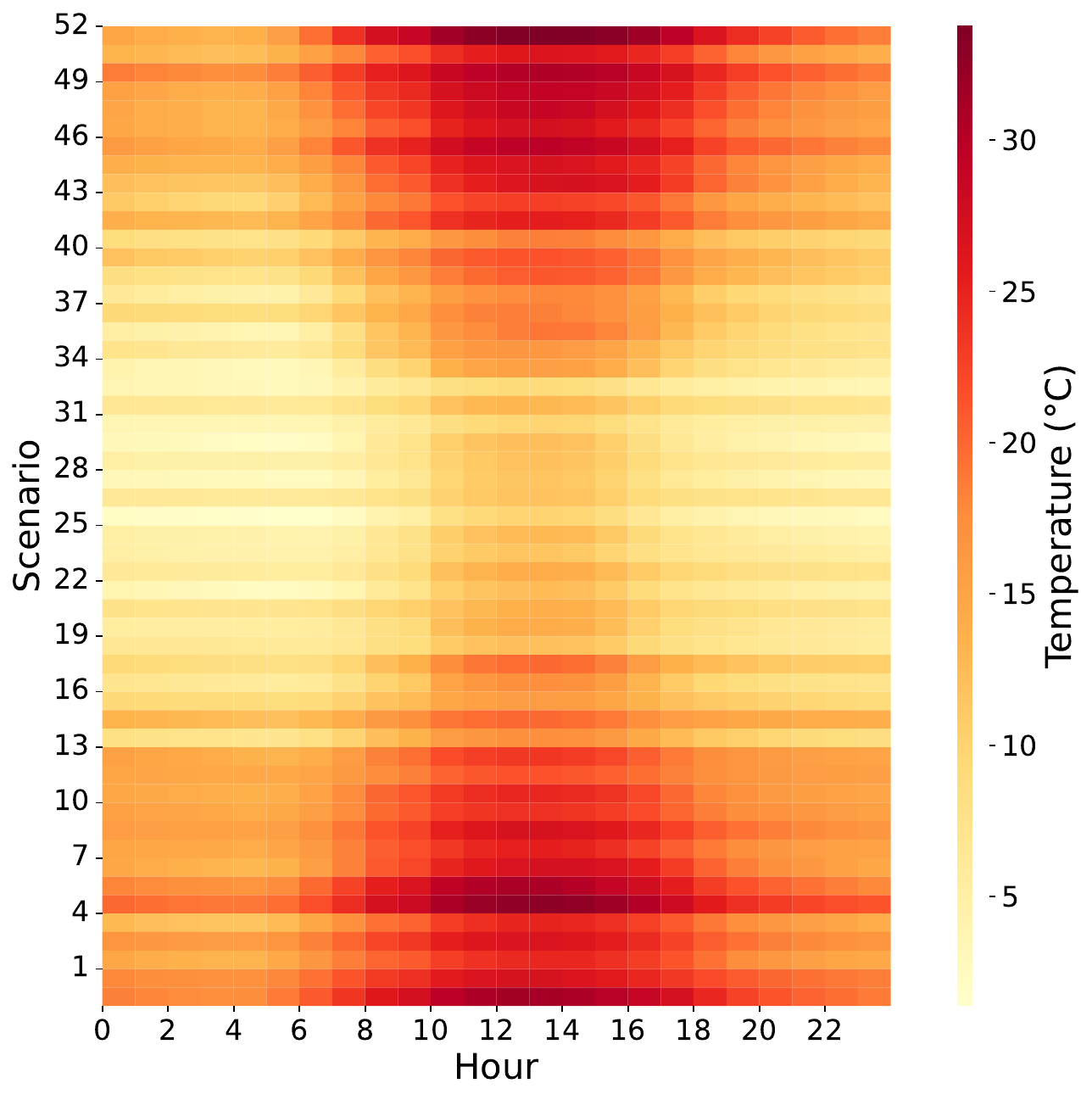}
        \caption{Canberra (Temperature variation)}
        \label{fig:temperature_variation_can}
    \end{subfigure}
    \hfill
    \begin{subfigure}{0.34\textwidth}
        \centering
        \includegraphics[scale=0.265]{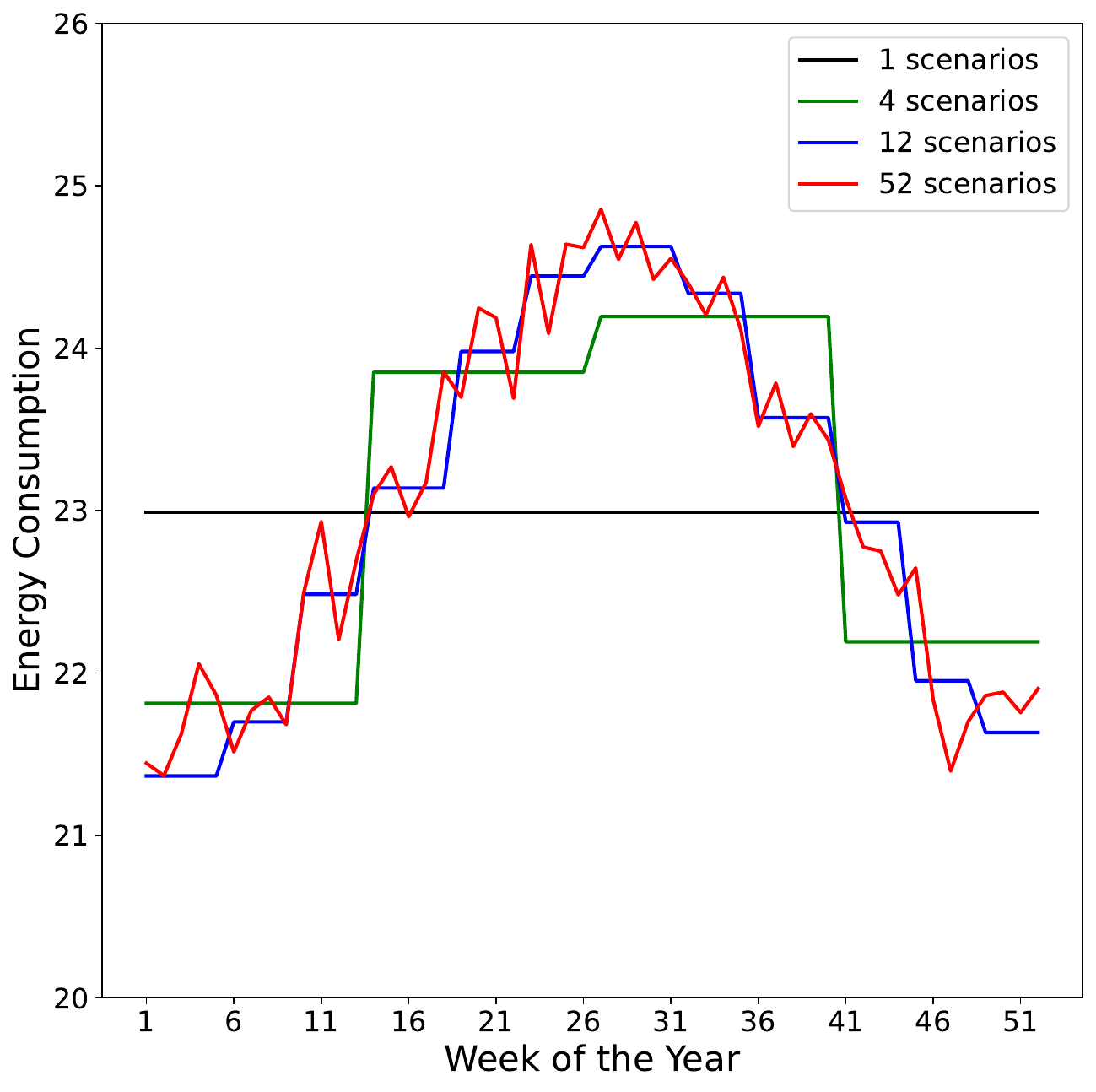}
        \caption{Canberra (Energy required for trip 307)}
        \label{fig:energy_requirement_can}
    \end{subfigure}
    \hfill
    \begin{subfigure}{0.32\textwidth}
        \centering
        \includegraphics[scale=0.265]{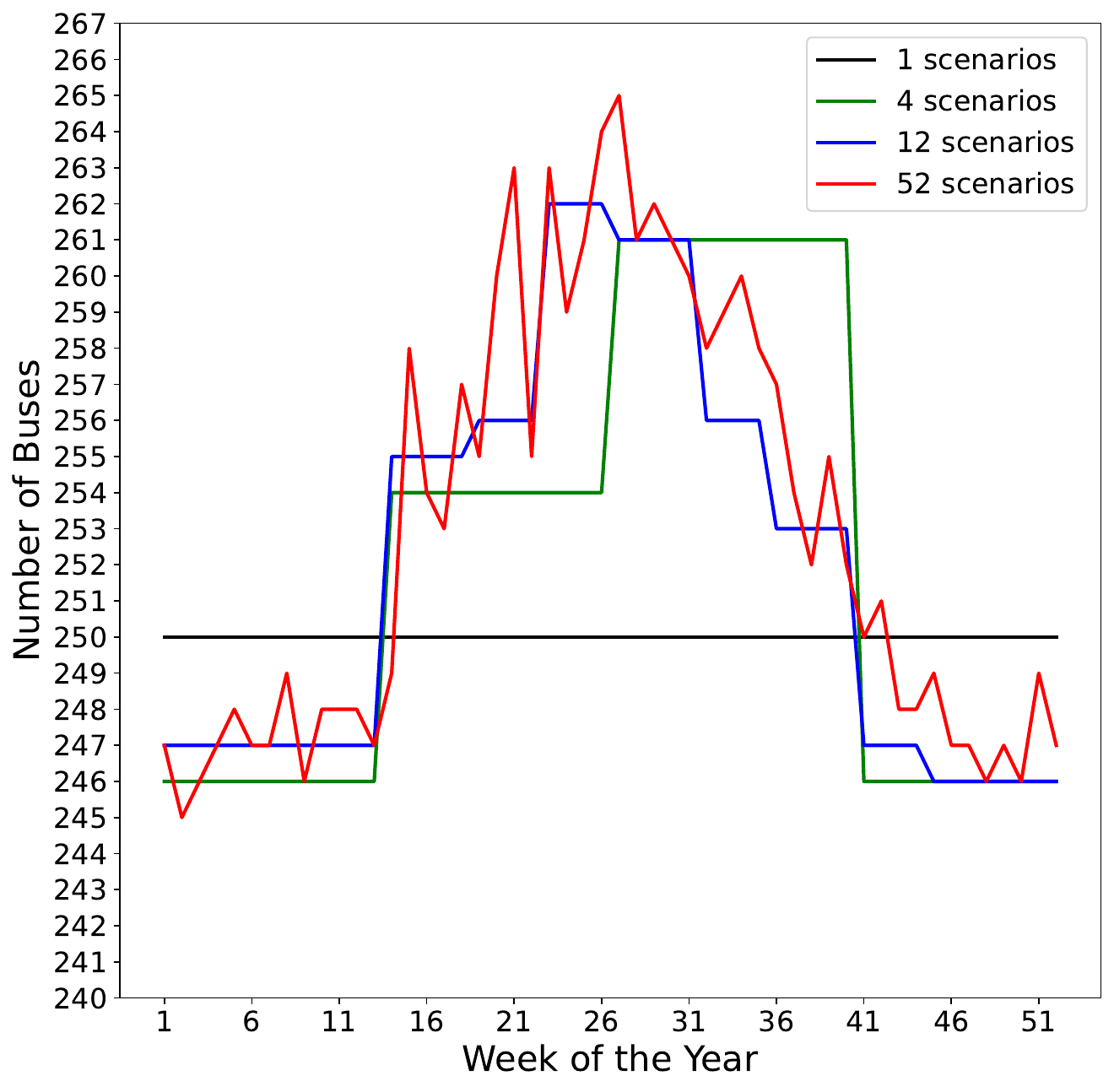}
        \caption{Canberra (Number of buses)}
        \label{fig:bus_requirement_can}
    \end{subfigure}
    \caption{Hourly temperature variation (left), trip energy requirement variations (middle), and total \textcolor{addedcolor}{number of} buses required (right) across different scenarios}
    \label{fig:temperature_variations}
\end{figure}

Variations in GTI values across scenarios and time-of-day for both networks are illustrated in Figures \ref{fig:GTI_variations} and \ref{fig:GTI__box_variations}, respectively. Note that seasonal changes (refer Figure \ref{fig:GTI_variations}) in the southern hemisphere are opposite to those in the northern hemisphere. The grid heat map (Figure \ref{fig:GTI_variations}) highlights variations across weeks/scenarios; box plots (Figure \ref{fig:GTI__box_variations}) indicate the quartile range variations \textcolor{addedcolor}{for each hour in a day}. Figure \ref{fig:temperature_variations} highlights the temperature variations in the two cities. Specifically, grid heat maps in Figures \ref{fig:temperature_variation_durham} and \ref{fig:temperature_variation_can} show the variations in temperature across scenarios (ranging from $-10$ to $25$ $^{\circ}$C for Durham and from $3$ to $33$ $^{\circ}$C for Canberra), which in turn result in different energy requirements for the same trip depending on the scenario. 

Figures \ref{fig:energy_requirement_durham} and \ref{fig:energy_requirement_can} show the variations in trip energy requirement when scenarios are built by averaging the temperature data over weekly, monthly, quarterly, or yearly time intervals. The energy needs for the longest trip can vary from 27 to 36 kWh for Durham (trip 181) and 21 to 25 kWh for Canberra (trip 307). \textcolor{addedcolor}{This affects both fleet size requirements (see Figures \ref{fig:bus_requirement_durham} and \ref{fig:bus_requirement_can}) and trip assignments across scenarios.} The number of required EBs was found to vary from 111 to 123 for Durham and 245 to 265 for Canberra. 

\subsubsection{Scenario-Based Results}
Table \ref{tab:run_time} shows the optimum objective values and computation times from our experiments. For small instances, \textcolor{addedcolor}{the} dual simplex \textcolor{addedcolor}{method} is quick. However, as the number of scenarios increases, Benders' decomposition solves the problem faster. Dual simplex could not find a feasible solution for Durham and Canberra networks with 52 scenarios in \textcolor{addedcolor}{36 hours}. The performance of Benders' decomposition depends on the tolerances of the optimality and feasibility cuts. \textcolor{addedcolor}{We use the default tolerances ($10^{-6}$) in CPLEX for all instances.} \textcolor{addedcolor}{The solver's implementation internally incorporates parallelization and efficient cut management. All reported solutions satisfy the solver’s global optimality gap; the Benders tolerance was used only for cut management.} For both networks, the objective values increase with the number of scenarios, \textcolor{addedcolor}{except for the 12-scenario case in Durham}, highlighting cost underestimation if sufficient scenarios are not considered. One of the reasons for this underestimation is that the model misses out on critical scenarios that influence long-term decision-making. 
\begin{table}[H]
    \centering
    \caption{\textcolor{addedcolor}{Summary of computational performance and cost breakdown for Durham and Canberra networks across scenarios. Bold values indicate the highest observed cost}}
    \begin{tabular}{m{4.5cm}rrrr|rrrr}
    \hline
    \textbf{Network} & \multicolumn{4}{c}{\textbf{Durham}} & \multicolumn{4}{c}{\textbf{Canberra}}\\
    \hline
        Scenarios & \textcolor{addedcolor}{1} & \textcolor{addedcolor}{4} & \textcolor{addedcolor}{12} & \textcolor{addedcolor}{52} & \textcolor{addedcolor}{1} & \textcolor{addedcolor}{4} & \textcolor{addedcolor}{12} & \textcolor{addedcolor}{52}\\
        \hline
        Variables (millions) & \textcolor{addedcolor}{0.27} & \textcolor{addedcolor}{1.07} & \textcolor{addedcolor}{3.21} & \textcolor{addedcolor}{13.99} & \textcolor{addedcolor}{0.41} & \textcolor{addedcolor}{1.66} & \textcolor{addedcolor}{4.99} & \textcolor{addedcolor}{21.78}\\
        Constraints (millions) & \textcolor{addedcolor}{0.29} & \textcolor{addedcolor}{1.17} & \textcolor{addedcolor}{3.50} & \textcolor{addedcolor}{15.22} & \textcolor{addedcolor}{0.31} & \textcolor{addedcolor}{1.27} & \textcolor{addedcolor}{3.81} & \textcolor{addedcolor}{16.57}\\
        Benders' runtime (minutes) & \textcolor{addedcolor}{73.13} & \textcolor{addedcolor}{112.73} & \textcolor{addedcolor}{173.29} & \textcolor{addedcolor}{360.11} & \textcolor{addedcolor}{71.22} & \textcolor{addedcolor}{382.54} & \textcolor{addedcolor}{1209.70} & \textcolor{addedcolor}{1446.50}\\
        Simplex runtime (minutes) & \textcolor{addedcolor}{2.55} & \textcolor{addedcolor}{18.74} & \textcolor{addedcolor}{261.96} & \textcolor{addedcolor}{2160} & \textcolor{addedcolor}{0.56} & \textcolor{addedcolor}{119.38} & \textcolor{addedcolor}{898.79} & \textcolor{addedcolor}{2160}\\
        Benders' objective (\$) & \textcolor{addedcolor}{1179.52} & \textcolor{addedcolor}{1363.54} & \textcolor{addedcolor}{1349.86} & \textcolor{addedcolor}{1406.48} & \textcolor{addedcolor}{5219.19} & \textcolor{addedcolor}{6585.90} & \textcolor{addedcolor}{6666.52} & \textcolor{addedcolor}{7190.09}\\
        Simplex objective (\$) & \textcolor{addedcolor}{1179.52} & \textcolor{addedcolor}{1363.54} & \textcolor{addedcolor}{1349.86} & \textcolor{addedcolor}{--} & \textcolor{addedcolor}{5219.19} & \textcolor{addedcolor}{6585.93} & \textcolor{addedcolor}{6666.52} & \textcolor{addedcolor}{--}\\
        \hline
        Contracted capacity cost (\$) & \textcolor{addedcolor}{143.57} & \textcolor{addedcolor}{\textbf{246.84}} & \textcolor{addedcolor}{206.16} & \textcolor{addedcolor}{230.95} & \textcolor{addedcolor}{0} & \textcolor{addedcolor}{446.69} & \textcolor{addedcolor}{526.65} & \textcolor{addedcolor}{\textbf{607.46}} \\
        Solar panel cost (\$) & \textcolor{addedcolor}{\textbf{456.15}} & \textcolor{addedcolor}{343.06} & \textcolor{addedcolor}{332.82} & \textcolor{addedcolor}{268.17} & \textcolor{addedcolor}{\textbf{1636.23}} & \textcolor{addedcolor}{746.12} & \textcolor{addedcolor}{674.99} & \textcolor{addedcolor}{511.13}\\
        BESS cost (\$) & \textcolor{addedcolor}{\textbf{46.21}} & \textcolor{addedcolor}{38.24} & \textcolor{addedcolor}{41.16} & \textcolor{addedcolor}{45.96} & \textcolor{addedcolor}{\textbf{3582.96}} & \textcolor{addedcolor}{500.31} & \textcolor{addedcolor}{401.90} & \textcolor{addedcolor}{152.66}\\
        Average operation cost (\$) & \textcolor{addedcolor}{533.59} & \textcolor{addedcolor}{735.40} & \textcolor{addedcolor}{769.72} & \textcolor{addedcolor}{\textbf{861.40}} & \textcolor{addedcolor}{0} & \textcolor{addedcolor}{5339.47} & \textcolor{addedcolor}{5062.98} & \textcolor{addedcolor}{\textbf{5918.84}} \\
        \hline
    \end{tabular}
    \label{tab:run_time}
\end{table}

When the number of scenarios is small, the installed solar panels dominate the charging station power requirements, and hence, the contracted capacity costs are low. Conversely, as the number of scenarios increases with a shorter time interval for each scenario, the grid power capacity requirement increases to handle the extreme cases of low solar energy generation, reducing the requirement for solar panels. Note that for the 4-scenario case in Canberra, the Benders' objective value is slightly \textcolor{addedcolor}{lower} than the dual simplex objective. This happens due to the tolerances used in the optimization solver.

\textcolor{addedcolor}{To further assess the computational behavior of Benders’ decomposition, we report convergence statistics for the 52-scenario case. The convergence behavior of the proposed algorithm is illustrated in Figure \ref{fig:gap_iterations} for both Durham and Canberra networks. In both cases, the optimality gap exhibits an overall decreasing trend as additional cuts are generated, eventually converging to the optimal solution. In the Durham network, the gap reduces rapidly in the initial iterations, then gradually tapers as the algorithm approaches optimality, a characteristic of decomposition-based methods. A similar trend is observed for the Canberra network, although convergence is comparatively slower due to the larger problem size, requiring more iterations. Note that the optimality gap increases slightly in a few iterations for both networks due to non-monotonic updates to the dual bounds, potentially caused by perturbations introduced by the solver to handle numerical issues.}

\textcolor{addedcolor}{In terms of cut generation, we observed that the algorithm relies predominantly on feasibility cuts, which constitute approximately 86.45\% and 80.21\% of the total cuts for the Durham and Canberra networks, respectively, according to the solver logs. Multiple feasibility/optimality cuts can be added in the same iteration. We noticed that feasibility plays a critical role in improving the dual bound, especially in the middle stages of the algorithm, whereas optimality cuts contribute more significantly during the initial and later stages.}

\begin{figure}[H]
    \centering
    \begin{subfigure}[b]{0.48\textwidth}
        \centering
        \includegraphics[width=\textwidth]{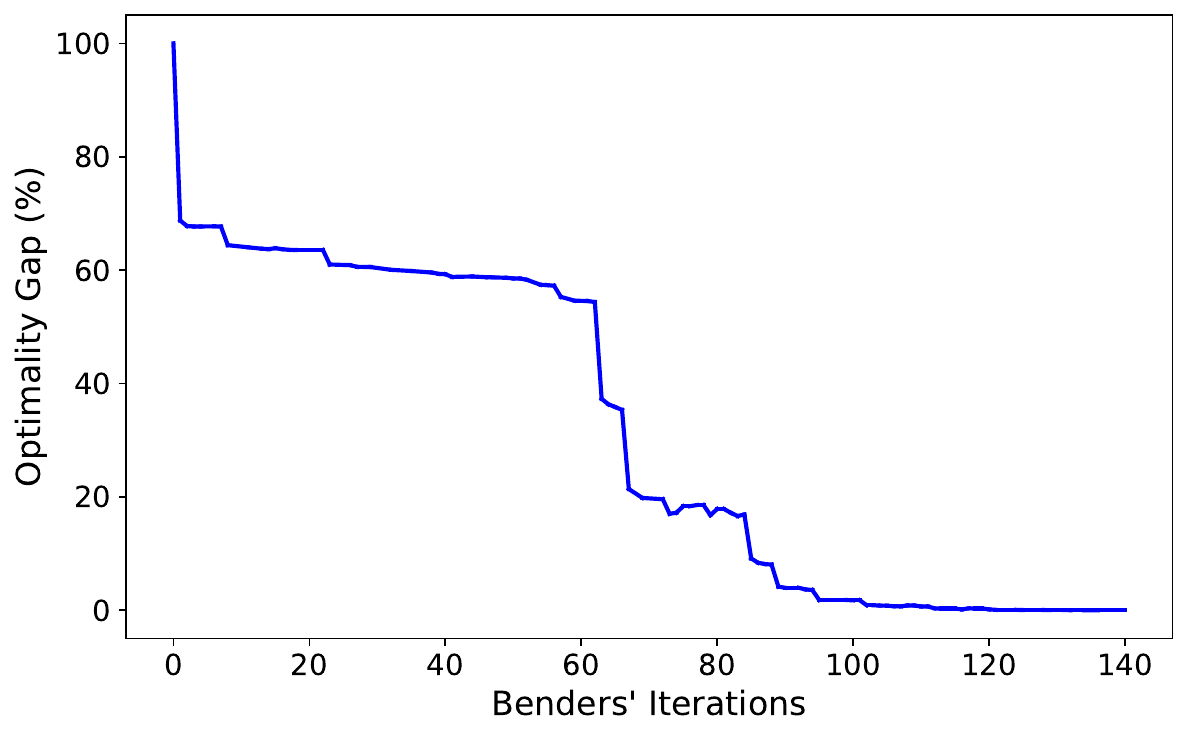}
        \caption{Durham network}
        \label{fig:durham_gap_iter}
    \end{subfigure}
    \hfill
    \begin{subfigure}[b]{0.48\textwidth}
        \centering
        \includegraphics[width=\textwidth]{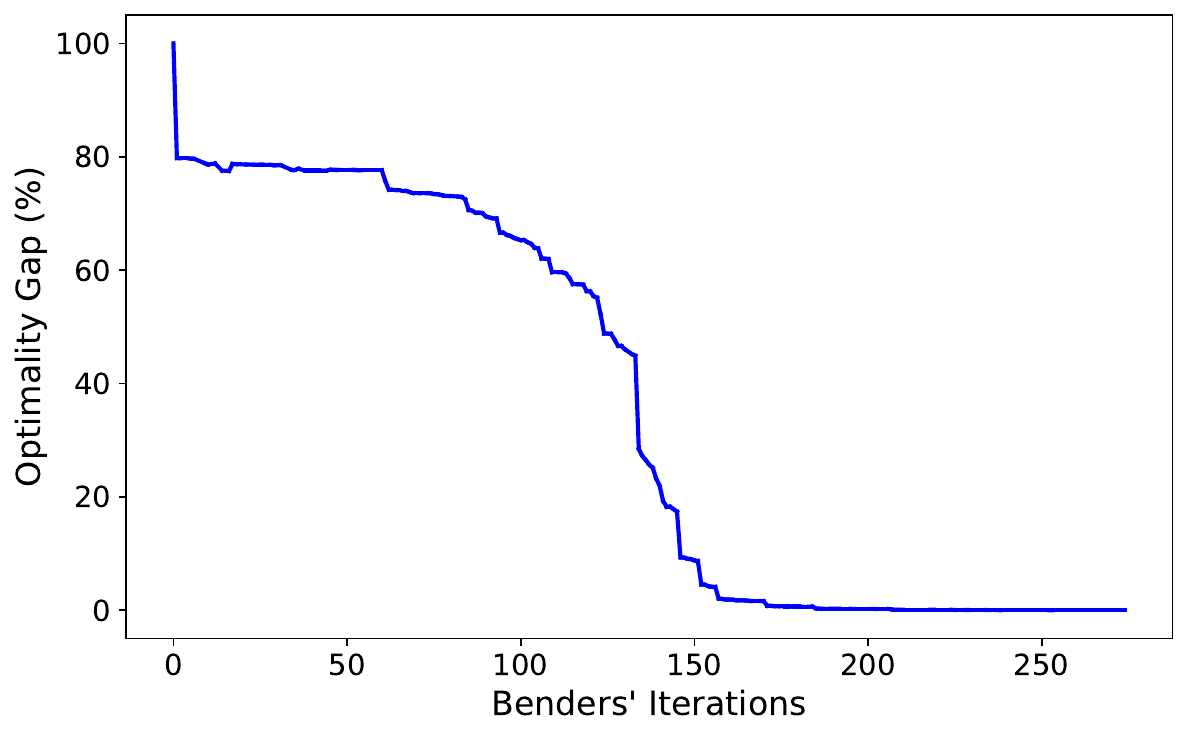}
        \caption{Canberra network}
        \label{fig:canberra_gap_iter}
    \end{subfigure}
    \caption{\textcolor{addedcolor}{Optimality gap (\%) versus the number of iterations for the two networks}}
    \label{fig:gap_iterations}
\end{figure}

\subsubsection{Effect of Introducing RES and Temperature-Based Trip Energy Estimates}

We also explored a 52-week scenario model with/without RES and with/without \textcolor{addedcolor}{temperature-based} variations in energy estimates. Table \ref{tab:renewables_and_temperature_results} summarizes the results. In both networks, considering renewables and BESS prove to be more cost-effective, as seen from the third and fourth columns of the table. Without RES, all energy is sourced from the grid, which is expensive. Note that the total energy requirements vary slightly between the two cases. This variation is due to the lower solar energy costs and deadheading of buses to charging stations within a cluster of terminals. When RES is available, buses may deadhead to a charging location that uses PVs and BESS, thereby avoiding direct charging from the grid at a higher electricity price. Due to RES, significant cost savings of \textcolor{addedcolor}{$9.72\%$ and $23.79\%$} are achieved for the Durham and Canberra networks, respectively. 

The cost savings for the Canberra network are much higher than those for the Durham network due to Canberra's greater potential for solar energy generation (see Figure \ref{fig:GTI_variations}), and its higher electricity prices (see Table \ref{tab:ToU_pricing}). The mean GTI across all locations for Durham and Canberra is $184.82$ W/\si{m^{2}} and $226.44$ W/\si{m^{2}}, respectively, indicating greater solar energy \textcolor{addedcolor}{availability} in Canberra. Without RES, the average operational costs constitute \textcolor{addedcolor}{$78.96\%$ and $89.70\%$} of the total costs for Durham and Canberra, respectively, suggesting a higher cost-saving potential for Canberra after integrating RES.

We also evaluated the optimal costs for a model with RES, ignoring the impacts of temperature on the trip energy estimates, and found that the overall cost is underestimated. This difference can be attributed to the temperature component in \eqref{trip_energy_consumption}, which is set to the optimal operating temperature in the no-temperature variation case, leading to an underestimation of the total energy requirements and, consequently, a reduction in operational costs. When temperature variations are neglected, we underestimate the total costs by \textcolor{addedcolor}{$24.71\%$ and $13.34\%$} for the two networks. The cost deviation is higher in the Durham network because the difference between the optimum working temperature and the ambient air temperature is greater in Durham than in Canberra.

\begin{table}[H]
    \centering
    \caption{\textcolor{addedcolor}{Comparison of system costs under RES and temperature-aware vs. temperature-agnostic trip energy modeling. Bold values indicate total daily costs, showing cost savings from RES and underestimation when temperature effects are ignored}}
    \begin{tabular}{c m{5.3cm} r r | r}
    \toprule
    \textbf{Network} & \textbf{Metrics} & \multicolumn{2}{c}{\textbf{With temperature effects}} & \multicolumn{1}{c}{\textbf{Ignoring temperature}} \\
    \cmidrule{3-4}\cmidrule{5-5} 
    & & \textbf{With RES} & \multicolumn{1}{c}{\textbf{Without RES}} & \textbf{With RES}\\
    \midrule
    \multirow{6}{4em}{Durham} 
        & Objective (\$) & \textcolor{addedcolor}{\textbf{1406.48}} & \textcolor{addedcolor}{\textbf{1557.85}} & \textcolor{addedcolor}{\textbf{1127.8}} \\
        \cline{2-5}
        & Average operational cost (\$) & \textcolor{addedcolor}{861.4} & \textcolor{addedcolor}{1230.14} & \textcolor{addedcolor}{652.81} \\
        & Contracted capacity cost (\$) & \textcolor{addedcolor}{230.95} & \textcolor{addedcolor}{327.71} & \textcolor{addedcolor}{164.65} \\
        & Solar panel cost (\$) & \textcolor{addedcolor}{268.17} & \textcolor{addedcolor}{0} & \textcolor{addedcolor}{271.01} \\
        & BESS cost (\$) & \textcolor{addedcolor}{45.96} & \textcolor{addedcolor}{0} & \textcolor{addedcolor}{39.33} \\ 
        & Total energy requirements (MWh) & \textcolor{addedcolor}{1046.11} & \textcolor{addedcolor}{1043.88} & \textcolor{addedcolor}{872.07} \\
    \midrule
    \multirow{6}{4em}{Canberra} 
        & Objective (\$) & \textcolor{addedcolor}{\textbf{7190.09}} & \textcolor{addedcolor}{\textbf{9434.91}} & \textcolor{addedcolor}{\textbf{6343.99}} \\
        \cline{2-5}
        & Average operational cost (\$) & \textcolor{addedcolor}{5918.84} & \textcolor{addedcolor}{8463.24} & \textcolor{addedcolor}{5132.81} \\
        & Contracted capacity cost (\$) & \textcolor{addedcolor}{607.46} & \textcolor{addedcolor}{971.67} & \textcolor{addedcolor}{561.41} \\
        & Solar panel cost (\$) & \textcolor{addedcolor}{511.13} & \textcolor{addedcolor}{0} & \textcolor{addedcolor}{499.35} \\
        & BESS cost (\$) & \textcolor{addedcolor}{152.66} & \textcolor{addedcolor}{0} & \textcolor{addedcolor}{150.42} \\  
        & Total energy requirements (MWh) & \textcolor{addedcolor}{2424.82} & \textcolor{addedcolor}{2420.73} & \textcolor{addedcolor}{2200.77} \\
    \bottomrule
    \end{tabular}
    \label{tab:renewables_and_temperature_results}
\end{table}

\textcolor{addedcolor}{To assess sensitivity to the input parameters, we perturbed the temperature-dependent parameter ($\alpha_4$) that governs the energy consumption. Table \ref{tab:temp_sensitivity} reports the resulting objective values and their cost components for both Durham and Canberra.}

\begin{table}[H]
\centering
\caption{\textcolor{addedcolor}{Sensitivity analysis to show variation of the total cost and cost components with respect to the temperature coefficient $\alpha_4$ for the 52-scenario case}}
\label{tab:temp_sensitivity}
\begin{tabular}{lrrrrr}
\toprule
\multicolumn{6}{c}{\textcolor{addedcolor}{\textbf{Durham}}} \\ \midrule \textcolor{addedcolor}{\textbf{Temperature coefficient ($\alpha_4$)}} & \textcolor{addedcolor}{0.002} & \textcolor{addedcolor}{0.004} & \textcolor{addedcolor}{0.006} & \textcolor{addedcolor}{0.008} & \textcolor{addedcolor}{0.010} \\ \midrule \textcolor{addedcolor}{Objective (\$)} & \textcolor{addedcolor}{1256.09} & \textcolor{addedcolor}{1303.36} & \textcolor{addedcolor}{1353.45} & \textcolor{addedcolor}{1406.48} & \textcolor{addedcolor}{1462.64} \\ \textcolor{addedcolor}{Contracted capacity cost (\$)} & \textcolor{addedcolor}{192.64} & \textcolor{addedcolor}{204.44} & \textcolor{addedcolor}{217.91} & \textcolor{addedcolor}{230.95} & \textcolor{addedcolor}{249.41} \\ \textcolor{addedcolor}{Solar panel cost (\$)} & \textcolor{addedcolor}{264.29} & \textcolor{addedcolor}{265.65} & \textcolor{addedcolor}{267.67} & \textcolor{addedcolor}{268.17} & \textcolor{addedcolor}{268.91} \\ \textcolor{addedcolor}{BESS cost (\$)} & \textcolor{addedcolor}{45.05} & \textcolor{addedcolor}{45.40} & \textcolor{addedcolor}{46.17} & \textcolor{addedcolor}{45.96} & \textcolor{addedcolor}{45.87} \\ \textcolor{addedcolor}{Average operation cost (\$)} & \textcolor{addedcolor}{754.11} & \textcolor{addedcolor}{787.87} & \textcolor{addedcolor}{821.70} & \textcolor{addedcolor}{861.40} & \textcolor{addedcolor}{898.45} \\ \midrule \multicolumn{6}{c}{\textcolor{addedcolor}{\textbf{Canberra}}} \\ \midrule \textcolor{addedcolor}{\textbf{Temperature coefficient ($\alpha_4$)}} & \textcolor{addedcolor}{0.002} & \textcolor{addedcolor}{0.004} & \textcolor{addedcolor}{0.006} & \textcolor{addedcolor}{0.008} & \textcolor{addedcolor}{0.010} \\ \midrule \textcolor{addedcolor}{Objective (\$)} & \textcolor{addedcolor}{6686.88} & \textcolor{addedcolor}{6849.67} & \textcolor{addedcolor}{7018.18} & \textcolor{addedcolor}{7190.09} & \textcolor{addedcolor}{7368.04} \\ \textcolor{addedcolor}{Contracted capacity cost (\$)} & \textcolor{addedcolor}{583.74} & \textcolor{addedcolor}{589.57} & \textcolor{addedcolor}{597.37} & \textcolor{addedcolor}{607.46} & \textcolor{addedcolor}{618.65} \\ \textcolor{addedcolor}{Solar panel cost (\$)} & \textcolor{addedcolor}{506.01} & \textcolor{addedcolor}{507.76} & \textcolor{addedcolor}{509.97} & \textcolor{addedcolor}{511.13} & \textcolor{addedcolor}{513.15} \\ \textcolor{addedcolor}{BESS cost (\$)} & \textcolor{addedcolor}{148.35} & \textcolor{addedcolor}{149.95} & \textcolor{addedcolor}{154.84} & \textcolor{addedcolor}{152.66} & \textcolor{addedcolor}{155.62} \\ \textcolor{addedcolor}{Average operation cost (\$)} & \textcolor{addedcolor}{5448.78} & \textcolor{addedcolor}{5602.39} & \textcolor{addedcolor}{5756.00} & \textcolor{addedcolor}{5918.84} & \textcolor{addedcolor}{6080.62} \\
\bottomrule
\end{tabular}
\end{table}
\textcolor{addedcolor}{For both networks, the total cost increases monotonically with $\alpha_4$, reflecting higher energy consumption under stronger temperature effects. For Durham, the total cost varies from \$1256.09 to \$1462.64, corresponding to a deviation of approximately $-10.69\%$ to $+3.99\%$ relative to the base case (\$1406.48). Similarly, for Canberra, the total cost ranges from \$6686.88 to \$7368.04, corresponding to a deviation of approximately $-7.00\%$ to $+2.48\%$ relative to the base case (\$7190.09). Analysis of the cost components shows that most of this variation is driven by the average operational cost, which increases from $\$754.11$ to $\$898.45$ for Durham and from $\$5448.78$ to $\$6080.62$ for Canberra as $\alpha_4$ increases. In comparison, the solar PV and BESS investment costs exhibit only minor changes across the tested parameter range, as the contracted-capacity cost increases. These results highlight the importance of accurately estimating temperature-dependent energy consumption, as it can have a significant impact on average operation and contracted capacity cost estimations.}

\subsection{Extended Analysis}

\subsubsection{\textcolor{addedcolor}{Infrastructure Decisions}}

While underestimating grid capacity can lead to power shortage/failure during operations, overestimating solar panel area can lead to excessive power generation that would be wasted if BESS capacity is insufficient. Figures \ref{fig:master_variables_durham} and \ref{fig:master_variables_canberra} display charging location capacities, solar panel areas, and battery storage capacities necessary for optimal operations across different scenarios for Durham and Canberra networks, respectively. Since Algorithm \ref{alg:cs} is a greedy approach and the second-stage decisions are optimized in the CSP, many charging locations (\textcolor{addedcolor}{55.56\%} for Durham and 12.5\% for Canberra, \textcolor{addedcolor}{considering 1, 4, 12, and 52 scenario cases}) obtained through the CS-based heuristic are redundant. \textcolor{addedcolor}{The redundancy arises because the CS heuristic finds feasible solutions that are not cost-optimal candidates, which results in selecting more charging locations than are actually necessary. The CSP subsequently filters from this candidate set during optimization. Incorporating cost-aware pruning or tighter candidate generation within CS could further reduce this redundancy, but would have little effect overall since the CSP fixes the inefficiencies in the solutions.} \textcolor{addedcolor}{In Figures \ref{fig:master_variables_durham} and \ref{fig:master_variables_canberra}}, the charging location labels on the x-axis represent the terminal bus stop indices used for the grid- or solar-powered system. Some stops are not shown in both figures because buses do not charge at those locations.

Interestingly, overnight charging requirements primarily dominate the contracted grid power capacity. Energy from the grid was drawn predominantly by the depots, as evident from the top panel in Figures \ref{fig:master_variables_durham} and \ref{fig:master_variables_canberra}. All opportunity charging at non-depots was done using BESS. This pattern can be attributed to continuous solar energy generation during the daytime, the lower solar energy costs compared to the ToU grid electricity prices, and the ability to store grid energy in the BESS at a lower cost. 

As discussed earlier, for both networks, as the number of scenarios decreases, the grid capacity gets underestimated, and the solar panel area is overestimated. Therefore, modeling decisions for shorter time intervals, i.e., incorporating more scenarios (52 in our case) to capture diverse climate conditions, is essential for obtaining a realistic picture of costs. We observed no discernible trends in battery storage capacity across the scenarios. However, it can be noted that the depots require higher BESS capacity than other terminals. The maximum grid capacity, panel area, and BESS capacity for Durham are \textcolor{addedcolor}{$319.55$ kW, $1590.78$ \si{m^{2}}, and $65.11$ kWh}, respectively. The corresponding \textcolor{addedcolor}{values} for the Canberra network are \textcolor{addedcolor}{$886.49$ kW, $1962.41$ \si{m^{2}}, and $152.41$ kWh}. These long-term decisions on grid and solar capacity requirements were found to be within practical limits. Additionally, for the 52-scenario case, the maximum energy transfer rates from the grid to BESS are \textcolor{addedcolor}{$5.33$ kWh/min and $14.78$ kWh/min} for the Durham and Canberra networks. These values are also within the limits prescribed in the literature \citep{bess_transfer_1, canadian_solar}.

\subsubsection{\textcolor{addedcolor}{Scenario-Based Schedules}}

For the 52-scenario model, we focus on the results of three distinct scenarios (Weeks 10, 26, and 50) to visualize the impact of incorporating renewables, BESS, and ToU electricity prices. Figure \ref{charge_schedule} illustrates the Gantt chart schedules for these weeks/scenarios. These charts indicate trips in blue, grid-powered charging events (direct or through the BESS) in red, solar-powered charging events through the BESS in green, and grid and solar-powered charging events that happen within the same minute in purple. We find the split of grid/solar charging from the BESS through a postprocessing step by assuming that energy sourced from the grid is spent first. Note that because we use a time discretization of one minute for short-term charging decisions, we do not split charging activities that draw energy from both sources within the one-minute interval. \textcolor{addedcolor}{Schedules consistently prioritize solar charging during daytime.} Additionally, the model optimizes the schedule based on peak-hour pricing, which varies across scenarios. (For Durham, peak and off-peak hours differ in summers and winters.)

\begin{figure}[h]
    \centering\includegraphics[scale=0.275]{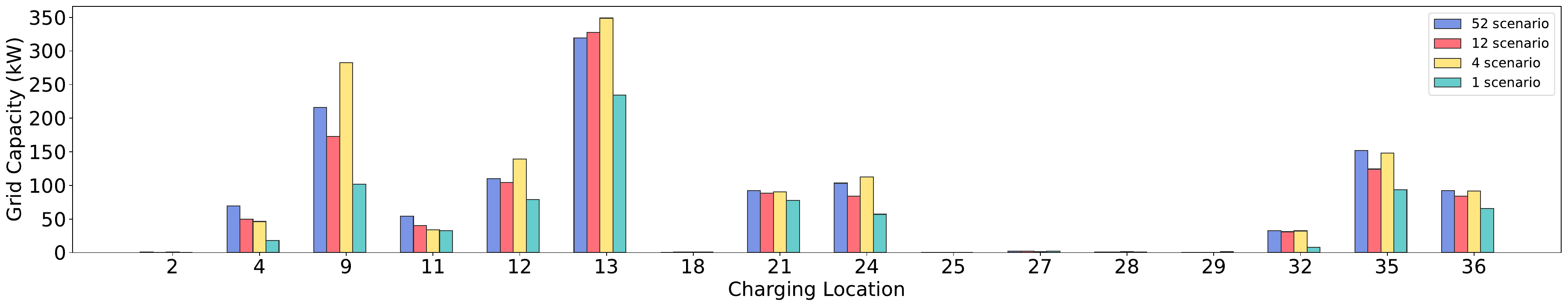}
        \includegraphics[scale=0.272]{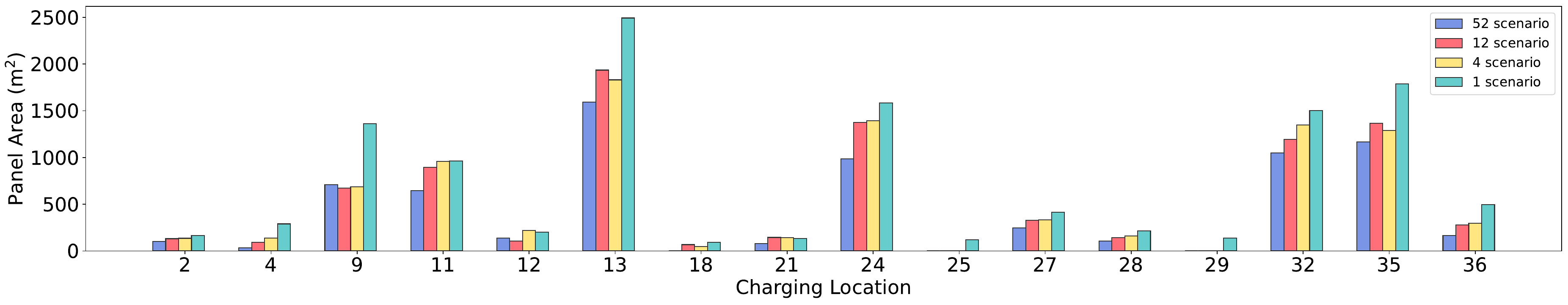}
        \includegraphics[scale=0.275]{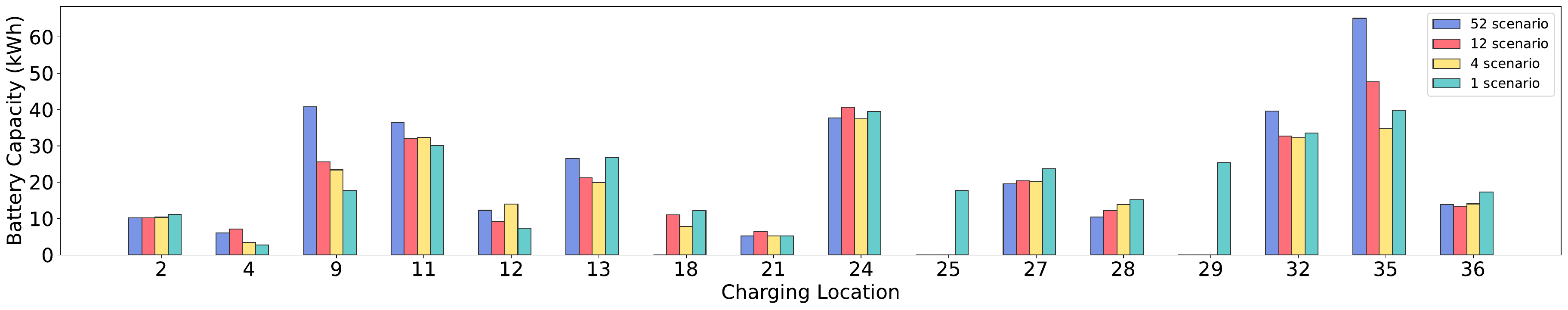}
 \caption{\textcolor{addedcolor}{Contracted grid capacity (top), solar panel area (middle), and BESS capacity (bottom) for the Durham network}}
 \label{fig:master_variables_durham}
\end{figure}    
\begin{figure}[H]
    \centering   
        
        \includegraphics[scale=0.275]{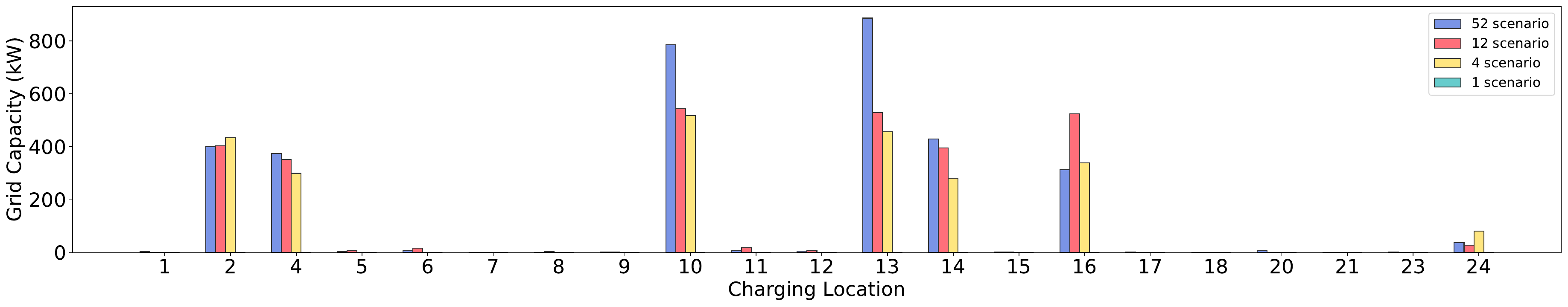}
        \includegraphics[scale=0.272]{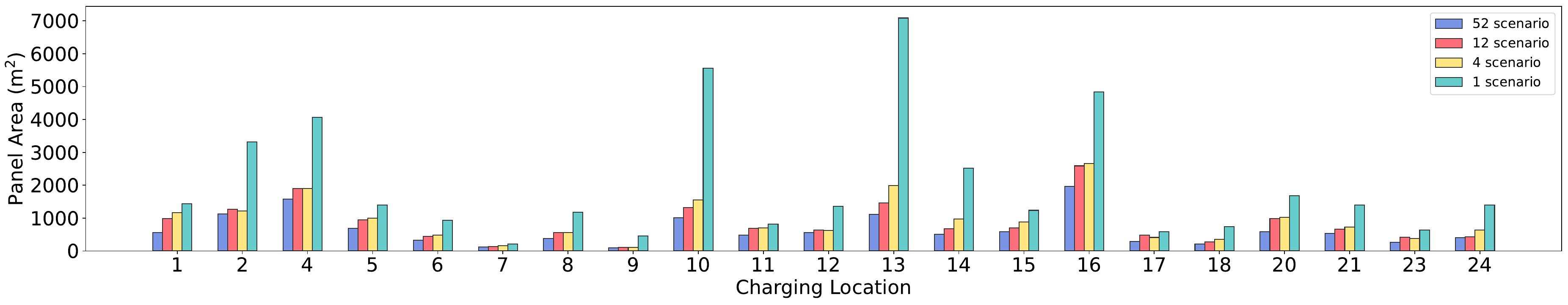}
        \includegraphics[scale=0.275]{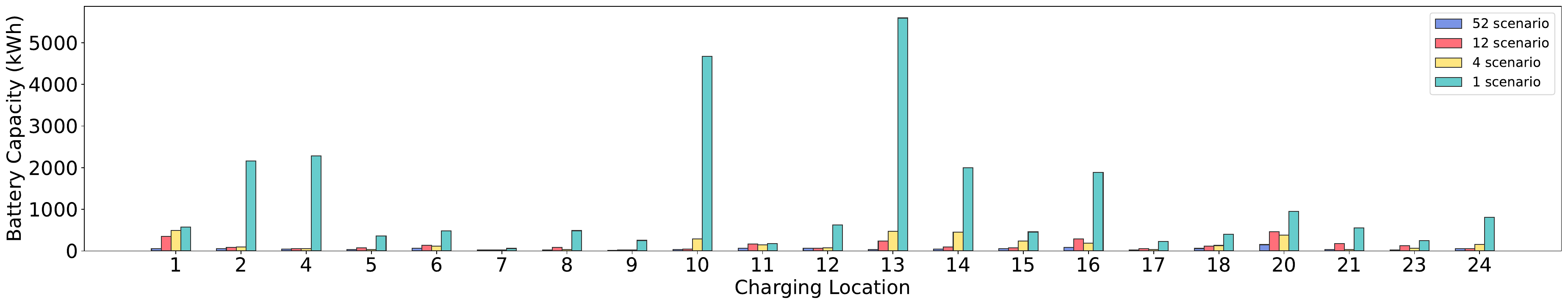}

\caption{\textcolor{addedcolor}{Contracted grid capacity (top), solar panel area (middle), and BESS capacity (bottom) for the Canberra network}} 
\label{fig:master_variables_canberra}
\end{figure}

\begin{figure}[H]
    \centering
     \begin{subfigure}{0.49\textwidth}
        \centering
        \includegraphics[scale=0.07]{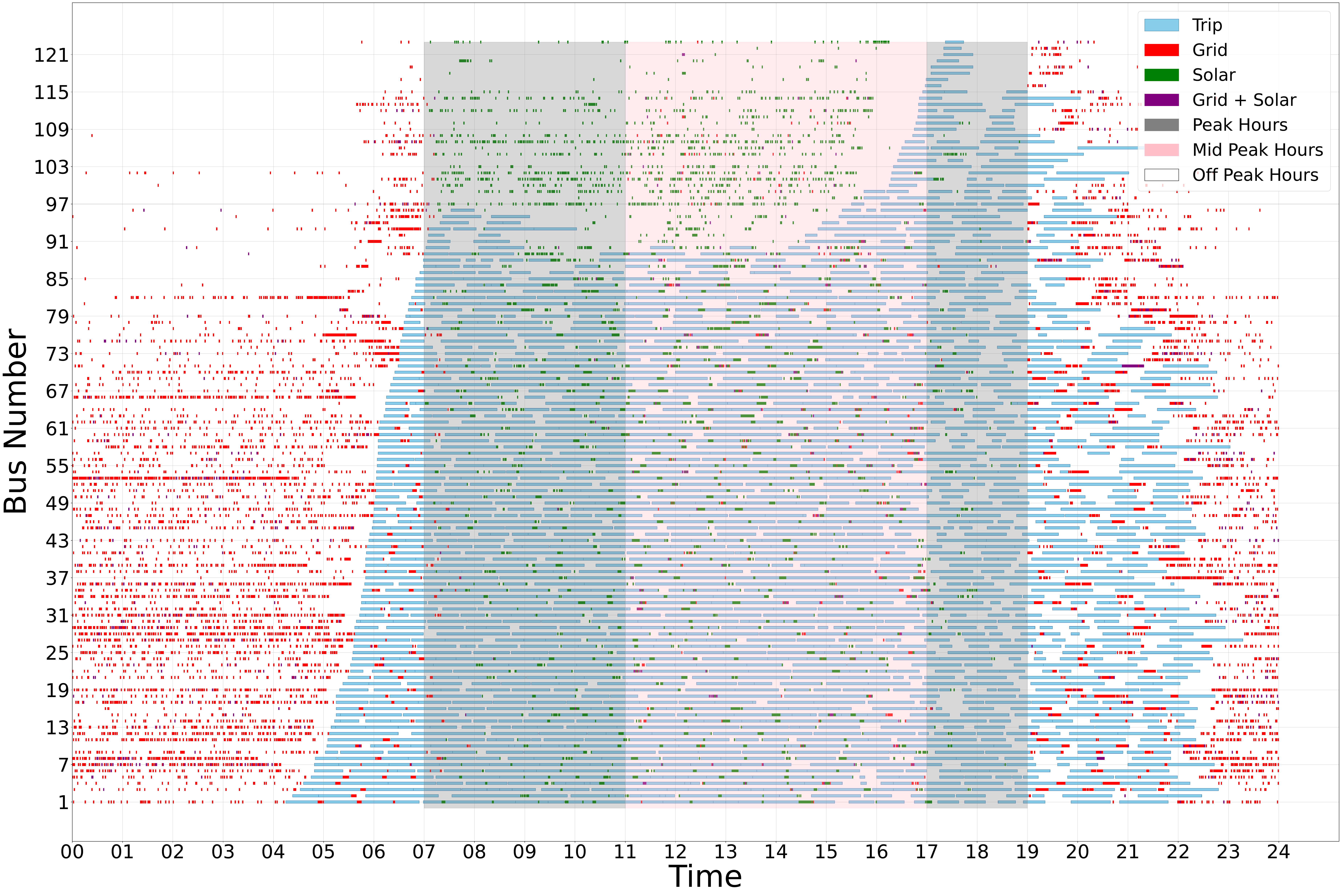}
        \caption{Week 10}
        \label{fig:week10_durham}
    \end{subfigure}
    \hfill
    \begin{subfigure}{0.49\textwidth}
        \centering
        \includegraphics[scale=0.07]{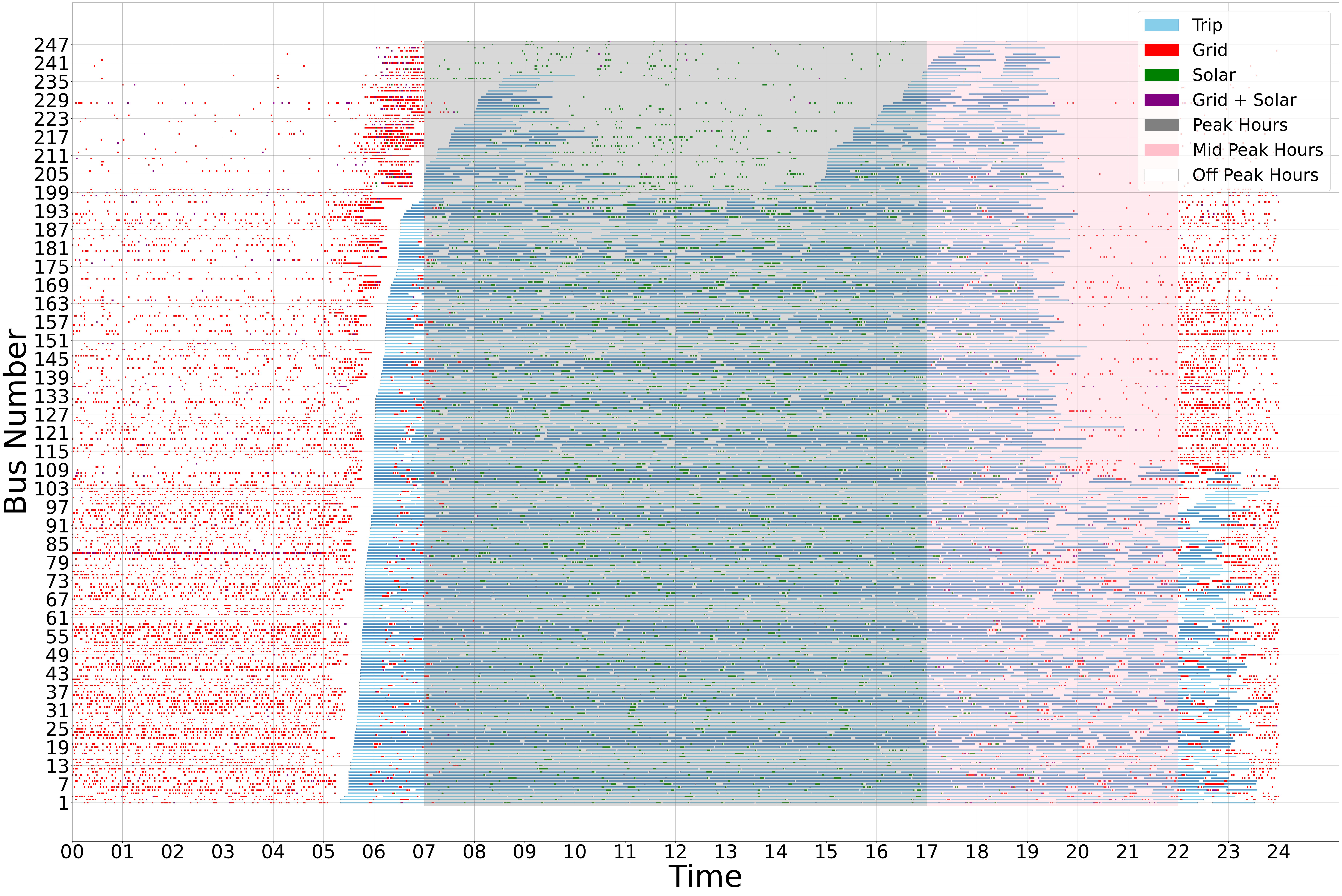}
        \caption{Week 10}
        \label{fig:week10_can}
    \end{subfigure}
    \hfill
        \begin{subfigure}{0.49\textwidth}
        \centering
        \includegraphics[scale=0.07]{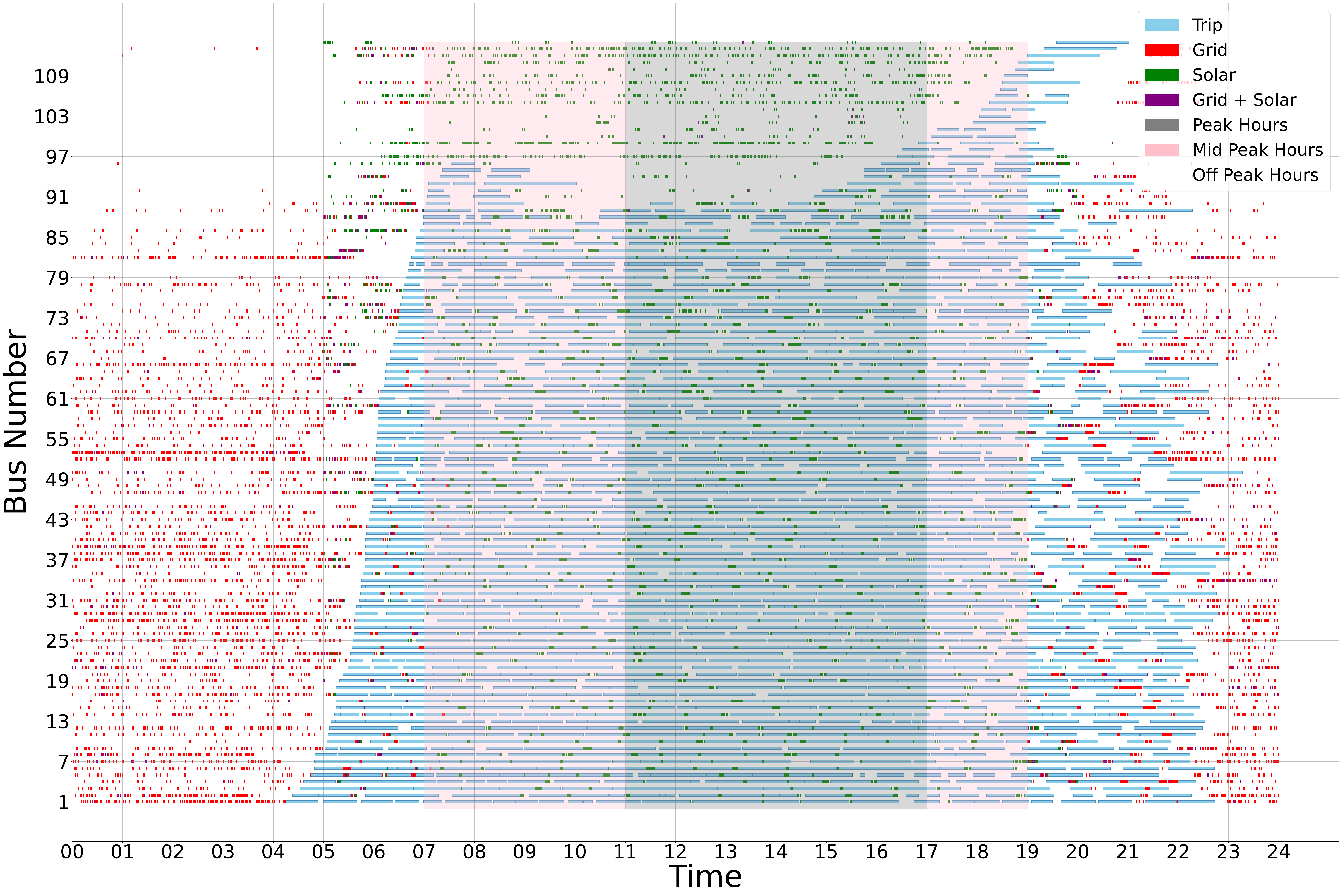}
        \caption{Week 26}
        \label{fig:week26_durham}
    \end{subfigure}
    \hfill
     \begin{subfigure}{0.49\textwidth}
        \centering
        \includegraphics[scale=0.07]{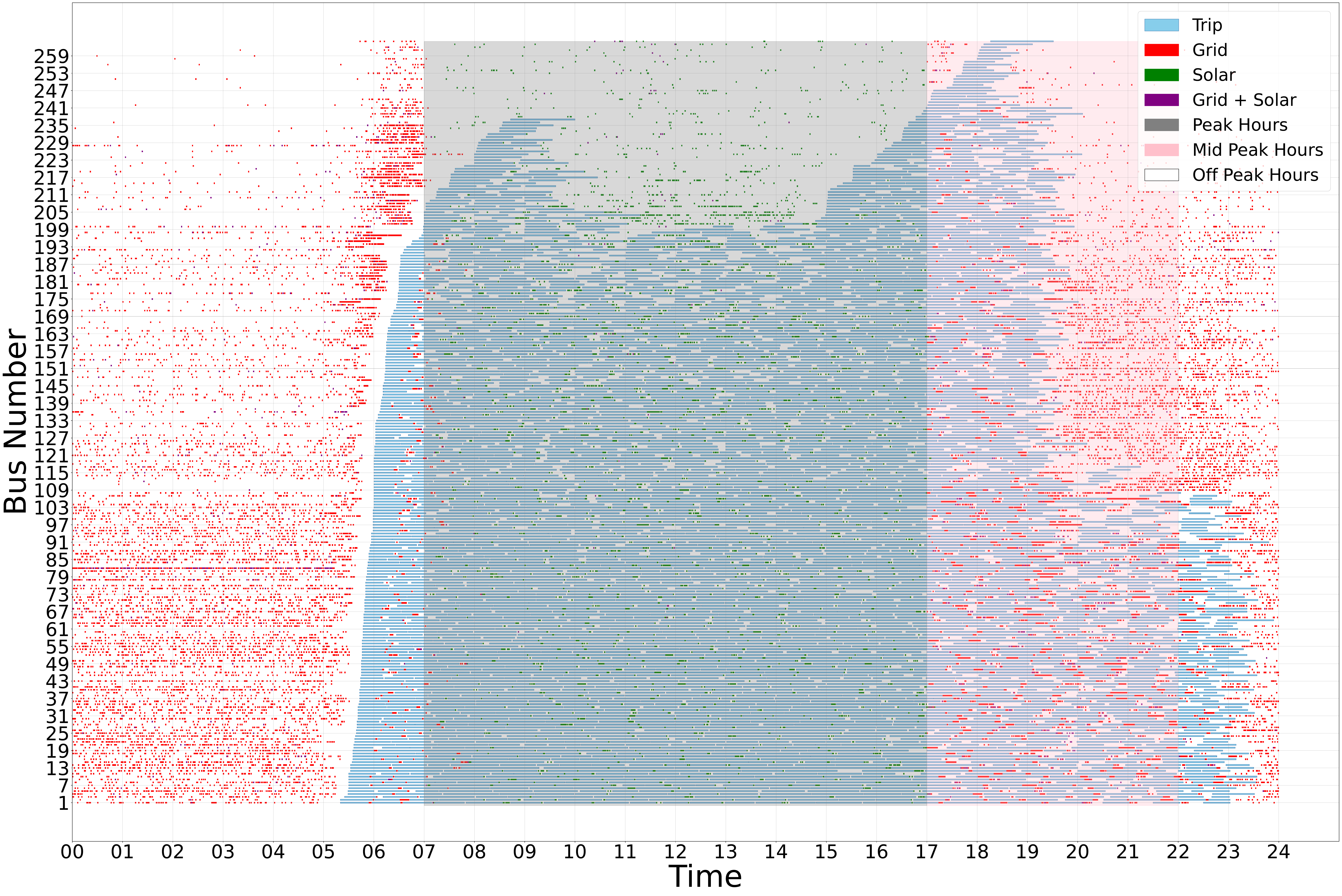}
        \caption{Week 26}
        \label{fig:week26_can}
    \end{subfigure}
    \hfill
    \begin{subfigure}{0.49\textwidth}
        \centering
        \includegraphics[scale=0.07]{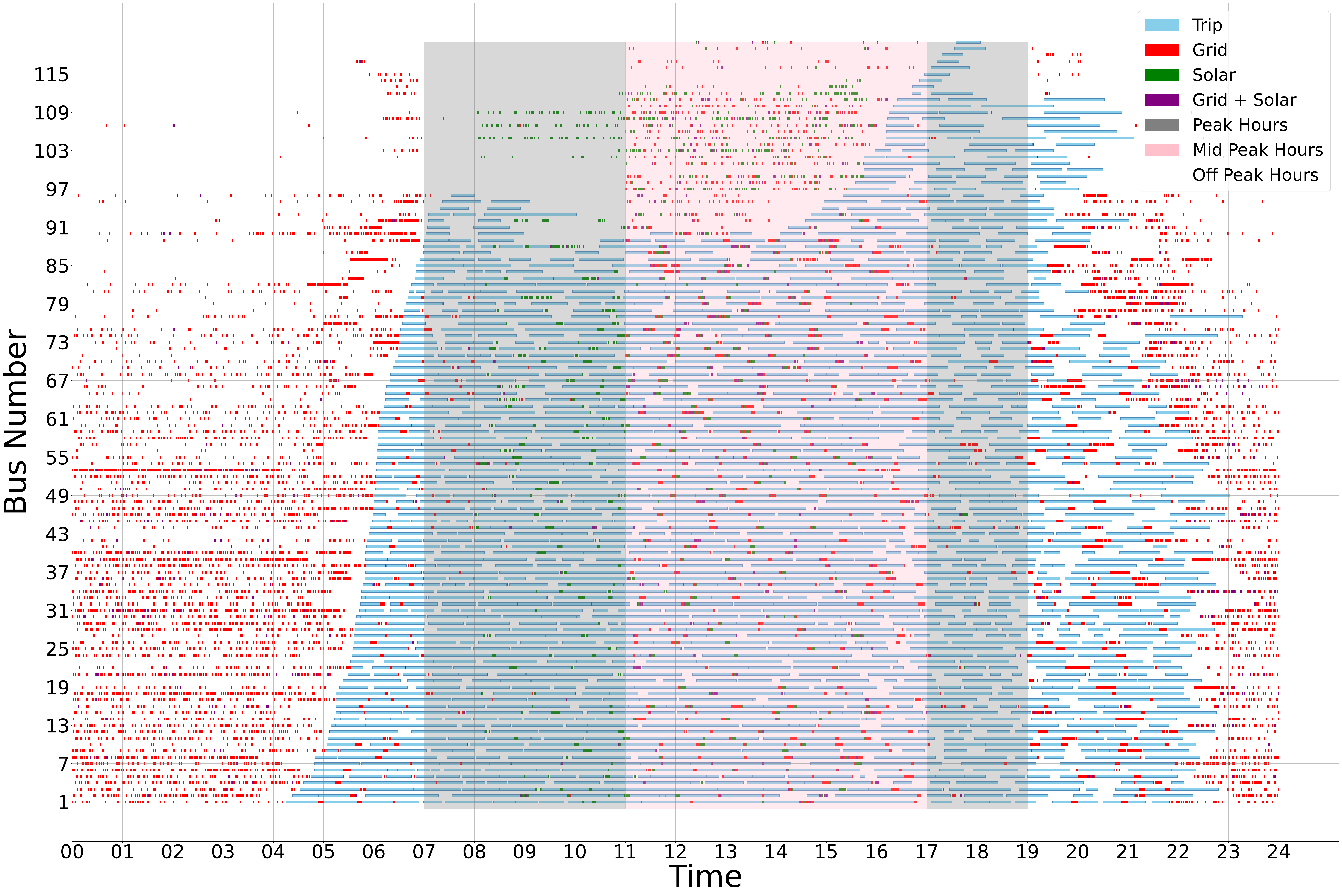}
        \caption{Week 50}
        \label{fig:week50_durham}
    \end{subfigure}
    \hfill
     \begin{subfigure}{0.49\textwidth}
        \centering
        \includegraphics[scale=0.07]{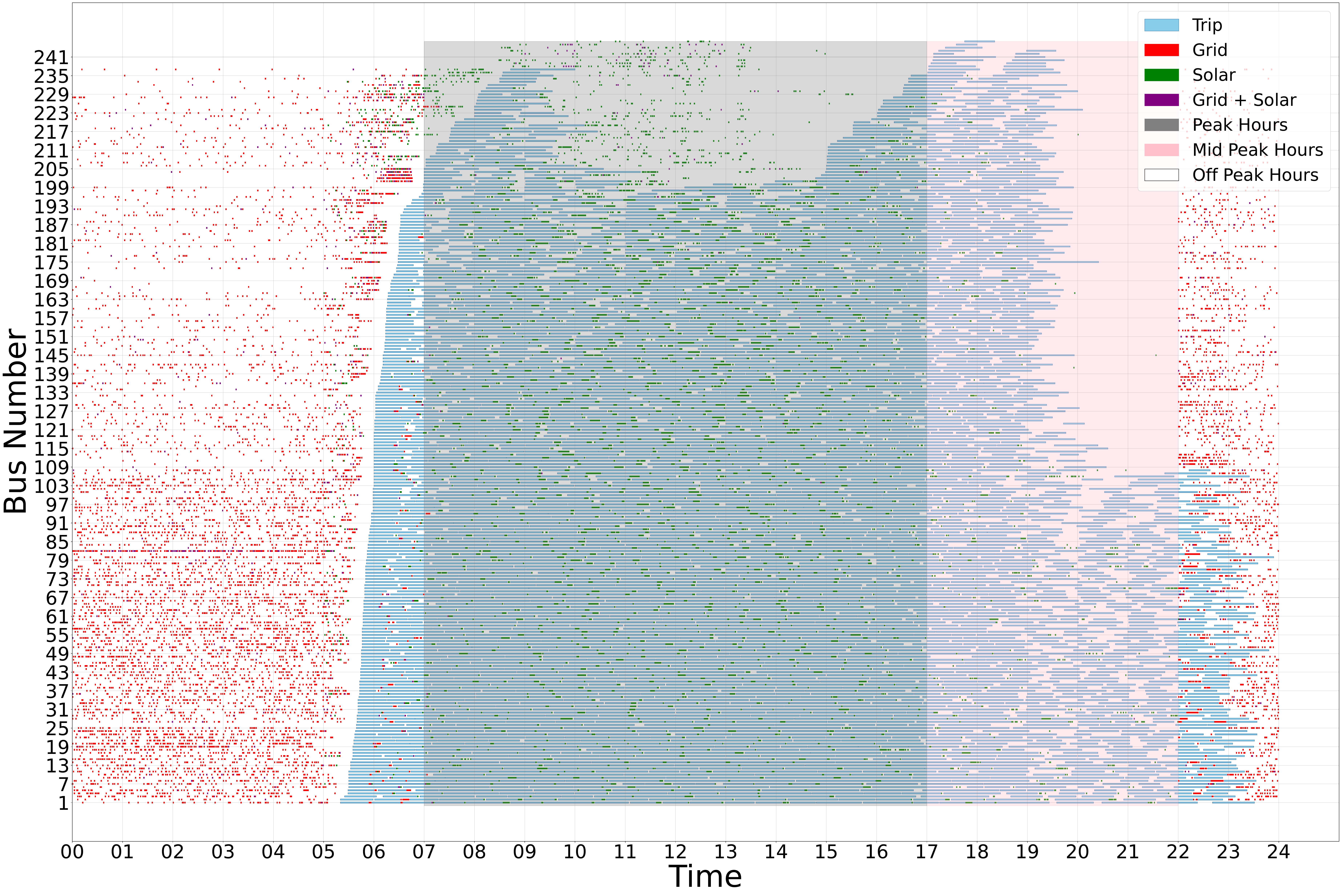}
        \caption{Week 50}
        \label{fig:week50_can}
    \end{subfigure}
    \caption{\textcolor{addedcolor}{Bus operation schedule for weeks  10, 26, and 50 (left: Durham, right: Canberra)}}
        \label{charge_schedule}
\end{figure}

For each network, the number of buses and bus-to-trip assignments \textcolor{addedcolor}{vary} depending on the type of scenario. For Canberra, schedule \ref{fig:week10_can} rarely relies on charging from the grid during peak hours, schedule \ref{fig:week26_can} uses grid charging during off-peak \textcolor{addedcolor}{and mid-peak} hours, and schedule \ref{fig:week50_can} relies entirely on solar energy during peak and \textcolor{addedcolor}{mid-peak hours}. A similar pattern, though less pronounced, can be observed in the Durham schedules due to reduced solar capacity and cheaper ToU pricing. 

\subsubsection{\textcolor{addedcolor}{Scenario-Based Energy Utilization}}

\begin{figure}[H]
    \centering
     \begin{subfigure}{0.49\textwidth}
        \centering
        \includegraphics[scale=0.31]{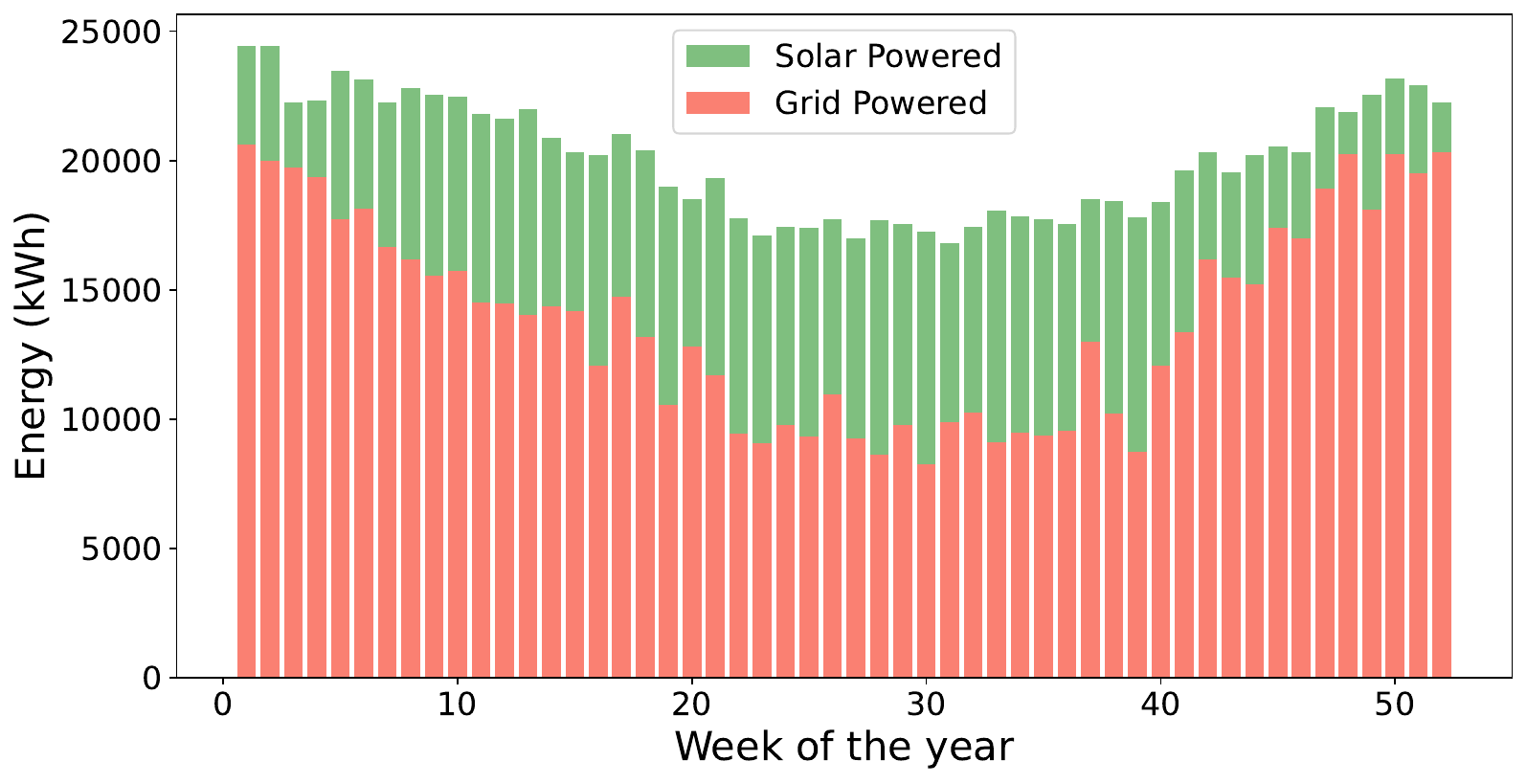}
        \caption{Durham}
        \label{fig:energy_distribution_durham}
    \end{subfigure}
    \hfill
    \begin{subfigure}{0.49\textwidth}
        \centering
        \includegraphics[scale=0.31]{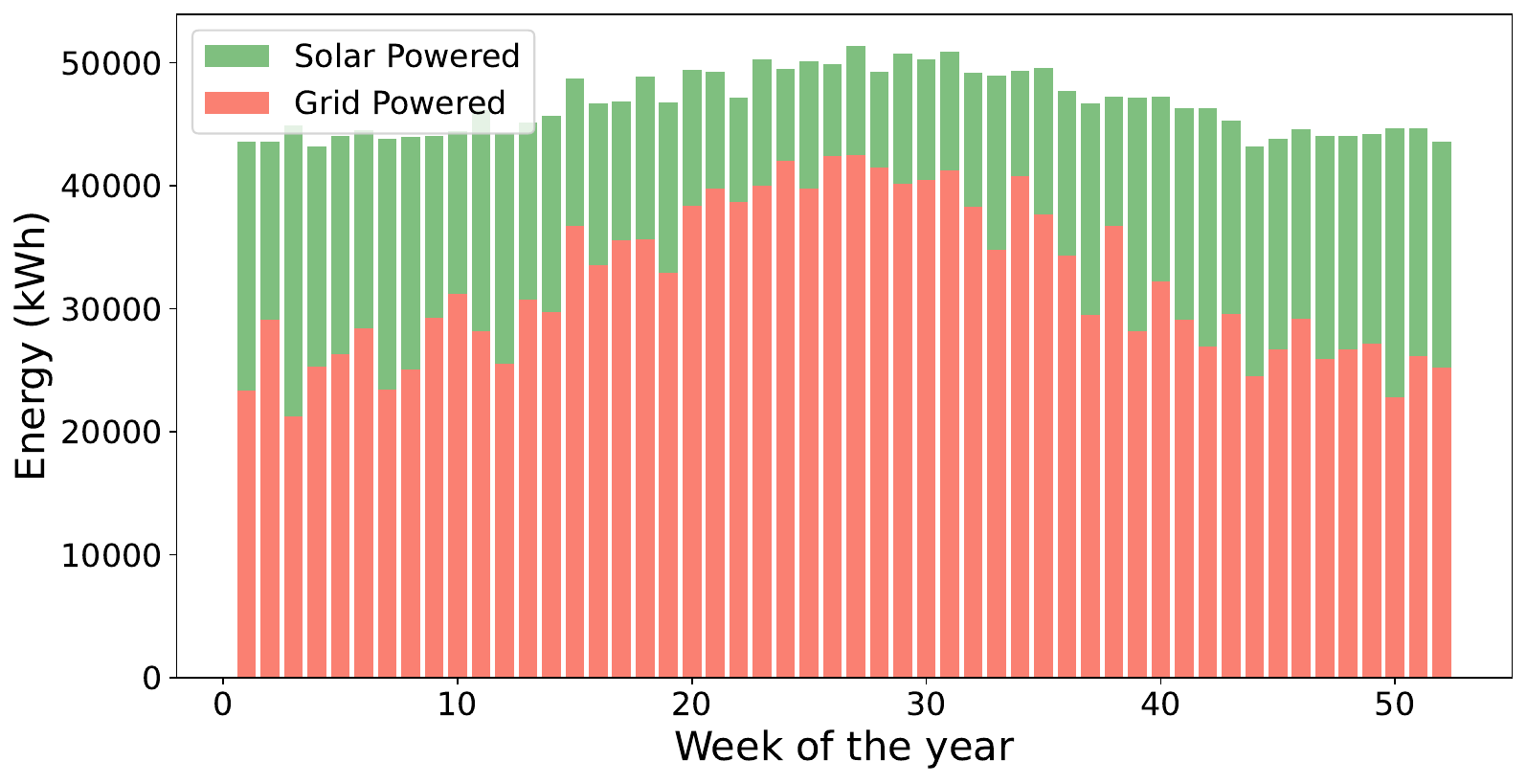}
        \caption{Canberra}
        \label{fig:energy_distribution_can}
    \end{subfigure}
    \caption{\textcolor{addedcolor}{Energy distribution from grid and RES across different scenarios}}
        \label{fig:combined_figures}
\end{figure} 

\begin{figure}[H]
    \centering
    \begin{subfigure}{0.49\textwidth}
        \centering
        \includegraphics[width=0.90\textwidth]{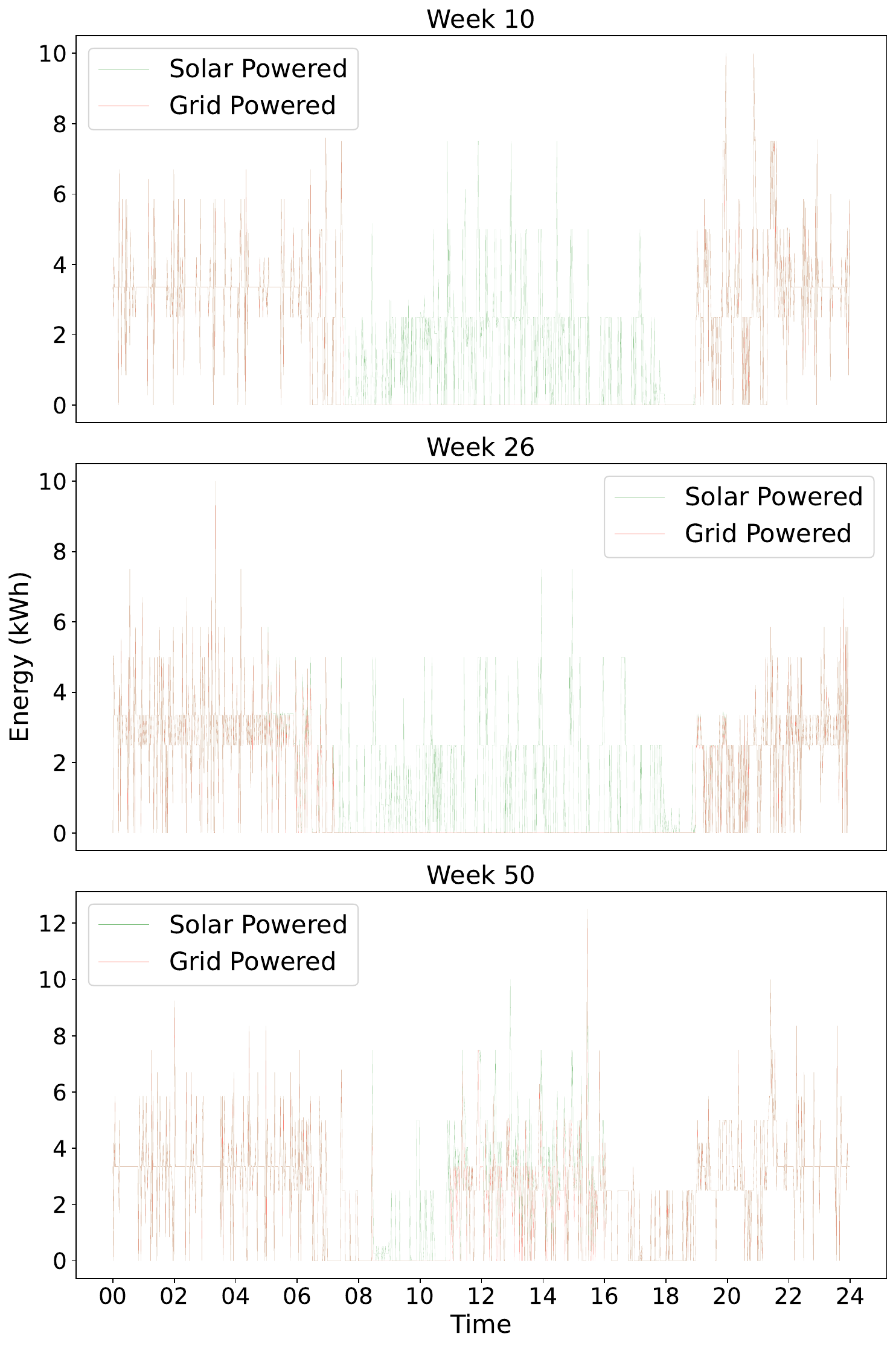}
        \caption{Durham (Location $9$)}
        \label{fig:energy_powered_durham}
    \end{subfigure}
    \hfill
    \begin{subfigure}{0.49\textwidth}
        \centering
        \includegraphics[width=0.90\textwidth]{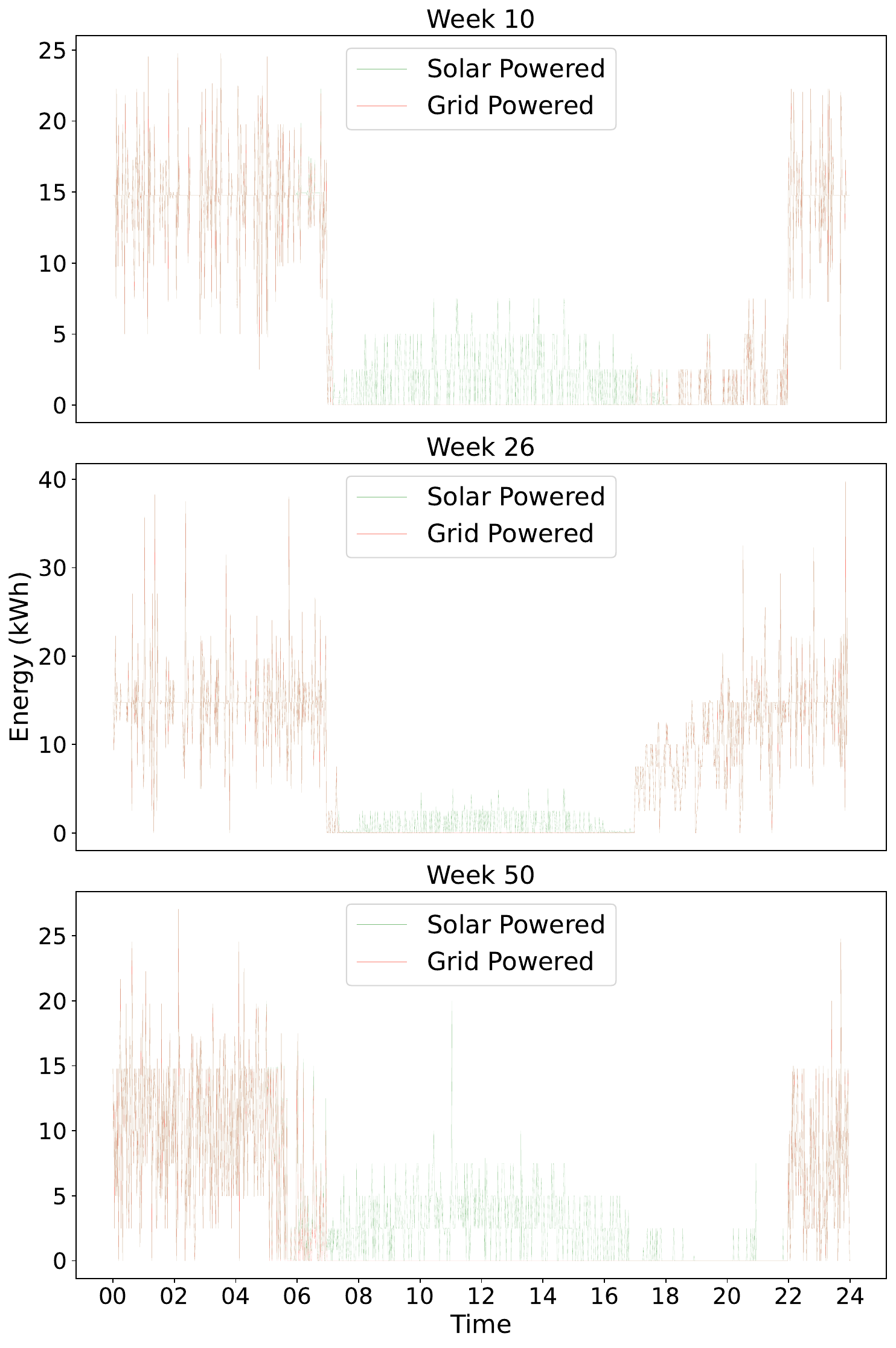}
        \caption{Canberra (Location $13$)}
        \label{fig:energy_powered_canberra}
    \end{subfigure}
    \caption{\textcolor{addedcolor}{Scenario-wise energy utilization for weeks 10, 26, and 50 for Durham (left) and Canberra (right)}}
        \label{fig:scenario_wise_power_utilization}
\end{figure}

Figures \ref{fig:energy_distribution_durham} and \ref{fig:energy_distribution_can} depict the system's total energy consumption across all 52 weeks/scenarios. As previously discussed, the energy requirements for both networks depend on the temperature variations across scenarios. Also, the energy drawn from the grid or BESS fluctuates based on each scenario's GTI and temperature characteristics. A high amount of solar energy is captured during the summer for Durham (mid-year) and Canberra (at the start and end of the year). The contribution of renewables to the total energy requirement was \textcolor{addedcolor}{$30.76\%$ and $31.13\%$} for the Durham and Canberra networks, respectively. Figures \ref{fig:energy_powered_durham}
and \ref{fig:energy_powered_canberra} \textcolor{addedcolor}{illustrate} the energy utilization at a given charging location for both networks across different scenarios. The trend is similar to the overall operations discussed in the Gantt charts. The maximum energy requirement in any time step is different across scenarios. These spikes can occur when the timetables and bus-to-trip assignments are such that multiple buses arrive at the charging station simultaneously. Figures \ref{fig:bess_energy_levels_durham} and \ref{fig:bess_energy_levels_can} show the variation of BESS energy levels across different times-of-the-day in week 50 for the Durham and Canberra networks, respectively. The blue lines indicate no change in the energy levels. Three energy transfer activities can occur every minute: from the solar panel to BESS (dark green), from the grid to BESS (light green), and from BESS to buses (red). Since we use a time discretization of one minute, for visualization purposes, we assume that these events occur sequentially, with each lasting 20 seconds, and show them separately in these figures. We notice that grid and solar power are used to charge the BESS in Durham for the selected location. However, for the Canberra network, solar power is \textcolor{addedcolor}{mainly} used to charge the BESS.

\begin{figure}[H]
    \centering
    \begin{subfigure}{\textwidth}
        \centering
        \includegraphics[scale=0.29]{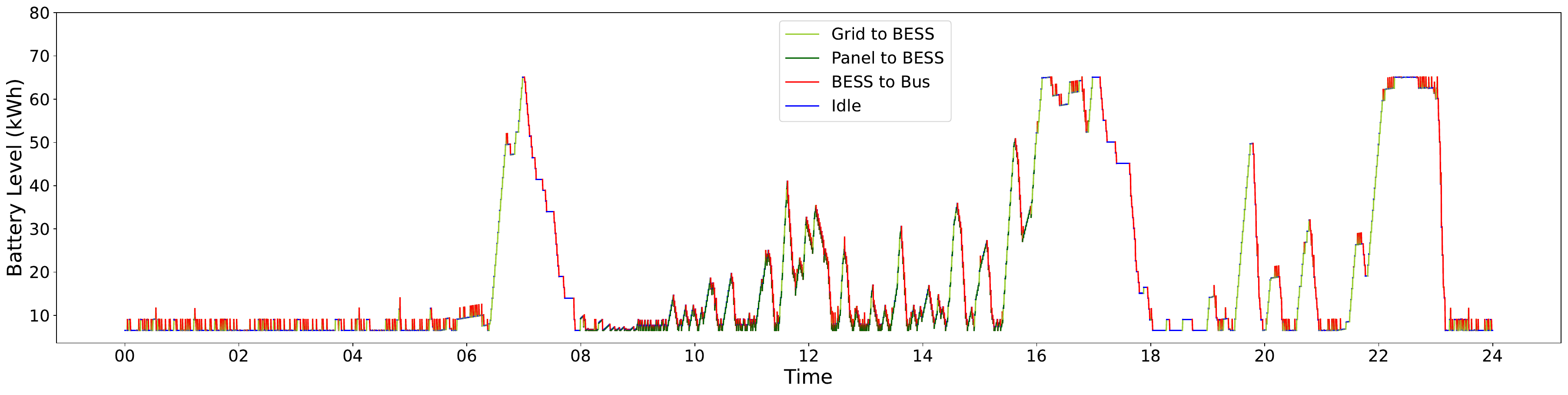}
        \caption{Durham (Location $35$)}
        \label{fig:bess_energy_levels_durham}
    \end{subfigure}
    
    \vspace{0.2cm} 

    \begin{subfigure}{\textwidth}
        \centering
        \includegraphics[scale=0.29]{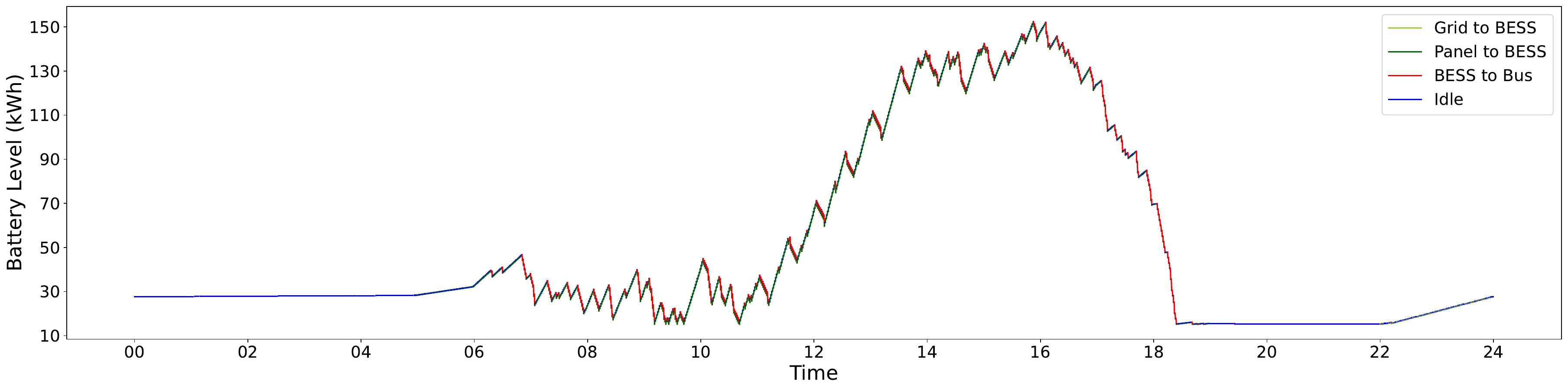}
        \caption{\textcolor{addedcolor}{Canberra (Location $20$)}}
        \label{fig:bess_energy_levels_can}
    \end{subfigure}
    
    \caption{\textcolor{addedcolor}{BESS energy levels for week 50 for Durham (top) and Canberra (bottom)}}
    \label{fig:bess_battery_level}
\end{figure}

\subsubsection{\textcolor{addedcolor}{Results from Local Search-Enhanced Initial Solution}}
\label{sec:local_search}

\textcolor{addedcolor}{The CS algorithm used to generate vehicle schedules is a heuristic. To analyze if the results of our models are sensitive to the bus-to-trip assignments, we generated an alternate vehicle schedule using the local search method proposed in \cite{nath2024impact}. The results obtained using the local search-based EVSP schedules are summarized in Table \ref{tab:local_search_results}, and we compare them with the base-case results reported in Section 5 for the 52-scenario case. For Durham, the total cost changes from \$1406.48 to \$1359.33 (a reduction of approximately 3.35\%), while for Canberra, it changes from \$7190.09 to \$6814.55 (a reduction of approximately 5.22\%). Despite these variations in total cost, the relative contributions of infrastructure and operational costs remain similar in both cases. In particular, the contracted capacity, solar panel, and BESS costs exhibit only minor deviations, indicating that the first-stage infrastructure decisions are largely unaffected due to the bus-to-trip assignments in these two networks. Furthermore, the cost savings from RES integration under the local search-based schedules (10.17\% for Durham and 26.92\% for Canberra) are consistent with those observed in the base case (9.72\% and 23.79\%, respectively), reinforcing the robustness of our key insights.}

\begin{table}[H]
\centering
\caption{\textcolor{addedcolor}{Results using an alternative schedule generated via a local search-based EVSP heuristic}}
\textcolor{addedcolor}{
\begin{tabular}{lrr}
\hline
\textbf{Metric} & \textbf{Durham} & \textbf{Canberra} \\
\hline
Objective (\$) & 1359.33 & 6814.55 \\
Contracted capacity cost (\$) & 217.92 & 558.99 \\
Solar panel cost (\$) & 272.60 & 574.27 \\
BESS cost (\$) & 42.40 & 150.01 \\
Average operation cost (\$) & 826.41 & 5531.28 \\
\hline
\end{tabular}}
\label{tab:local_search_results}
\end{table}

\subsubsection{\textcolor{addedcolor}{Effect of Cluster Radius}}
\label{sec:clustering_sensitivity}

\textcolor{addedcolor}{Our model generates candidate charging locations using a cluster radius of 500 m. The cluster radius defines the maximum distance within which nearby candidate charging locations are aggregated into a single representative site. To evaluate the sensitivity of the daily amortized cost and infrastructure planning decisions to this assumption, Table \ref{tab:cluster_750_results} presents the results obtained by increasing the cluster radius from 500\,m (base case) to 750\,m for the 52-scenario setting. The base-case objective values are \$1406.48 for Durham and \$7190.09 for Canberra. For the 750\,m cluster radius case, the total cost changes marginally to \$1405.96 for Durham and \$7261.96 for Canberra. Despite this change in spatial configuration, the infrastructure cost components (contracted capacity, solar panel area, and BESS capacity) exhibit only minor variations, indicating that the first-stage planning decisions are robust to clustering assumptions. Additionally, the cost savings from RES integration under the 750\,m clustering case (9.73\% for Durham and 23.66\% for Canberra) remain consistent with those observed in the base case.}

\begin{table}[H]
\centering
\caption{\textcolor{addedcolor}{Results for alternative charging location configuration (cluster radius = 750\,m, 52-scenario case)}}
\textcolor{addedcolor}{
\begin{tabular}{lrr}
\toprule
\textbf{Metric} & \textbf{Durham} & \textbf{Canberra} \\
\midrule
Objective (\$) & 1405.96 & 7261.96 \\
Contracted capacity cost (\$) & 230.95 & 612.47 \\
Solar panel cost (\$) & 268.17 & 510.09 \\
BESS cost (\$) & 45.96 & 141.21 \\
Average operation cost (\$) & 860.88 & 5998.19 \\
\bottomrule
\end{tabular}
}
\label{tab:cluster_750_results}
\end{table}

\subsubsection{\textcolor{addedcolor}{Worst-Case Scenario Analysis}}
\label{sec:worst_case_analysis}

\textcolor{addedcolor}{To further investigate the effect of conservative planning, we considered two benchmark cases derived from the weekly scenarios. In the first benchmark, the model was solved using a single scenario corresponding to the week with the highest total trip energy consumption. In the second benchmark, the model was solved using a single scenario corresponding to the week with the lowest solar energy generation. These benchmarks represent planning decisions based on a single extreme operating condition, whereas the proposed scenario-based approach accounts for seasonality across the full year.}

\textcolor{addedcolor}{Table \ref{tab:worst_case_results} presents the worst-case scenario analysis for the 52-scenario setting. The base-case objective values are \$1406.48 for Durham and \$7190.09 for Canberra. For Durham, the total cost increases to \$1755.27 under the highest energy consumption scenario (Week 1), representing an increase of approximately 24.80\%, and to \$1703.87 under the lowest solar generation scenario (Week 52), corresponding to an increase of about 21.14\%. The cost variations for Canberra are less pronounced. The total cost increases to \$7310.91 under the highest energy consumption scenario (Week 27), an increase of approximately 1.68\%, and to \$7411.41 under the lowest solar generation scenario (Week 24), an increase of about 3.08\%. These results highlight the importance of accounting for seasonality during planning, as extreme scenarios can lead to substantially higher costs.}

\begin{table}[H]
\centering
\caption{\textcolor{addedcolor}{Worst-case scenario analysis (52-scenario setting)}}
\textcolor{addedcolor}{
\begin{tabular}{lcc}
\toprule
\textbf{Scenario} & \textbf{Durham} & \textbf{Canberra} \\
\midrule
Base case (52 scenario) & 1406.48 & 7190.09 \\
\midrule
Worst-case (highest energy demand) & 1755.27 (Week 1) & 7310.91 (Week 27) \\
Worst-case (lowest solar generation) & 1703.87 (Week 52) & 7411.41 (Week 24) \\
\bottomrule
\end{tabular}}
\label{tab:worst_case_results}
\end{table}

\subsubsection{\textcolor{addedcolor}{Bus Acquisition Costs}}
\label{sec:bus_acquisition_costs}

\textcolor{addedcolor}{The results presented thus far account for infrastructure and operational costs but exclude bus acquisition costs. Since the required fleet size varies across scenarios, we evaluate how including bus acquisition costs affects the total daily amortized system cost. Assuming a bus acquisition cost of \$381,500 per EB \citep{dirks2022integration}, Table \ref{tab:fleet_size_sensitivity_scenarios} summarizes the number of buses required under each scenario, along with the corresponding total amortized costs if we were to include bus acquisition costs to the previously calculated daily amortized cost.}
\begin{table}[H]
\centering
\caption{\textcolor{addedcolor}{Updated daily amortized costs considering different fleet sizes across scenarios}}
\label{tab:fleet_size_sensitivity_scenarios}

\small
\begin{tabular}{l 
  >{\raggedright\arraybackslash}m{1.5cm} 
  >{\raggedright\arraybackslash}m{1.5cm} 
  >{\raggedright\arraybackslash}m{3cm} 
  >{\raggedright\arraybackslash}m{3.22cm} 
  >{\raggedright\arraybackslash}m{3.22cm}}
\hline
\textcolor{addedcolor}{\textbf{Network}} &
\textcolor{addedcolor}{\textbf{Scenarios}} &
\textcolor{addedcolor}{\textbf{No. of buses}} &
\textcolor{addedcolor}{\textbf{Daily amortized bus costs (\$)}} &
\textcolor{addedcolor}{\textbf{Previous daily amortized cost (\$)}} &
\textcolor{addedcolor}{\textbf{Updated daily amortized cost (\$)}} \\

\hline

\multirow{4}{*}{\textcolor{addedcolor}{Durham}}
& \textcolor{addedcolor}{1}  & \textcolor{addedcolor}{119} & \textcolor{addedcolor}{12,871.28} & \textcolor{addedcolor}{1,179.52} & \textcolor{addedcolor}{14,050.80} \\
& \textcolor{addedcolor}{4}  & \textcolor{addedcolor}{119} & \textcolor{addedcolor}{12,871.28} & \textcolor{addedcolor}{1,363.54} & \textcolor{addedcolor}{14,234.82} \\
& \textcolor{addedcolor}{12} & \textcolor{addedcolor}{122} & \textcolor{addedcolor}{13,195.76} & \textcolor{addedcolor}{1,349.86} & \textcolor{addedcolor}{14,545.62} \\
& \textcolor{addedcolor}{52} & \textcolor{addedcolor}{123} & \textcolor{addedcolor}{13,303.92} & \textcolor{addedcolor}{1,406.48} & \textcolor{addedcolor}{14,710.40} \\
\hline

\multirow{4}{*}{\textcolor{addedcolor}{Canberra}}
& \textcolor{addedcolor}{1}  & \textcolor{addedcolor}{250} & \textcolor{addedcolor}{27,040.50} & \textcolor{addedcolor}{5,219.19} & \textcolor{addedcolor}{32,259.69} \\
& \textcolor{addedcolor}{4}  & \textcolor{addedcolor}{261} & \textcolor{addedcolor}{28,230.28} & \textcolor{addedcolor}{6,585.90} & \textcolor{addedcolor}{34,816.18} \\
& \textcolor{addedcolor}{12} & \textcolor{addedcolor}{262} & \textcolor{addedcolor}{28,338.44} & \textcolor{addedcolor}{6,666.52} & \textcolor{addedcolor}{35,004.96} \\
& \textcolor{addedcolor}{52} & \textcolor{addedcolor}{265} & \textcolor{addedcolor}{28,662.93} & \textcolor{addedcolor}{7,190.09} & \textcolor{addedcolor}{35,853.02} \\
\hline
\end{tabular}
\end{table}

\textcolor{addedcolor}{
As shown in Table \ref{tab:fleet_size_sensitivity_scenarios}, the number of buses required increases as the number of scenarios increases. For Durham, the required fleet size increases from 119 buses under the 1- and 4-scenario cases to 123 buses for the 52-scenario case. Similarly, for Canberra, the required fleet size increases from 250 buses (one scenario) to 265 buses (52 scenarios). Consequently, the total daily amortized bus acquisition costs and the overall daily amortized system costs also increase with the number of scenarios considered.}

\textcolor{addedcolor}{
These results reinforce our observation that using more scenarios better captures seasonality. In contrast, considering only a small number of scenarios may underestimate the fleet size requirements, operational costs, and, consequently, the total system cost. The 52-scenario case, therefore, provides a realistic estimate of the total system cost by accounting for variability throughout the year.
}

\section{Conclusions and Future Work}
\label{sec:conc}

This study introduced a two-stage multi-scenario \textcolor{addedcolor}{LP} model to determine the charging schedules of EBs using energy from both a conventional grid and a solar-powered system. We addressed the seasonality in solar energy generation and variations in trip energy consumption due to changes in ambient temperature within a day. The first stage, or planning-level decisions, involved determining the optimal power capacity for charging stations, the area of solar panels to be installed, and the capacity of the BESS at the charging stations. \textcolor{addedcolor}{The second-stage decisions are scenario-specific and determine time-dependent energy flows between the grid, BESS, and buses.} The BESS can store energy from the solar panels or the grid. A concurrent scheduler-based algorithm was used to obtain scenario-specific bus-to-trip \textcolor{addedcolor}{assignments}, which were used as an input to the CSP. A Charge-And-Go priority scheme was assumed, where buses preferred charging at the current trip's end-stop rather than \textcolor{addedcolor}{at} the next trip's start-stop. Charging rates were assumed to be piecewise constant across time steps, and partial charging was permitted. \textcolor{addedcolor}{To ensure computational scalability, we adopt a formulation that does not explicitly enforce charging continuity, discrete outlet allocation, or mutually exclusive BESS charging/discharging, and these are interesting extensions that could be addressed in the future.} Case studies conducted on two real-world bus transit networks, Durham and Canberra, considering RES and seasonal ToU electricity prices, demonstrated significant cost savings of \textcolor{addedcolor}{$9.72\%$ and $23.79\%$}, respectively. Due to the block diagonal structure of the problem, we solved the CSP using Benders' decomposition and compared its performance with that of the dual simplex method.

Our scenario-based LP model effectively captured the seasonal variations in RES, offering a more accurate representation of real-world conditions. \textcolor{addedcolor}{An additional extension could be to jointly optimize the battery capacity of the fleet, thereby capturing trade-offs between battery sizing, charging infrastructure, and operational costs across scenarios. Since vehicle Original Equipment Manufacturers (OEMs) typically offer only a limited number of battery-capacity configurations, one could also solve the proposed model for each available configuration and select the one that minimizes the overall system cost.} As more EB and BESS data become available, future research could explore RES integration with non-linear charging profiles, bus battery drive cycles, BESS calendar and cycle aging, and online decision-making with delays in trip schedules. \textcolor{addedcolor}{Future work can also incorporate integer decision variables to explicitly model these operational constraints, albeit at the cost of increased computational complexity.} \textcolor{addedcolor}{Other promising future directions include (1) extending the model to a multi-day framework with explicit representation of overnight operations, (2) incorporating bidirectional charging, and (3) joint optimization of vehicle assignment and charging decisions.} Nonetheless, our study provides strong support for integrating renewable energy and adapting to power generation fluctuations in the charge scheduling of electric bus fleets. 


\section*{Acknowledgments}
The authors thank the Ministry of Education for supporting this study through the Scheme for Transformational and Advanced Research in Sciences (STARS2/2023-0252) project titled \textit{Efficient Algorithms for Large-Scale Public Transit Planning}.

\textbf{{\large CRediT authorship contribution statement}}

\textbf{Rito Brata Nath, Madhusudan Baldua} Methodology, Software, Writing -- Original Draft, Data Curation, Investigation, Visualization; \textbf{Vivek Vasudeva:} Software, Investigation, Writing - Review \& Editing; \textbf{Tarun Rambha:} Conceptualization, Methodology, Investigation, Writing - Review \& Editing, Supervision. 

\appendix
\renewcommand{\thesection}{Appendix \Alph{section}}

\section{Glossary}
\label{sec:abbreviations}

\begin{table}[H]
 \small
  \centering
  \caption{List of abbreviations}
\begin{tabular}{p{3.2cm}p{6.0cm}} 
\hline
       \textbf{{Acronym}} & \textbf{{Description}}  \\ 
  \hline
    EBs & Electric Buses\\
    EVs & Electric Vehicles\\
    EVSP & Electric Vehicle Scheduling Problem\\
    CSP & Charge Scheduling Problem\\
    MDVSP & Multi-Depot Vehicle Scheduling Problem\\
    ToU & Time-of-Use\\
    RES & Renewable Energy Sources\\
    PVs & Photo-Voltaics\\
    GTI & Global Tilted Irradiance\\
    BESS & Battery Energy Storage System\\
    LP & Linear Programming\\
    CS & Concurrent Scheduler\\
    NREL & National Renewable Energy Laboratory\\
    MILPs & Mixed-Integer Linear Programs\\
    MDP & Markov Decision Process \\
    CAG & Charge-and-Go\\
    GTFS & General Transit Feed Specification\\
    DoD & Depth-of-Discharge\\
    SoC & State-of-Charge\\
    \hline
\end{tabular}
\end{table}

\bibliography{references}
\appendix

\end{document}